\documentclass[10pt]{article}
\usepackage{flafter,amsmath,amssymb,latexsym,psfrag,graphicx,color,indentfirst}
\usepackage[margin=2.5cm]{geometry}

\usepackage[colorlinks,
linktocpage=true,
linkcolor=black,
citecolor=black,
urlcolor=black,
anchorcolor=black
]{hyperref}

\usepackage[numbers,sort&compress]{natbib}
\usepackage{appendix}
\usepackage{multirow}
\usepackage{caption}
\usepackage{subfigure}
\usepackage{graphicx}
\usepackage{mathtools}
\usepackage{esint}
\usepackage{ulem}

\DeclareGraphicsExtensions{.eps,.pdf,.jpg,.png}
\allowdisplaybreaks

\newtheorem{theorem}{Theorem}[section]

\newtheorem{lemma}[theorem]{Lemma}

\numberwithin{equation}{section}
\numberwithin{figure}{section}

\begin{document}
\setlength\arraycolsep{2pt}
\date{\today}
\title{Minnaert Frequency and Simultaneous Reconstruction of the Density, Bulk and Source in the Time-Domain Wave Equation}
\author{Soumen Senapati and Mourad Sini 
\\
\\
RICAM, Austrian Academy of Sciences, A-4040, Linz, Austria
\\
E-mail: soumen.senapati@ricam.oeaw.ac.at and mourad.sini@oeaw.ac.at}

\maketitle
\begin{abstract}
We deal with the inverse problem of reconstructing acoustic material properties or/and external sources for the time-domain acoustic wave model. The traditional measurements consist of repeated active (or passive) interrogations, as the Dirichlet-Neumann map, or point sources with source points varying outside of the domain of interest. It is reported in the existing literature, that based on such measurements, one can recover some (but not all) of the three parameters: mass density, bulk modulus or the external source term. In this work, we first inject isolated small-scales bubbles into the region of interest and then measure the generated pressure field at a {\it{single point}} outside, or at the boundary, of this region. Then we repeat such measurements by moving the bubble to scan the region of interest. Using such measurements, we show that
\begin{enumerate}
\item If either the mass density or the bulk modulus is known then we can simultaneously reconstruct the other one and the source term. 

\item If the source term is known at the initial time, precisely we assume to know its first non vanishing time-derivative, at the initial time, then we reconstruct simultaneously the two parameters, namely the mass density with the bulk modulus and eventually the source function.
\end{enumerate}
Here, the source term, which is space-time dependent,  can be active (and hence known) or passive (and unknown). It is worth mentioning that in the induced inverse problem, we use measurements with $4=3+1$ dimensions ($3$ in space and $1$ in time) to recover $2$ coefficients of $3$ spatial dimensions, i.e. the mass density and the bulk modulus and the
4 = 3 + 1 dimensional source function. In addition, the result is local, meaning that we do reconstruction in any subpart, of the domain of interest, we want.


\bigskip

\end{abstract}


\newpage

\section{Introduction and statement of the results}\label{sec_intro}
\subsection{Introduction}\label{introduction}

We consider a bounded domain $\Omega \subset \mathbb{R}^3$ with a smooth boundary $\partial\Omega$ and $T>0$. Let us denote the bulk modulus and density of a medium which occupies the region $\Omega$ by smoothly varying positive functions $k_0(\cdot)$ and $\rho_0(\cdot)$ respectively, satisfying $k_0(x) = \rho_0(x) = 1$ in $\mathbb{R}^3 \setminus \overline\Omega$ and $ c^{-1}\le k_0(x),\, \rho_0(x) \le c,$ in $\Omega$ for some  constant $c>0$. The acoustic wave speed is then given by 
\begin{align}\label{speed}
    c_0(x) = \sqrt{\frac{k_0(x)}{\rho_0(x)}}, \quad x \in \mathbb{R}^3
\end{align}
which, in essence, we assume to be identically one outside the medium $\Omega$. Next, we consider the following initial value problem, IVP in short, in the presence of the source term $J(\cdot, \cdot)$, 
\begin{align}\label{IVP_v}
    \begin{cases}
    k_0^{-1}(x) v_{tt} - \textnormal{div}\left(\rho_0^{-1}(x) \nabla v\right) = J(x,t), \,  \textnormal{ in } \mathbb{R}^3 \times \mathbb{R}_{+}, \\
    v(\cdot,0) = v_t(\cdot,0) = 0, \, \textnormal{ in } \mathbb{R}^3.
   \end{cases}
\end{align}

\medskip

In accordance with the IVP \eqref{IVP_v}, one typically formulates an inverse problem which amounts to recover the functions $k_0(\cdot), \, \rho_0(\cdot)$ and the source term $J(\cdot,\cdot)$ from the medium response observed on the lateral boundary i.e. $v(\cdot,\cdot)\vert_{\partial\Omega \times (0,T)}$. This inverse problem is of importance in many applications using sonar interrogations. Note here that the IVP is formulated in the whole space rather than in a bounded domain setting. In the later case, one typically prescribe the Neumann derivative of $v$ on the lateral boundary. The literature pertaining to the recovery of acoustic properties or the source function in time-domain wave equation is quite rich and the subject garnered  attention and extensive mathematical treatment during last several decades from many authors. The use of Carleman estimates in this direction was introduced by Bukhge\u{\i}m and Klibanov in their seminal work \cite{Bukhgeim_Klibanov_1981} to address the uniqueness question for a multi-dimensional coefficient identification problem from a single measurement. The method can be also adapted to the case when one tries to determine the source function (in a bounded domain formulation) from the Neumann observation of $v$ on the lateral boundary. This approach led to subsequent development and improvement by Bellassoued, Imanuvilov, Yamamoto et al. in their works \cite{Bellassoued_2004, Bellassoued_Yamamoto_2006_JMPA, Immanuvilov_Yamamoto_2001_CPDE, Immanuvilov_Yamamoto_2001_IP, Yamamoto_1999_JMPA}.  In this approach, the medium properties $k_0(\cdot)$ and $\rho_0(\cdot)$ are always known a-priori and the goal is to recover $J(\cdot,\cdot)$. This method also presumes that the source function is incompletely separated into its spatial and temporal components with the later being completely known in $Q$ and belongs to certain admissible class of functions. To be precise, one considers $J(x,t) = f(x)\, h(x,t)$ where $f(\cdot)$ is the only unknown to be recovered from the boundary response and $h(\cdot,0)$ admits a strict sign condition throughout $\Omega$.  In connection to the coefficient identification problem, the later condition translates to non-vanishing assumption of initial velocity which, for obvious reason, is not suitable in the physical context. Nonetheless, Rakesh and Salo in \cite{Rakesh_Salo_1, Rakesh_Salo_2} avoided such assumptions by working with more realistic sources such as plane wave or point sources. Therefore, in terms of our model \eqref{IVP_v}, these results can be rephrased to recover one of the coefficients $\rho_0(\cdot)$, $k_0(\cdot)$ or $J(\cdot,\cdot)$. Taking the parameter count into account, we note that these are formally-determined inverse problems meaning that the measured data and unknown coefficients share the same number of variables. For a detailed discussion on this approach, we refer to the monographs \cite{Bellassoued_Yamamoto_book, Klibanov_Timonov_book}.  Based on this approach, the authors in \cite{Klibanov_Romanov_2023} recently derived a H\"older stability result for the determination of wave-speed from the data coming from a delta type plane wave source.

 In the context of Riemannian wave equation, let us mention that there is yet another powerful tool, known as the Boundary Control (BC) method, pioneered by Belishev \cite{Belishev_1987, Belishev_1997} which is subsequently adapted and developed by Katchalov, Kurylev, Lassas, see \cite{K-K-L-2001}, and Oksanen et. al. in numerous different settings \cite{Belishev_Kurylev_1992, Kian_Kurylev_Lassas_Oksanen_2019, Lassas_Oksanen_2010, Lassas_Oksanen_2014}. The method is reconstructive and aims at determining the metric structure or lower order coefficients from the hyperbolic DN map which records the Neumann output produced by the medium when exposed to  infinitely many Dirichlet sources on the lateral boundary. The approach, taking advantage of its hyperbolic feature, draws important connection between controllability of waves and the theory of multi-dimensional inverse problems. We refer to \cite{Belishev_2007} and \cite{K-K-L-2001} for an overview.

\smallskip

In this work, we focus on developing a method to reconstruct the medium properties i.e. $k_0(\cdot), \, \rho_0(\cdot)$ and internal source function $J(\cdot,\cdot)$ simultaneously from some measurements made on the lateral boundary i.e. $\partial\Omega \times (0,T)$. To define these measured data, let us first introduce an IVP which models the wave-field generated by the medium when we have injected a micro-scaled contrast agent (bubbles in this case) within the medium at a point location $z \in  \Omega$ and then employ the same source function as before i.e. $J(\cdot,\cdot)$. Let us emphasize on the point that the contrast agent has low bulk modulus and density enjoying a critical scale with respect to its size. To illustrate these ideas, let us consider a bounded $C^2$-smooth domain $B$, containing the origin of coordinates and having a volume $\vert B\vert$ of order $1$, and a small parameter $0 <\epsilon \ll 1$ (as compared to $\vert B\vert$). We denote by $D(z) : = z + \epsilon B$ the region where we inject the contrast agent. Therefore, the parameter $\epsilon$ refers to the smallness aspect of the volume which the contrast agent possess. As a result of injecting these contrast agents, we presume that the medium properties undergo the following alterations:
\begin{align}\label{scales of bulk modulous and density}
    k(x) = 
\begin{cases}
 k_1,\quad \textnormal{ in } D(z),\\
k_0(x) ,\quad \textnormal{ in } \mathbb{R}^3\setminus D(z),
\end{cases}
 \quad \quad \rho(x) = 
\begin{cases}
 \rho_1 ,\quad \textnormal{ in } D(z),\\
\rho_0(x) ,\quad \textnormal{ in } \mathbb{R}^3\setminus D(z),
\end{cases}
\end{align}
where $k_1, \, \rho_1$ are positive constants satisfying
\begin{equation}\label{scales-of-contrast}
k_1 \simeq \epsilon^2 ~\mbox{ and }~ \rho_1\simeq\epsilon^2.
\end{equation}
Such scales are natural for different gases enjoying contrasting properties with liquids (or other backgrounds) of the form $\frac{k_1}{k_0}$ and $\frac{\rho_1}{\rho_0}$ very small. Then from these small scales, we design the radius $\epsilon$ for which we have (\ref{scales-of-contrast}).
Due to this change, the wave speed inside $D(z)$ becomes $c_1=\sqrt{\frac{k_1}{\rho_1}}$. Therefore, the resulting wave-field satisfies the IVP 
\begin{align}\label{IVP_u}
    \begin{cases}
    k^{-1}(x) u_{tt} - \textnormal{div}\left(\rho^{-1}(x) \nabla u\right) = J(x,t),  \,  \textnormal{ in } \mathbb{R}^3 \times \mathbb{R}_{+},\\  
    u\vert_{+} = u\vert_{-} , \quad \rho_1^{-1} \partial_\nu u\vert_{+} = \rho^{-1}_0 \partial_\nu u\vert_{-} \quad \textnormal{ on } \partial D(z)\times\mathbb{R}_{+}, \\
    u(\cdot,0) = u_t(\cdot,0) = 0, \, \textnormal{ in } \mathbb{R}^3,
   \end{cases}
\end{align}
where $\nu(x)$ and $\partial_{\nu(x)}$ denote the normal vector and the normal derivative respectively at $x\in \partial D(z)$. With $f\vert_{+}(x)$ and $f\vert_{+}(x)$, we indicate the limiting value of some function $f$ along the normal vector $\nu(x)$ from inside and outside respectively. The appearance of transmission condition in the IVP \eqref{IVP_u} is quite natural as the density function $\rho(\cdot)$ experiences jump on the interface $\partial D(z)$.  Although, the modified wave speed inside the bubble $D$, i.e. $c_1$, is moderate when compared to the wave speed of the background, there is a large jump in the transmission condition across $\partial D(z)$. These properties are pivotal for our analysis in the current work. In terms of notations, one should introduce an extra argument in the wave-field which is $u(x,t;z)$, to underline its dependence with respect to the position of injected bubble. As our analysis largely depends on deriving an asymptotic profile of $u(x,t;z)$ with regard to the scale $\epsilon$ with $z$ being fixed, we merely denote the latter wave-field by $u(x,t)$. Similarly, we use the notation $D$ instead of $D(z)$.  

Now, we are in a position to elaborate on the measured data to be considered in our setting. At a fixed point $x_0 \in \partial\Omega$, we record the medium response, i.e. the pressure, before and after injecting the contrast agent at different positions. That is, we allow $z$ to vary within $\Omega$. Hence, the data consists of the functions 
\begin{align*}
   v(x_0,t)\big\vert_{t\in(0,T)} \, \textnormal{ and, } \, u(x_0,t;z)\big\vert_{t\in(0,T), \, z\in\Omega}.
\end{align*}
From these measurements, we aim at reconstructing the material properties $k_0(\cdot)$ and $\rho_0(\cdot)$ and the source function $J(\cdot, \cdot)$. Similar to many standard coefficient identification problems pertaining to linear PDEs, we also notice that this inverse problem is also non-linear, as we aim to reconstruct some unknowns (i.e. coefficients and source term) by means of inverting the map 
\begin{align*}
    \mathcal{F}: (\rho_0, \, k_0,\, J) \to \left[v(x_0,t), \, u(x_0,t;z)\right]_{t\in(0,T), \, z\in\Omega},
\end{align*}
which is non-linear. In view of the parameter count, we notice that our measurements depend on 4 parameters (1 for time and 3 for injected agents) and also the unknown functions (more specifically, the source term) are four dimensional. In this regard, the inverse problem under consideration here is not over-determined. In addition, although we need to move the position of the contrast agent, the position of the receiver $x_0 \in \partial \Omega$ is kept fixed. The needed time $T$ of measurement can be estimated as $T \ge \frac{\textnormal{diam}(\Omega)}{\inf_{x \in \Omega}c_0(x)}$. Therefore knowing an a-priori lower estimate of the speed of propagation $c_0(\cdot)$, we can estimate the needed measurement time $T$. 

\smallskip

The use of micro-scaled agents is already in practice  in various imaging techniques or therapeutic applications, see for instance \cite{C-F-Q2009, F-M-S2003, I-I-F2018, Q2007, Z-L-L-S-W19, L-C2015, Q-C-F2009}. The principal idea in this approach is to supply large contrast to the material properties by using small-scaled agents that serve the role of  resonators. This consideration may sidetrack the instability issues present in many traditional inverse problems where the medium is probed by a specific type of waves and the interaction of the medium is measured away from region of interest. The instability phenomenon there is hidden within the highly smoothing character of the forward problem. To remedy this situation, let us also point out that a different set of ideas like the regularization approach \cite{E-H-N1996} or hybrid imaging techniques \cite{B2013, B-U2013}, see also the monograph \cite{Ammari-imaging-book} for other methodologies, are well-known in the existing literature. For the case of  injecting small-scaled contrast agents, the scale between their size and contrasts are quite critical as they help generate the resonating frequencies for these inclusions as eigenvalues of integral operators present in the Lippmann Schwinger representation which governs the interaction of the medium with the applied source. A key feature of the current work is that the injected bubbles have contrasting properties in both mass density and bulk modulus as compared to the background medium. With such contrast, it is known that these bubbles resonate at the single Minnaert frequency (as a subwavelentgh resonance). Indeed, with such scales, it is shown, first in \cite{Ammari-et-al-Minnaert}, that the related (indirect) system of boundary integral equations is not injective at this frequency. Afterwords, it is proved in \cite{Li_Sini_Resolvent-Anal} that this frequency is a true mathematical resonance of the natural Hamiltonian given by the wave operator, i.e. the Minnaert frequency is a pole of the meromorphic extension of this operator's resolvent. In addition, the wave propagation effects governed by Minnaert resonance in the time domain is analyzed in \cite{ Li_Sini_Large time} where  the resonant expansion of the wave-field is derived, for a large time, and it is confirmed that the real part of the Minnaert frequency characterizes its period while its imaginary part fixes its life time. These last cited works are dealt for a homogeneous background. But as it shown in \cite{D-G-S2021}, for the time-harmonic setting, and in the current work, for the time-domain setting, this resonance appears in a variable background as well. Using this resonant effect, we could extract from the measured data, the internal values of the travel-time function and the internal values of the free wave-field generated in the absence of the bubble. From the travel-time function, we recover the speed of propagation, via the Eikonal equation, and then from the free wave-field, we recover the mass density and the source terms.

The mathematical attributes of using small-scaled contrast agents in inverse problems and imaging has already appeared in optics and photo-acoustic context \cite{G-S_ICCM2022, G-S_JMAA2022, G-S_JDE2022} and acoustic time harmonic \cite{D-G-S2021} and time domain \cite{Senapati_Sini_Wang_23, Sini_Wang_2022} setting along with the elastic time harmonic case \cite{C-G-S2023}. Our problem being in the time-domain acoustic setting, we highlight the key aspects of the current work and its distinguishing features from that of \cite{Senapati_Sini_Wang_23} and \cite{Sini_Wang_2022}. While the works \cite{Senapati_Sini_Wang_23} and \cite{Sini_Wang_2022} employ droplets enjoying very low bulk modulus and moderate mass density, we use bubbles in the present article which have both the bulk and mass density very small as compared to the background. Such contrast in material properties justifies the appearance of Newtonian operator in \cite{Senapati_Sini_Wang_23, Sini_Wang_2022} and the surface Double layer operator (i.e. the Neumann-Poincar\'e operator) in the current setting. In the process, a special sequence of eigenvalues for the Newtonian operator are generated in the first case, whereas we could only excite an isolated eigenvalue for the Double layer operator for the later (i.e. the so-called Minnaert frequency mentioned above). In \cite{Senapati_Sini_Wang_23}, using droplets as contrast agents, we could recover the bulk modulus and the source term assuming that the mass density is a known constant. We do believe that such restriction is not only technical but it is related to the fact that the droplet enjoy contrast in its bulk modulus but not in its mass density. In addition, the spherical-shape of the droplets was needed as we heavily used the asymptotic behavior of the sequence of the eigenvalues of the related Newtonian operator to generate a useful Riesz basis. In the present work, such a restriction on the shape is not needed and we can recover both the mass density, the bulk modulus and the source term. Therefore, we can say that the gas bubbles are superior as compared to the liquid droplets as far as the inverse acoustic problem is concerned.

Finally, let us mention that our work discusses reconstruction aspects to a non-linear parameter (and, source) identification problem for a linear acoustic model which is the linearized version of the bubbles model described in \cite{C-M-P-T1985, C-M-P-T1986}. Stating the corresponding inverse problems for the original nonlinear bubbles model is highly challenging but worth trying.

\subsection{Statement of the results}
Consider $w$ to be the contrast between the wave-fields after and before introducing the bubble within the medium at a position $z$, i.e. $w(x,t;z)= u(x,t;z) - v(x,t)$. To simplify the presentation, we again omit the argument $z$ from $w$ and $u$. Now, from \eqref{IVP_v} and \eqref{IVP_u}, it is immediate that $w$ satisfies the IVP 
\begin{align}\label{diff pde}
    \begin{cases}
     k_0^{-1}(x) w_{tt} - \textnormal{div}\left(\rho_0^{-1}(x) \nabla w\right) = \left(\frac{1}{k_0(x)} - \frac{1}{k(x)}\right) u_{tt} - \textnormal{div}\left(\left(\frac{1}{\rho_0(x)} - \frac{1}{\rho(x)}\right) \nabla u\right), \, \textnormal{ in } \, \mathbb{R}^3\times \mathbb{R}_+, \\
     w(\cdot,0) = w_t(\cdot,0) = 0, \textnormal{ in } \mathbb{R}^3,
    \end{cases}   
\end{align}
The non-constant wave-speed $c_0$ induces a Riemannian metric with arc-length $d\zeta$ given by 
\begin{equation}\label{metric}
d\zeta = \left( \sum_{i=1}^3 c_0^{-2}(x) (dx_i)^2 \right)^{1/2}.
\end{equation}
It is closely connected to the notion of travel-time function introduced in Subsection \ref{Preliminaries_Green func}. In all of the article, we presume the following regularity assumption of the geodesics given by the metric \eqref{metric} which ensures
\begin{align}\label{geodesic assumption}
   & \textit{Any pair of points } x,\,y \textit{ in } \overline\Omega \textit{ can be connected with a single geodesic line } \Gamma(x,\,y) \textit{ of the metric } \nonumber\\ 
   & \eqref{metric} \textit{ and the domain }  \Omega \textit{ is convex with respect to geodesics connecting two boundary points.}
\end{align}
This regularity of geodesics is an important assumption to be used in all of our analysis. We need this for solving an Eikonal equation of travel-time function given in \eqref{Eikonal eqn} along its rays which are nothing but the geodesics to the metric \eqref{metric}. Furthermore, the assumption \eqref{geodesic assumption} is key to derive a progressive wave expansion, which we discuss in Subsection \ref{Preliminaries_Green func}. A detailed account on this condition can be found in \cite[Chapter 3]{Romanov_book_1987}. For our analysis, the source function $J(\cdot,\cdot)$ should belong to suitable function spaces. For $p\in\mathbb{Z}_+$ and $T \in \mathbb{R}\cup\{\infty\}$, we introduce the usual real valued Sobolev space denoted by 
\begin{align*}
    H^p_0(0, T) = \left\{ f \vert_{(0,T)}; \, f \in H^p(\mathbb{R}) \textnormal{ and } f\vert_{(-\infty,0)}=0\right\},
\end{align*}
which can be similarly extended in the context of a Banach space say $\textnormal E$ by
\begin{align*}
    H^{p}_{0,\sigma} \left(0,T; \textnormal{E}\right):= \left\{f ; \, \, e^{-\sigma t}f\in\mathcal{S}^{'}_+(\textnormal{E}), \textnormal{ and } \|f\|_{H^{p}_{0,\sigma} \left(0,T; \textnormal{E}\right)} := \sum_{k=0}^{p} \int_0^T e^{-2\sigma t}\|\partial_t^{k} f(\cdot,t)\|^2_{\textnormal{E}} \, dt < \infty \right\}
\end{align*}
where $\sigma>0$ and $ \mathcal{S}^{'}_+(\textnormal{E}) $ denotes the space of $\textnormal{E}$-valued tempered distribution in $\mathbb{R}$ having support in $[0,\infty)$. 

\begin{theorem}\label{asymptotic expansion_w}
    With the condition (\ref{geodesic assumption}) being fulfilled, let us choose $x_0\in \partial\Omega$ and $T<\infty$. For $k_0 ,\, \rho_0 \in C^{17}(\mathbb{R}^{3})$ with $k_0 = \rho_0 = 1$ in $\mathbb{R}^3 \setminus \Omega$ and $J \in H^{7}_{0,\sigma}(\mathbb{R}_+;L^2(\mathbb{R}^3))$, the wave-field $w$ admits the expansion
   \begin{align}\label{asymptotic expansion_2}
   \nonumber  w(x_0,t) & = \frac{|D|}{2\pi A_{\partial D}} \frac{m(x_0,z)}{\rho_0(z) \zeta(x_0,z)} v(z,t-\zeta(x_0,z)) + \frac{2|D|}{\rho_0(z) A_{\partial D}} \int_\mathbb{R} g(x_0,t-\tau;z) v(z,\tau)\, d\tau \\
   \nonumber  &  - \frac{\sqrt{k_1}|D|}{\pi \sqrt{2} A^{3/2}_{\partial D}} \frac{m(x_0,z)}{\rho^{3/2}_0(z)\zeta(x,z)} \int_0^{t-\zeta(x_0,z)} \sin\left(\omega_M(z)(t-\zeta(x_0,z)-\tau)\right) v(z,\tau) \, d\tau\, \\
      &  - \frac{2\sqrt{2k_1}|D| }{A^{3/2}_{\partial D}} \frac{1}{\rho^{3/2}_0(z)} \, \int_{\mathbb{R}} g(x_0,t-\tau;z) \int_{0}^{\tau} \sin\left(\omega_M(z)(\tau-s)\right) v(z,s) \, ds\, d\tau + O(\epsilon^2) ,
   \end{align}
   which holds point-wise in $t\in(0,T)$. Here, $\omega_M(z)$ is given in \eqref{Minnaert} and, the functions $g$ and $m$ are defined in \eqref{defn_g_m}.
\end{theorem}
\bigskip

The quantity $A_{\partial D}$ is a positive geometric constant defined as 
$$A_{\partial D}:=  \frac{1}{\vert \partial D\vert}\int_{\partial D} \int_D \frac{dx}{|x-y|} \, dS_ydS_x \simeq \epsilon^2,$$
see Section \ref{theo_expansion}.  Consequently, the first four terms in \eqref{asymptotic expansion_2}, having order $\epsilon$ (in pointwise sense w.r.t time), are dominant parts of $w(x_0,t)$. To underline this aspect of the wave-field at $x_0$, we equivalently express the asymptotic expansion \eqref{asymptotic expansion_2} as 
\begin{align}\label{asymptotic expansion_2-equiv}
    w(x_0,t) = w_d(x_0,t) + O(\epsilon^2), 
\end{align}
to be understood point-wise in time, where the function $w_d(x_0,\cdot)$ precisely consists of first four terms from the right hand side of \eqref{asymptotic expansion_2}. Here, the Minnaert frequency for the bubble $D$ is given by 
\begin{align}\label{Minnaert}
    \omega_M(z):= \sqrt{\frac{2k_1}{\rho_0(z)A_{\partial D}}}
\end{align}
which crucially depends on the mass density of the background at $z$ and also we have $\omega_M(z) \simeq 1$, as $\epsilon \ll 1$. From \cite{D-G-S2021}, it is also clear, in the context of time harmonic setting, that $\omega_M(z)$ is encoded in the contrasted scattered waves. Although, our discussion is presented in a way that $x_0\in\partial\Omega$, we can modify our argument as in \cite{Mukherjee_Sini_acoustic cavitation} to encompass the case when $x_0\in\Omega$ and $d(x_0, D)$ is at least of the order $\epsilon$. Such estimates of the pressure wave-fields very near to the bubble (as a contrast agent) has several medical applications such as in sonar therapy where one needs to create enough (but not too much) amount of pressure only near the injected bubble, see \cite{Mukherjee_Sini_acoustic cavitation} for more discussion on this topic. Let us mention here that the characterization of the Minnaert resonance was first given in \cite{Ammari-et-al-Minnaert} when the bubble is embedded in a homogeneous background. Its dependency in terms of the local heterogeneity ($\rho(z)$) was first observed in \cite{D-G-S2021} (when the bubble is embedded in a heterogeneous background). These two works were stated in the time-harmonic regime. As we can see, this characterization plays a key role in the time-domain ultrasound imaging using bubbles. 
\bigskip

Now, we address the inverse problem of determining the medium properties and source function from the measured wave-field $w$ collected at a fixed boundary point $x_0$ and for large enough interval of time $(0,T)$.

\begin{theorem}\label{inverse problem}
    Under the same assumptions as in Theorem \ref{asymptotic expansion_w}, we have
    \begin{align}\label{recovering v}
         v(z,t) = \mathbb{A} w(x_0,\cdot)(t+\zeta(x_0,z)) + O(\epsilon)
    \end{align}
    to be understood in point-wise sense in $t\in(0,T)$, where $\mathbb{A}: L^2(0,T) \to L^2(0,T)$ is an invertible operator depending on the coefficients $\rho_0(\cdot)$ and $k_0(\cdot)$. For $p\ge7$ and $J\in H_{0,\sigma}^{p+1}\left(\mathbb{R}_+; L^2(\mathbb{R}^3)\right) \cap C^{p+3}(\mathbb{R}^3\times[0,\infty))$,  the function $w_d(x_0,\cdot)$ defined in \eqref{asymptotic expansion_2-equiv} inherits the property  
   \begin{align}\label{dominant term}
       w_d(x_0,t) = \begin{cases}
           0, &\textnormal{ for } t<\zeta(x_0,z), \\
           \simeq \left(t-\zeta(x_0,z)\right)^{p+3}, &\textnormal{ for } t\ge\zeta(x_0,z),
          \end{cases} 
   \end{align} 
   under the assumption that  $\partial_t^{p+1} J(z,0) \neq 0$.
\end{theorem}

\medskip
 
    The operator $\mathbb{A}$ is inverse to a second kind Volterra operator in $L^2(0,T)$, i.e. $\mathbb{A}:= \left(\alpha(x_0,z)\mathbb{I} + \mathbb{K}\right)^{-1}$. For a precise definition of the quantity $\alpha(x_0,z)$ and the operator $\mathbb{K}$, we refer to \eqref{alpha}. The operator $\mathbb{A}$ emerges quite naturally out of the expansion \eqref{asymptotic expansion_2}. Due to the scale present in the operator $\mathbb{A}$ (see \eqref{alpha} and \eqref{kernel}), we notice that the first term in right hand side of \eqref{recovering v} is of order one (in point-wise sense in time) and therefore, this is the dominant term of the initial wave-field $v(z,\cdot)$. As a consequence, we are able to recover the wave-field in a point-wise manner (modulo an error $O(\epsilon)$) once we know the coefficients $\rho_0$ and $k_0$ in that matter. This is because the operator $\mathbb{A}$ depends on medium properties (unlike the case in \cite{Senapati_Sini_Wang_23}). Here, we would like to underline that the operator $\mathbb{K}$ can be also inverted in regular spaces such as $C[0,T]$; see \cite{Senapati_Sini_Wang_23}. In light of Theorem \ref{inverse problem}, let us interpret the travel-time function. Taking only the dominant part of the measurement into account, we identify the quantity $\zeta(x_0,z)$ as the time level where it experiences jump  (in point-wise sense) for the first time. However, the function $w(x_0,\cdot)$ vanishes trivially before $t_*:=\inf_{y\in\overline D}\zeta(x_0,y)$, which readily follows when we invoke Theorem \ref{Green function} to the representation \eqref{representation of w} at $x=x_0$. In view of \eqref{dominant term} and the fact that $t_*= \zeta(x_0,z) + O(\epsilon)$, it is apparent that $w(x_0,\cdot)$ admits a strict sign after the time-level $t_{**}= \zeta(x_0,z)+ \epsilon^{1/(p+3)}$. Consequently, we see that the measured data i.e. $w(x_0,\cdot)$ determines the quantity $\zeta(x_0,z)$ upto an error of order $\epsilon^{1/(p+3)}$. After that, we use the Eikonal equation \eqref{Eikonal eqn} to determine $c_0(z)$ from $\zeta(x_0,\cdot)$ near $z$. Here we remark that the integer $p$ in Theorem \ref{dominant term} is to be understood as the first integer for which $\partial^{p+1}_t J(z,0)\neq 0$. However, this requirement only needs to be satisfied by the source function initially at the point $z$. In the context of the inverse problem, we assume the same $p\in\mathbb{N}$ works for all point $z\in\Omega$. In equivalent terms,  $\partial^{p+1}_t J(\cdot,0)$ is assumed to be non-zero throughout $\Omega$. To recover the medium properties, we consider two cases. The first one assumes either one of $k_0$ and $\rho_0$ is known, whereas, in the second case, the function $\partial^{p+1}_t J(\cdot,0)\vert_{\Omega}$ is assumed to be known. Indeed, the later case is more general compared to the first one. After recovering the material parameters, we use \eqref{recovering v} to recover $v(\cdot, \cdot)$ and then the IVP \eqref{IVP_v} to determine $J(\cdot,\cdot)$ in $\Omega \times (0, T)$. In the following, we present a reconstruction scheme for the general case.    
\bigskip

\textbf{Reconstruction scheme.} We reconstruct the medium properties and source via the following steps:

\vspace{1mm}

\textbf{Step 1:} Consider the data as the difference of wave-fields measured in time at $x_0\in\partial\Omega$ before and after injecting the bubble at $z\in\Omega$ and identify the time level when the data experiences jumps for the first time. This time level determines the travel time function $\zeta(x_0,z)$.

\vspace{1mm}

\textbf{Step 2:} Inject the bubble at different points near to $z$ and perform Step 1 to extract information the function $\zeta(x_0,\cdot)$ in a neighbourhood of $z$. Determine $c_0(z)$ i.e. the wave speed at $z$ by substituting $\zeta(x_0,\cdot)$ into the Eikonal equation \eqref{Eikonal eqn}. Likewise, determine the wave speed for other points in $\Omega$ and hence the function $c_0(\cdot)$, which determines the travel-time function $\zeta(\cdot,\cdot)$ from \eqref{Eikonal eqn} and therefore the function $\eta(\cdot,\cdot)$ from \eqref{Normal coordiante}. 

\vspace{1mm}

\textbf{Step 3:} Update the data by considering its dominant part that is $w_d(x_0,\cdot)$. From the following equality, determine the coefficient $\alpha(x_0,z)\, m(z,z)$ by considering the right-sided limit in $(p+3)$-th time-derivative of $w_d(x_0,t)$ at the time $t=\zeta(x_0,z)$,
\begin{align*}
    \lim_{t\to\zeta(x_0,z)+}\partial^{p+3}_t w_d(x_0,t) =  \frac{\alpha(x_0,z)\, m(z,z)}{4\pi}\,\partial^{p+1}_t J(z,0) \,\lim_{t\to 0+}\partial_t^{p+3} \int_{\Omega_t} \frac{\left(t-\zeta(z,y)\right)^{p+1}}{\zeta(z,y)} \, dy.
\end{align*}
 With the help of \eqref{alpha and m}, recover the bulk modulus $k_0(\cdot)$ from the quantity $\alpha(x_0,\cdot)\, m(\cdot,\cdot)$ from \eqref{alpha and m} where we utilize the knowledge of $\eta(\cdot,\cdot)$, and then the matrix $\nabla_x \eta(x_0,\cdot)$, derived in Step 2. This thereafter leads to the determination of $\rho_0(\cdot)$ as we have already determined the wave speed. At this step, we crucially use prior knowledge of the non-vanishing function $\partial^{p+1}_t J(\cdot,0)$.

\vspace{1mm}

\textbf{Step 4:} From the knowledge of $k_0$ and $c_0$, compute the operator $\mathbb{A}$ and employ it to find the wave-field $v$ from the equation \eqref{recovering v}. Inserting $k_0, \, \rho_0$ and $v$ in \eqref{IVP_v}, determine the source function $J$.
\bigskip

Taking the dimensionality into consideration, we have seen earlier that the inverse problem is not over-determined, yet we can recover all the coefficients in PDE \eqref{IVP_v}. To the best of our knowledge, simultaneous recovery of two medium properties and source function in the time-domain setting has not been achieved so far in the literature which is the principal theme of this work. In addition to this, we are able to develop a reconstruction procedure pertaining to this inverse problem, which indeed is another appealing aspect of our result. In our approach, the prior knowledge of the source function at initial time is needed to split up the reconstruction process of two medium properties. Due to this reason, we believe that our result for the inverse problem under consideration here is optimal.  Furthermore, we notice that our reconstruction method has a {\it local} feature embedded in it. It indicates that one needs to perform the experiment for points close enough to $z$ to extract the information of unknown function at $z$. This local-feature is often desired in many practical applications. At this point, we mention that we can easily adapt our reconstruction scheme to discuss an imaging technique where one aims to determine material properties by probing the medium with an active or known source. Such an active source is applied before and after injecting each bubble and then  measure, for each applied source, the generated pressure field at the single point on $\partial \Omega$. With regard to this, we notice that Step 4 in our reconstruction scheme is not needed as it was meant to recover the source which, in this case, is assumed to be known.. We discuss more on this in Section \ref{active source}.

Let us briefly expound on the mathematical aspects of the current work. To deal with the inverse problem of reconstructing the medium properties and source function, we consider, as measurement, the wave-field generated by the medium in the presence of source function when a small-scaled contrast agent (bubbles in our case) is injected to the medium. To derive the asymptotic profile of this wave-field, we us time-dependent integral equation methods and key properties of the Newtonian as well as the Magnetization-type operators appearing in the related Lippmann-Schwinger equation. Note that the whole asymptotic analysis, apart from the first a-priori estimate, is done in the time domain. This avoids the systematic use of Fourier-Laplace transform which imposes rather formal computations as we have no control of the high frequency components. In performing the asymptotic analysis in time domain, the spectral decomposition of magnetization type operator is crucial to deduce a-priori estimate for the wave-field. Furthermore, the structure of Green's function to the hyperbolic operator in \eqref{IVP_v} in conjugation with the travel time function has been heavily exploited in our analysis. To carry out these steps, we allow the medium properties or the source to be enough regular. With the consideration of time-jump of the data, we reduce the multidimensional inverse problem to an one dimensional inverse problem which we address by analyzing a second kind Volterra integral operator. Such reduction is possible when we have already identified the medium properties. To this effect, the knowledge of the source function at the initial time is of importance which we presume to be known a-priori. As mentioned in our discussion above, this assumption looks to be optimal.      

Although our proposed method can simultaneously determine two coefficients appearing in the leading order coefficients a PDE and a causal source function (under an assumption), we require sufficient regularity for those unknowns; see Theorem \ref{asymptotic expansion_w} and \ref{inverse problem}. In view of applications, this seems to be quite constraining. To underline the need of such assumption in our work, let us mention that the high regularity of the medium properties is used to guarantee the required smoothness needed, in our context, for the {\it regular} part of Green's function, which we discuss in Subsection \ref{Preliminaries_Green func}. In addition, we require the causal source function to be sufficiently regular and concurrently non-vanishing at initial time. This prompts us to the consideration $p\ge 7$ in Theorem \ref{inverse problem}. The high time-regularity of the source is utilized to obtain some a-priori estimate for the auxiliary wave-field and its time derivatives, whereas the non-vanishing character of the source at initial time is reflected in the measured data leading us to determination of the travel time function. Here we also remark that we reconstruct the wave speed by differentiating an Eikonal equation satisfied by the travel-time function which is recovered upto an error in Step 1 of the reconstruction scheme given before. Similarly in Step 3, we considered the time gradient of dominant part of measurement as a part of data. Needless to say, in both these steps, one encounters numerical error.

Let us finish this introduction by briefly discussing the practicability of the proposed framework. In order to align this framework to more realistic experiments of introducing microscaled agents to the medium, we might need some adjustments while considering the model described by IVP \eqref{IVP_u}. It might appear that the IVP \eqref{IVP_u} requires us to reproduce the same internal source function $J(\cdot,\cdot)$ every time we insert micro-scaled bubbles into the medium. This is a similar situation as in inverse problems based on the dynamic Dirichlet-Neumann map where one needs to repeat collecting the boundary data keeping the same space-time variable coefficients. Such assumptions might be supported by practical considerations. In our case, another attempt to deal with this issue is to assume the source function to be repeating in time with period $T$ which is large enough (but finite). We briefly discuss on this now. In this setting, we first fix $x_0\in\partial\Omega$ and record the interaction of the medium with the source $J$ for the interval $[0,T]$. This is given by the function $v(x_0,t)\big\vert_{t\in(0,T)}$, where $v$ is given by the IVP \eqref{IVP_v}. At $t=T$, we introduce small-scaled bubble at $z\in\Omega$ and the perturbed medium has the properties as in \eqref{scales of bulk modulous and density}. Now we observe the generated wave at $x_0\in\partial\Omega$ for the interval $(T, 2T)$ i.e. $\mathfrak{u}(x_0,t;z)\big\vert_{t\in(T,2T)}$, where $\mathfrak{u}$ solves the following IVP
\begin{align}\label{IVP_u_new}
    \begin{cases}
    k^{-1}(x) \partial_t^2\mathfrak{u} - \textnormal{div}\left(\rho^{-1}(x) \nabla \mathfrak{u}\right) = J(x,t),  \,  \textnormal{ in } \mathbb{R}^3 \times (T, \infty),\\  
    \mathfrak{u}\vert_{+} = \mathfrak{u}\vert_{-} , \quad \rho_1^{-1} \partial_\nu \mathfrak{u}\vert_{+} = \rho^{-1}_0 \partial_\nu \mathfrak{u}\vert_{-} \quad \textnormal{ on } \partial D(z)\times(T, \infty), \\
    \mathfrak{u}(\cdot,T) = v(\cdot,T), \quad \mathfrak{u}_t(\cdot,T) = v_t(\cdot,T), \, \textnormal{ in } \mathbb{R}^3,
   \end{cases}
\end{align}
Having defined $\mathfrak{w}(x,t) := \mathfrak{u}(x,t)- v(x,t)$ and introducing the translation operator 
$$\widetilde f(t) = f(t+T)$$ 
for a time dependent function $f$, we recast the PDE \eqref{IVP_u_new} into
\begin{align*}
\begin{cases}
    k_0^{-1}(x) \partial_t^2\widetilde{\mathfrak{w}} - \textnormal{div}\left(\rho_0^{-1}(x) \nabla \widetilde{\mathfrak{w}}\right) = \left(\frac{1}{k_0(x)} - \frac{1}{k(x)}\right) \widetilde{\mathfrak{u}}_{tt} - \textnormal{div}\left(\left(\frac{1}{\rho_0(x)} - \frac{1}{\rho(x)}\right) \nabla \widetilde{\mathfrak{u}}\right), \, \textnormal{ in } \, \mathbb{R}^3\times \mathbb{R}_+,\\
    \widetilde{\mathfrak{w}}(\cdot,0) = \partial_t \widetilde{\mathfrak{w}}(\cdot,0)=0
\end{cases}
\end{align*}
which is similar to \eqref{diff pde}. Under the assumption that we can get rid of the micro-scaled agent at $z$ to initiate the process of introducing the agent at some different point, we see, the dataset $\left[v(x_0,t), \, \mathfrak{u}(x_0,t+T;z)\right]_{t\in(0,T),\, z\in\Omega}$ is equivalent to the dataset $\left[v(x_0,t), \, u(x_0,t;z)\right]_{t\in(0,T),\, z\in\Omega}$. Here the time-periodicity of $J$ is used. Therefore, our arguments for the inverse problem based on the IVPs \eqref{IVP_v} and \eqref{IVP_u} can be applied to deal with the current framework that describes more realistic experiments. Furthermore, it is worthwhile to consider the inverse problem for the case when we inject multiple or a cluster of bubbles, at once, and then record the wave-field generated by the interaction of this perturbed medium with the source. Such models are recently considered in \cite{Mukherjee_Sini_multiple inclusion} and \cite{Mukherjee_Sini_dispersive case} and can avoid the assumption which requires repeating the experiments of injecting bubbles discretely at different points.

Throughout the article, we follow the notation, $f(x) = O(\epsilon^r)$ to imply that there is $C>0$ independent of $\epsilon$ so that $|f(x)| \le C \epsilon^r$ holds true in the domain of definition of $f$. Also, we use the notation $ a\preceq b$ (or, $a \succeq b$) and $a \simeq b$ to imply $|a| \le C \, b$ (or, $ |a| \ge C\,b$) and $|a| = C b$ respectively for some $C>0$ independent of $a, \, b$.


\section{Preliminaries}

\subsection{Structure of the Green's function}\label{Preliminaries_Green func}

In this section, we present a brief account on the structure of Green's function related the hyperbolic operator consider in \eqref{IVP_v}. The structure is important in our analysis. The exposition here is mainly borrowed from \cite{Romanov_book_1987} that is based on progressive wave expansion. Let us mention that the condition \eqref{geodesic assumption} is needed for the proof of the singularity structure of the Green's function. 
\bigskip

For $(y,\tau)\in \mathbb{R}^n \times \mathbb{R}$, let us introduce $G(\cdot,\cdot;y,\tau)$ which solves 
\begin{equation}\label{v_Green_H}
\begin{cases}
  k_0^{-1}(x) G_{tt} - \textnormal{div}\left(\rho_0^{-1}(x) \nabla G\right) = \delta(t-\tau,x-y),  &  \mathrm{in}\;\mathbb R^3\times \mathbb R,\\
  G|_{t<\tau}=0  &  \mathrm{in}\;\mathbb R^3.
\end{cases}
\end{equation}
The next result is the singularity structure of $G(x,t;y,\tau)$ which follows from \cite{Romanov_book_1987, Romanov_}. At this stage, let us stress that the use of PDE \eqref{v_Green_H} or integrating some expressions relating to the distribution $G$ in our various parts of analysis should be taken into consideration, as usual, by deleting a small neighborhood around the source $(y,\tau)$ followed by a limiting process. Without repeating these comments, we will always follow this interpretation throughout this article, as it is quite standard.
\begin{theorem}\label{Green function}
     Assume that the condition \eqref{geodesic assumption} is satisfied. For $k_0,\, \rho_0 \in C^{17}(\mathbb{R}^3)$, we have
    \begin{align*}
        G(x,t;y,\tau) = \frac{m(x,y)}{4\pi\zeta(x,y)}\,\delta(t-\tau-\zeta(x,y)) + g(x,t-\tau;y) , \quad \textnormal{ for } t > \tau,
    \end{align*}
    for some $\zeta^2(\cdot,\cdot)\in C^{16}(\Omega \times \Omega), \, m(\cdot,\cdot)\in C^{15}(\Omega \times \Omega)$ and $g(x,t-\tau;y)=0$ for $t < \tau + \zeta(x,y)$. In the conoid $t \ge \zeta(x,y)$, we have continuity of the function $\partial^\alpha_x \partial^\beta_t\partial^\gamma_y g(x,t;y)$, where $\alpha$ and $\gamma$ are multindices satisfying $|\alpha|+\beta+|\gamma| \le 2$.
\end{theorem}
 Here, $\zeta(x,y)$ denotes the travel-time function which measures the time needed for a signal originating at $y$ to reach $x$  in the medium $\Omega$ having sound speed $c_0(\cdot)$. It can be also viewed as the Riemannian distance function induced by the metric \eqref{metric}. The functions $m(x,y)$ and $g(x,t-\tau;y)$ are defined in \eqref{defn_g_m}.

\medskip

\textbf{Proof:} 
 We outline the proof here. Let us get started with the function $H(\cdot,\cdot;y,\tau)$ which solves the IVP 
\begin{align}\label{v_Green_G}
\begin{cases}
    H_{tt} - \textnormal{div}_x\left(c_0^2(x) \nabla_x H\right) + \rho_0^{-1}(x) \nabla_x k_0(x)\cdot\nabla_x H= \delta(t-\tau,x-y), &  \mathrm{in}\;\mathbb R^3\times \mathbb R,\\
  H|_{t<\tau}=0  &  \mathrm{in}\;\mathbb R^3.
\end{cases}
\end{align}
From the definition \eqref{speed}, it is obvious to note that the function $H(\cdot,\cdot;y,\tau)$ is introduced in a way that  
\begin{align*}
   G(x,t;y,\tau) = k_0(y) \, H(x,t;y,\tau)
\end{align*}
which can be justified by multiplying both sides of the equation \eqref{v_Green_H} by the bulk modulus resulting into \eqref{v_Green_G}. Therefore, it suffices to derive a structure for $H(\cdot,\cdot;y,\tau)$ which we borrow from \cite{Romanov_}. See also \cite[Theorem 4.1]{Romanov_book_1987}. The travel time function $\zeta(x,y)$ solves the Eikonal equation
\begin{align}\label{Eikonal eqn}
   \begin{cases}
      c_0^2(x)\,|\nabla_x \zeta(x,\,y)|^2 = 1, \\
      \zeta(x,\,y) = O(|x-y|), \quad \textrm{as}\; x\to y.
   \end{cases}
\end{align}
 For the solvability of \eqref{Eikonal eqn} in $\Omega$ and the progressive wave expansion of $G$, we need the regularity assumption of geodesics in \eqref{geodesic assumption}. For a detailed discussion, we ask the reader to consult \cite[Chapter 3 and 4]{Romanov_book_1987}.

 With the regularity assumption on geodesics in \eqref{geodesic assumption}, we can assign Riemannian coordinates $\xi= (\xi_1,\xi_2,\xi_3)$ to each point $x$ for fixed $y$ which can be computed from the knowledge of the travel-time function. In this regard, we recall from \cite{Romanov_} that 
\begin{align}\label{Normal coordiante}
    \xi = - \frac{1}{2} c_0^{2}(y) \, \nabla_y \zeta^2(x,y) := \eta(x,y),
\end{align}
and the function $\xi = \eta(x,y)$ admits the inverse say $x= f(\zeta, y)$ which enables us to locate $x$ in terms of the coordinate $\xi$. Furthermore, following the arguments from Step I in the proof of Lemma \ref{H^2 estimate} or \cite[Lemma 3.3]{Romanov_book_1987}, we can derive the relation 
\begin{align*}
  \eta(x,y) = x-y + O(|x-y|^2), \quad  \textrm{as}\; x\to y, 
\end{align*}
which by means of the Eikonal equation \eqref{Eikonal eqn} and \eqref{Normal coordiante} translates to 
\begin{align}\label{reln_travel time}
   \zeta(x,y) = c_0^{-1}(y) |\eta(x,y)| = c_0^{-1}(y)|x-y| + O(|x-y|^2).
\end{align}

With reference to \cite{Romanov_book_1987}, we express the ansatz for $H(x,t;y,\tau)$ as 
\begin{align*}
    H(x,t;y,\tau) = \frac{\sigma(x,y)}{2\pi c_0^3(y)}  \delta\left((t-\tau)^2 - \zeta^2(x,\,y)\right) + h(x,t-\tau;y), 
\end{align*}
where the function $h(x,\,t-\tau;\, y)$ is supported for points satisfying $t>\tau+\zeta(x,y)$, and 
\begin{align*}
    \sigma(x,y) & = |\nabla_x \eta(x,y)|^{1/2}\exp\left(\frac{1}{2}\int_{\Gamma(x,y)} \sum_{i=1}^3 c_0^{-2}(\xi)\,\rho^{-1}_0(\xi)\,\partial_i k_0(\xi)\, d\xi_i \right) \\
    & = |\nabla_x \eta(x,y)|^{1/2}\exp\left(\frac{1}{2}\int_{\Gamma(x,y)} \langle \nabla_\xi \log k_0(\xi), d\xi\rangle \right) \\
    & = |\nabla_x \eta(x,y)|^{1/2}\exp\left(\log \sqrt{\frac{k_0(y)}{k_0(x)}}\right)  = |\nabla_x \eta(x,y)|^{1/2} \sqrt{\frac{k_0(y)}{k_0(x)}}.
\end{align*}
As we are interested in $t\ge\tau\ge0$, we find
\begin{align*}
    \delta\left((t-\tau)^2 - \zeta^2(x,\,y)\right)  = \frac{1}{2 \zeta(x,y)} \Big(\delta(t-\tau-\zeta(x,y)) + \delta(t-\tau+\zeta(x,y)\Big)  =  \frac{1}{2 \zeta(x,y)} \delta(t-\tau-\zeta(x,y)),
\end{align*}
implying, for $t\ge\tau$,
\begin{align} \label{G_expression}
   \nonumber G(x,t;y,\tau) = k_0(y) H(x,t;y,\tau)  & = \frac{1}{{4\pi}} \frac{\sigma(x,y)k_0(y)}{c_0^3(y)\zeta(x,y)} \delta(t-\tau-\zeta(x,y)) + k_0(y) \, h(x,t-\tau;y) \\
    & = \frac{m(x,y)}{4\pi\zeta(x,y)} \delta(t-\tau-\zeta(x,y)) + g(x,t-\tau;y),
\end{align}
where we denote 
\begin{align}\label{defn_g_m}
    g(x,t-\tau;y) := k_0(y)h(x,t-\tau;y), \, \textnormal{ and } \, m(x,y) := \sigma(x,y)k_0(y) c_0^{-3}(y).
\end{align}
Note that, the support condition is translated to the function $g$.  Here, we would like to point out that the singularity structure of $H$ is also valid for medium properties which have lesser regularity as compared to the ones mentioned in Theorem \ref{Green function}. See \cite[Theorem 4.1]{Romanov_book_1987}. To meet our particular interest for this work, we require $g(x,t-\tau;y)$ and therefore $h(x,t-\tau;y)$ to be Lipschitz in the $y$ variable and twice differentiable with respect to $x$ variable in the interior of the characteristic conoid $t>\tau+\zeta(x,y)$. For this reason, we need $k_0, \, \rho_0 \in C^{17}(\mathbb{R}^3)$. For more details on this, we refer to \cite[Theorem 4.1]{Romanov_book_1987} for the case when $s=3$ and $l=13$.

\subsection{A-priori estimates}

It is evident that the function $u$ solving the IVP \eqref{IVP_u} depends on the scaling parameter $\epsilon$. As a first step, we need to quantify this dependency through a-priori estimates. These a-priori estimates are crucial in deriving the asymptotic expansion for our measured data or, the wave-field $u$. Our arguments are somewhat along the lines of \cite{Mukherjee_Sini_acoustic cavitation} and \cite{Senapati_Sini_Wang_23} but with the difference that both the bulk modulus and the density functions here are non-constants within $\Omega$. This inhomogeneity present in medium properties requires us for a local comparison of the travel-time function with its Euclidean analog, which crucially relies on the regularity assumption \ref{geodesic assumption}. Our arguments for obtaining the a-priori estimates are based on Laplace transform technique. However, the estimates which we obtain initially seem insufficient to derive Theorem \ref{asymptotic expansion_w}. To improve those, we use of the spectral theory of the magnetization operator defined in \eqref{magnetization}.

\begin{theorem}\label{gradient estimate}
    For  $k_0,\, \rho_0 \in C^{6}(\mathbb{R}^3)$, and $J\in H^p_{0,\sigma}(\mathbb{R}_+;L^2(\mathbb{R}^3)$, we have $u\in H^{p+1}_{0,\sigma}(\mathbb{R}_+; L^2(\mathbb{R}^3))$ solving \eqref{IVP_u} which satisfies the estimate
    \begin{align*}
       \|u\|_{H^{p}_{0,\sigma}\left(\mathbb{R}_+; H^1(D)\right)} \preceq \epsilon\, \|J\|_{H^{p}_{0,\sigma}\left(\mathbb{R}_+; L^2(\mathbb{R}^3)\right)}.
    \end{align*}
    Furthermore, for $p \ge 7$, the solution $u$ admits the bounds
    \begin{align}\label{a-priori estimate}
        \|u(\cdot,t)\|_{H^1(D)} + \|\partial_t^k\Delta u(\cdot,t)\|_{L^2(D)} \preceq \epsilon^{3/2}, \quad \|\partial_t^k\nabla u(\cdot,t)\|_{L^2(D)} \preceq \epsilon^{5/2}, \quad \|\partial^l_t\partial_\nu u(\cdot,t)\|_{L^2(\partial D)} \preceq \epsilon^2,
    \end{align}
    for a.e. $t \in (0,T)$, $k\in\{0,1,2\}$ and $l\in\{0,1,2,3\}$.
\end{theorem}
\textbf{Proof.} 
For better readability, we divide the proof into three steps. In Step I and II, we first apply the Lax-Milgram theorem and Laplace transform techniques to obtain a preliminary result on a-priori estimates of $u$. In Step III, we show how they can be improved by considering the magnetization operator \eqref{magnetization}.
\medskip

\textbf{Step I:}
Let us consider the elliptic problem of finding $\tilde u \in H^1(\mathbb{R}^3)$ such that
\begin{align}\label{elliptic pde}
    -\textnormal{div}(\rho^{-1}(x)\nabla \tilde u) + \frac{s^2}{k(x)} \tilde u  = \tilde{J}(x,s), \quad s= \sigma + i\mu, \mbox{ with } \sigma>0.
\end{align}
As one might expect, the PDE \eqref{elliptic pde} is nothing but the Laplace transformed version of the problem \eqref{IVP_u}. The variational formulation of \eqref{elliptic pde} is presented by introducing the sesquilinear form 
\begin{align*}
    \mathcal{B}(m,n) = \int_{\mathbb{R}^3} \rho^{-1}(x) \nabla m(x) \cdot \overline{\nabla n(x)} \, dx \, + s^2 \int_{\mathbb{R}^3} k^{-1}(x) m(x)\overline{n(x)}\, dx, \quad \textnormal{ for } m, n\in H^1(\mathbb{R}^3).
\end{align*}
The quadratic form $\mathcal{B}(\cdot,\cdot)$ is coercive for $s$ with to small $\mu$, but it is not necessarily in general. To remedy this situation, we rather work with $\mathcal{B}(\cdot,s\,\cdot)$ which exhibits coercivity for any $s$ as follows
\begin{align*}
     \mathcal{B}(m,sm) = \overline{s} \int_{\mathbb{R}^3} \rho^{-1}(x) |\nabla m(x)|^2 \, dx \, + s |s|^2 \int_{\mathbb{R}^3} k^{-1}(x) |m(x)|^2 \, dx
\end{align*}
which yields
\begin{align*}
     \mathfrak{Re}\left(\mathcal{B}(m,sm)\right) \succeq \sigma \min\{1,\sigma^2\} \|m\|^2_{H^1(\mathbb{R}^3)}.
\end{align*}
Therefore, we define $\tilde u$ to be a solution to \eqref{elliptic pde} if it satisfies 
\begin{align}\label{variational formulation}
    \mathcal{B}(\tilde u(\cdot,s),\, s\phi) = \overline{s} \int_{\mathbb{R}^3} \tilde J(x,s) \, \overline{\phi}(x), \quad \forall \phi\in H^1(\mathbb{R}^3).
\end{align}
 We further note that, the transmission conditions in \eqref{IVP_u} are encoded within the variational characterization of $\tilde u$. The existence and uniqueness of $\tilde u(\cdot,s)\in H^1(\mathbb{R}^3)$ satisfying \eqref{variational formulation} can be established by virtue of the Lax-Milgram theorem which also implies the equality
\begin{align}\label{variational bound}
   \overline{s} \int_{\mathbb{R}^3} \rho^{-1}(x) |\nabla \tilde{u}(x,s)|^2 \, dx \, + s |s|^2 \int_{\mathbb{R}^3} k^{-1}(x) |\tilde{u}(x,s)|^2 \, dx = \mathcal{B}(\tilde u,s\tilde u) = \overline{s} \int_{\mathbb{R}^3} \tilde{J}(x,s) \overline{\tilde u}(x,s) \, dx. 
\end{align}
Now, we use Cauchy-Schwarz inquality to obtain
\begin{align*}
    \left\vert \overline{s} \int_{\mathbb{R}^3} \tilde{J}(x,s) \overline{\tilde u}(x,s) \, dx \right\vert \preceq \|\tilde J(\cdot,s)\|_{L^2(\mathbb{R}^3)} \left( |s|^2 \int_{\mathbb{R}^3} k^{-1}(x) |\tilde{u}(x,s)|^2 \, dx \, + \int_{\mathbb{R}^3} \rho^{-1}(x) |\nabla \tilde{u}(x,s)|^2 \, dx \right)^{1/2}
\end{align*}
which is employed to the real parts of \eqref{variational bound} to estimate
\begin{align*}
    |s|^2 \int_{\mathbb{R}^3} k^{-1}(x) |\tilde{u}(x,s)|^2 \, dx \, + \int_{\mathbb{R}^3} \rho^{-1}(x) |\nabla \tilde{u}(x,s)|^2 \, dx \preceq \frac{1}{\sigma^2} \int_{\mathbb{R}^3} |\tilde J(x,s)|^2 \, dx.
\end{align*}
Due to the scales of $\rho$ and $k$ given in \eqref{scales of bulk modulous and density}, we find 
\begin{align}\label{estimates for laplace transformed solution}
    \|\tilde u(\cdot,s)\|_{L^2(D)} \preceq \frac{\epsilon}{\sigma|s|^2} \|\tilde J(\cdot,s)\|_{L^2(\mathbb{R}^3)}, \quad \|\nabla\tilde u(\cdot,s)\|_{L^2(D)} \preceq \frac{\epsilon}{\sigma} \|\tilde J(\cdot,s)\|_{L^2(\mathbb{R}^3)}.
\end{align}
\medskip

Let us now define 
\begin{align}\label{laplace trans soln}
    u(x,t) = \int_{\sigma-i\infty}^{\sigma+i\infty} e^{st} \tilde u(x,s) \, ds = e^{\sigma t}\int_\mathbb{R} e^{i\mu t} \tilde u(x,\sigma+i\mu) \, d\mu.
\end{align}
The definition \eqref{laplace trans soln} is independent of the choice of $\mathfrak{Re}(s)=\sigma>0$ and the function $u$ becomes causal, meaning $u(x,t) = 0$ for $t\le 0$. This can be proved by standard argument involving contour integration techniques. For details, we refer the reader to \cite{Rudin_func anal_book} and \cite{Sayas_book}. It is straightforward to check that $u(x,t)$ defines a weak solution for the IVP \eqref{IVP_u}, since for a.e. $t\in\mathbb{R}_{+}$ and $\psi\in H^1(\mathbb{R}^3)$, the formulation \eqref{variational formulation} yields
\begin{align*}
    & \int_{\mathbb{R}^3} k_0^{-1}(x) \partial_t^2 u(x,t) \overline{\psi}(x) \, dx +  \int_{\mathbb{R}^3} \rho_0^{-1}(x) \nabla u(x,t) \cdot \overline{\nabla \psi}(x) \, dx \\
    & = \int_{\sigma+i\mathbb{R}}^{} e^{st} \int_{\mathbb{R}^3} \left(k_0^{-1}(x) \partial_t^2 \tilde{u}(x,t) \overline{\psi}(x) \, dx +  \int_{\mathbb{R}^3} \rho_0^{-1}(x) \nabla\widetilde{u}(x,t) \cdot \overline{\nabla \psi}(x) \, dx \right) \, ds \\
    & = \int_{\sigma+i\mathbb{R}}^{} e^{st} \int_{\mathbb{R}^3} \tilde J(x,s) \overline{\psi}(x) \, dx \, ds
    \\
    & = \int_{\mathbb{R}^3} J(x,t) \overline{\psi}(x) \, dx.
\end{align*}
It is evident from \eqref{laplace trans soln} that
\begin{align*}
    \mathcal{F}_{t \to \mu} \left(e^{-\sigma t} \partial^k_t u(x,t)\right) = s^k \tilde{u}(x,s), \quad s = \sigma + i \mu,
\end{align*}
where $\mathcal{F}_t$ denotes Fourier transform w.r.t time variable. Consequently, we deduce 
\begin{align*}
    \|u\|^2_{H^{p+1}_{0,\sigma}\left(\mathbb{R}_+; L^2(D)\right)} := \int_0^{\infty} e^{-2\sigma t} \sum_{k=0}^{p+1} \|\partial^k_t u(\cdot,t)\|^2_{L^2(D)} \,dt 
  & \preceq \sum_{k=0}^{p +1} \int_{\sigma+i\mathbb{R}} |s|^{2k} |\tilde{u}(\cdot,s)|^2_{L^2(D)} \, ds \\
  & \preceq \epsilon^2 \|J\|^2_{H^{p}_{0,\sigma}\left(\mathbb{R}_+; L^2(\mathbb{R}^3)\right)},
\end{align*}
where the estimate \eqref{estimates for laplace transformed solution} is used. Continuing the same argument, we can further conclude 
\begin{align}\label{Basic-estimate}
    \|u\|_{H^{p}_{0,\sigma}\left(\mathbb{R}_+; H^1(D)\right)} \preceq \epsilon\, \|J\|_{H^{p}_{0,\sigma}\left(\mathbb{R}_+; L^2(\mathbb{R}^3)\right)}.
\end{align}

\medskip

\textbf{Step II:}
In order to improve the estimate (\ref{Basic-estimate}), we consider the PDE for the contrast $w:= u -v$. Similar to Step I, we  analyze the problem in the Laplace domain and then return back to the time domain. Borrowing the notations from Step I, it evidently follows from \eqref{IVP_v} and \eqref{IVP_u} that
\begin{align*}
    -\textnormal{div}\left(\rho^{-1}(x) \nabla \tilde w\right) + \frac{s^2}{k(x)} \tilde w  = F(x,s):= s^2\left(k_0^{-1} - k^{-1}\right) \tilde v - \textnormal{div}\left( \left(\rho_0^{-1} - \rho^{-1}\right) \nabla \tilde v\right).
\end{align*}
We identify $\tilde w$ as the unique solution satisfying the variational formulation
\begin{align*}
   \overline{s} \int_{\mathbb{R}^3} \rho^{-1}(x) \nabla \tilde w(x) \cdot \overline{\nabla \phi(x)} \, dx \, + s |s|^2 \int_{\mathbb{R}^3} k^{-1}(x) \tilde w(x)\overline{\phi(x)}\, dx = \langle F(\cdot,s), \phi\rangle_{H^{-1}(\mathbb{R}^3), H^1(\mathbb{R}^3)}
\end{align*}
for all $\phi \in \mathbb{R}^3$. We also notice that
\begin{align}\label{H^{-1} bound}
   \nonumber  \mathfrak{Re}\left(\langle F(\cdot,s), \phi\rangle_{H^{-1}(\mathbb{R}^3), H^1(\mathbb{R}^3)}\right) & = \mathfrak{Re} \left(s^2 \int_{D} \left(k_0^{-1}(x) - k^{-1}(x)\right) \tilde v(x,s) \overline{\phi}(x) \, dx\right) \\
   \nonumber & \quad -  \mathfrak{Re} \left(\int_{D} \left(\rho_0^{-1}(x) - \rho^{-1}(x)\right) \nabla \tilde v(x,s) \cdot \overline{\nabla \phi}(x) \, dx\right) \\
    & \preceq \epsilon^{-2} \max\{\|\nabla \tilde v(\cdot,s)\|_{L^\infty(D)}, \, |s|^2\|\tilde v(\cdot,s)\|_{L^\infty(D)}\} |D|^{1/2} \|\phi\|_{H^1(D)}.
\end{align}
The functions $\tilde v(\cdot,s)$ and $\nabla\tilde v(\cdot,s)$ being bounded in $D$,  the variational formulation, in combination with the estimate \eqref{H^{-1} bound}, gives us
\begin{align*}
   \mathfrak{Re}\left( s |s|^2 \int_{D} k^{-1}(x) |\tilde{w}(x)|^2 \, dx \, + \overline{s} \int_{D} \rho^{-1}(x) |\nabla \tilde{w}(x)|^2 \, dx \right) \preceq \epsilon^{-1/2} \|\tilde w(\cdot,s)\|_{H^1(D)}
\end{align*}
which then implies 
\begin{align}\label{estimate for contrast} 
    \|\tilde w(\cdot,s)\|_{H^1(D)} \preceq \frac{\epsilon^{3/2}}{\sigma\min\{1,\sigma^2\}}.
\end{align}
As we have $\tilde w=\tilde u-\tilde v$ and, $\|\tilde v\|_{H^1(D)} \preceq \epsilon^{3/2}$, the estimate \eqref{estimate for contrast} results into $\|u(\cdot,t)\|_{H^1(D)} \preceq \epsilon^{3/2}$
and similarly as before, $\|\Delta u(\cdot,t)\|_{L^2(D)} \preceq \epsilon^{3/2}$ for a.e $t\in(0,T)$. After scaling to $\partial B$, using trace theorem there and then scaling back to $\partial D$, we obtain 
\begin{align*}
    \|\partial_\nu u(\cdot,t)\|_{L^2(\partial D)} \preceq \left( \epsilon^{1/2} \|\Delta u(\cdot,t)\|_{L^2(D)} +  \epsilon^{- 1/2} \|\nabla u(\cdot,t)\|_{L^2(D)} \right) \preceq \epsilon.
\end{align*}
For future use, the preceding estimate on $\partial_\nu u$ needs to be improved to the order $\epsilon^2$. However, the scaled trace theorem readily ensures this once we manage to establish $\|\nabla u(\cdot,t)\|_{L^2(D)}\preceq \epsilon^{5/2}$. The following step addresses this issue by illustrating this gradient estimate.
\medskip

\textbf{Step III:} From \eqref{diff pde}, we take note of the Lippmann-Schwinger representation  
\begin{align}
   \nonumber u(x,t) & = v(x,t) + \int_\mathbb{R} \int_{\mathbb{R}^3} G(x,t;y,\tau) \left(\frac{1}{k_0(y)} - \frac{1}{k(y)}\right)\ u_{tt}(y,\tau) \, dy \, d\tau \, \\
   \label{Lippmann-Schwinger 1} & \quad \quad - \textnormal{div} \int_\mathbb{R} \int_{\mathbb{R}^3} G(x,t;y,\tau) \left(\frac{1}{\rho_0(y)} - \frac{1}{\rho(y)}\right) \nabla u(y,\tau) \, dy \ d\tau.
\end{align}
By applying the space gradient in \eqref{Lippmann-Schwinger 1}, it reduces to the integral equation
\begin{align}\label{gradient Lippmann-Schwinger_1}
    \nabla u(x,t) - \nabla \textnormal{div} \mathcal{N}_G\left[\alpha\nabla u\right](x,t) = \nabla v(x,t) - \nabla\mathcal{N}_G\left(\beta \partial_t^2 u\right)(x,t)
\end{align}
which will be analyzed by employing the eigensystem of certain integral operators. Keeping this in mind, let us introduce some notations and define some integral operators
\begin{align}
 \label{contrast}   & \alpha(y) = \rho^{-1}(y) - \rho^{-1}_0(y), \quad \beta(y) =k^{-1}(y)- k^{-1}_0(y), \quad \gamma(y) = \beta(y) - \alpha(y) \frac{\rho(y)}{k(y)}, \\
 \label{volume op_variable}   & \mathcal{N}_{G}(f)(x,t) = \int_{\mathbb{R}^3 \times \mathbb{R}} G(x,t;y,\tau) f(y,\tau) \, dy\, d\tau, \quad (x,t)\in \mathbb{R}^n\times \mathbb{R}, \\
 \label{single layer op}   &  \mathcal{S}_Gf(x,t) = \int_{\partial D \times \mathbb{R}} G(x,t;y,\tau) f(y,\tau) \, dS_y \, d\tau, \quad (x,t)\in \partial D \times \mathbb{R}, \\
 \label{magnetization}   & \nabla M_\epsilon[g](x) = \nabla \int_{D} \nabla_y \left(\frac{1}{4\pi|x-y|}\right) \cdot g(y) \, dy, \quad x\in D.
\end{align}

Here, we do not intend to provide explicit mapping properties and other exotic features of these integral operators. Our discussion is limited to the current set-up where the functions are sufficiently regular. With a slight abuse of notation, let us denote $\nabla M_0$ as the integral operator which can be defined in identical terms to \eqref{magnetization} but in $B$. Let us note an important observation which showcases the scale invariance of the eigenvalues of $\nabla M_\epsilon$. This can be proved by a simple scaling argument. Indeed, if $\tilde e_n$ denotes an normalized eigenfunction corresponding to an eigenvalue $\tilde\lambda_n$ for $\nabla M_0$, then the function 
\begin{align*}
    e_n(x) = \frac{1}{\epsilon^{3/2}} \tilde e_n\left(\frac{x-z}{\epsilon}\right), \quad x\in D
\end{align*}
is normalized in $L^2(D)$ and defines an eigenfunction of $\nabla M_\epsilon$ for the same eigenvalue i.e. $\tilde\lambda_n$. We recall the following orthogonal decomposition of $\mathbb{L}^2(D)$ which is the space of square integrable vector functions in $D$,
\begin{align*}
    \mathbb{L}^2(D) = \mathbb{H}_0(\textnormal{div}, 0) \, \oplus \mathbb{H}_0(\textnormal{curl}, 0) \, \oplus \nabla\mathbb{H}_{\textnormal{arm}}
\end{align*}
where $\mathbb{H}_0(\textnormal{div},\, 0)$ is the $L^2$-closure of the compactly supported divergence-free smooth functions in $D$ and
\begin{align*}
    \mathbb{H}_0(\textnormal{curl},\, 0) =   \{\nabla \phi, \, \phi\in H^1_0(D)\},  \quad
    \nabla\mathbb{H}_{\textnormal{arm}} =   \{\nabla \psi, \, \psi\in H^1(D) \textnormal{ and } \Delta \psi = 0 \textnormal{ in } D\}.  
\end{align*}
The operator $\nabla M_\epsilon$ acts as  zero and identity maps on the subspace $\mathbb{H}_0(\textnormal{div}\,0)$ and $\mathbb{H}_0(\textnormal{curl}, 0)$ respectively. In addition, it has a sequences of eigenvelues $(\lambda_n^{(3)})_{n\in \mathbb{N}}$ for which the corresponding eigenfunctions  $(e_n^{(3)})_{n\in \mathbb{N}}$ span $\nabla\mathbb{H}_{\textnormal{arm}}$, see \cite{Raevskii} for instance. In view of spectral decomposition of $\mathbb{L}^2(D)$ for the operator $\nabla M_\epsilon$, let us take any orthonormal basis $(e^{(1)}_n)_{n\in \mathbb{N}}$ of $\mathbb{H}_0(\textnormal{div}\,0)$ and any orthonormal basis $(e^{(1)}_n)_{n\in \mathbb{N}}$ of $\mathbb{H}_0(\textnormal{curl}, 0)$. Therefore $\{0, e^{(1)}_n\}_{n=1}^{\infty}$, $\{1, e^{(3)}_n\}_{n=1}^{\infty}$ and $\{\lambda^{(3)}_n, e^{(3)}_n\}_{n=1}^{\infty}$ denote normalized eigensystems which form an orthonormal basis for $\mathbb{L}^2(D)$, when considered together. Additionally, we have $\lambda^{(3)}_n \to 1$ as $n\to \infty$.
The integral equation \eqref{gradient Lippmann-Schwinger_1} reduces to   
\begin{align}\label{gradient Lippmann-Schwinger_2}
    \nabla u(x,t) + \alpha(z)\rho_0(z) \nabla M_\epsilon\left[\nabla u\right](x,t) = \nabla v(x,t) - \nabla\mathcal{N}_G\left(\beta \partial_t^2 u\right)(x,t) + \textnormal{err.}(x,t)
\end{align}
To make the presentation less cumbersome, we avoid expressing the term `err.' in its original form. We will be explicit about it in a subsequent step. For the time being, we continue the discussion regarding it an expression consisting of lower order terms which comes up in $\nabla \textnormal{div} \mathcal{N}_G\left[\alpha\nabla u\right](x,t)$.

\medskip

Making use of the orthonormality of the basis, we write
\begin{align}\label{decomposition w.r.t eigenfuncs}
    \|\nabla u\|^2_{L^2(D)} = \sum_{n=1}^{\infty} \langle\nabla u, e^{(1)}_n\rangle^2 +  \sum_{n=1}^{\infty} \langle \nabla u, e^{(2)}_n\rangle^2 +  \sum_{n=1}^{\infty} \langle \nabla u, e^{(3)}_n\rangle^2.
\end{align}
It is quite straightforward to notice that $\langle\nabla u, e^{(1)}_n\rangle = 0,\, \forall n\in \mathbb{N}$. Therefore, we focus on the remaining terms of \eqref{decomposition w.r.t eigenfuncs}. From \eqref{gradient Lippmann-Schwinger_2}, we see 
\begin{align*}
    \langle \nabla u, e^{(2)}_n\rangle = \frac{1}{1+\alpha(z)\rho_0(z)} \left(\langle\nabla v, e^{(2)}_n\rangle - \langle\nabla \mathcal{N}_G\left(\beta \partial_t^2 u\right),e^{(2)}_n\rangle + \langle \textnormal{err.},e^{(2)}_n\rangle \right), \quad n\in\mathbb{N}, \\
    \langle \nabla u, e^{(3)}_n\rangle = \frac{1}{1+\lambda^{(3)}_n\alpha(z)\rho_0(z)} \left(\langle\nabla v, e^{(3)}_n\rangle - \langle\nabla \mathcal{N}_G\left(\beta\partial_t^2 u\right),e^{(3)}_n\rangle + \langle \textnormal{err.},e^{(3)}_n\rangle \right), \quad n\in\mathbb{N}.
\end{align*}
At this point, we recall the scale invariance property of the eigenvalues of $\nabla M_\epsilon$. Now, the initial wave-field $v(\cdot,\cdot)$ being smooth in space variables and the volume of $D$ being of order $\epsilon^3$, we arrive at the estimation
\begin{align*}
   \frac{1}{\left( 1+\alpha(z)\rho_0(z)\right)^2} \sum_{n=1}^\infty \left\langle\nabla v, e^{(2)}_n\right\rangle^2 + \sum_{n=1}^\infty \frac{1}{\left( 1+\lambda_n^{(3)}\alpha(z)\rho_0(z)\right)^2} \left\langle\nabla v, e^{(3)}_n\right\rangle^2 \, \preceq \epsilon^4 \, \|\nabla v\|^2_{L^2(D)} \, \preceq \epsilon^{7},  
\end{align*}
where we have used the fact that $\alpha(z)\rho_0(z) \simeq \epsilon^{-2}$ and $\lambda_n \to 1$ as $n\to\infty$. It is also apparent from the preceding argument that
\begin{align*}
   & \frac{1}{\left( 1+\alpha(z)\rho_0(z)\right)^2} \sum_{n=1}^\infty \left\langle\nabla\mathcal{N}_G\left(\beta \partial_t^2 u\right), e^{(2)}_n\right\rangle^2 + \sum_{n=1}^\infty \frac{1}{\left( 1+\lambda_n^{(3)}\alpha(z)\rho_0(z)\right)^2}\left \langle\nabla\mathcal{N}_G\left(\beta\partial_t^2 u\right), e^{(3)}_n\right\rangle^2 \\
   & \preceq \epsilon^4 \|\mathcal{N}_G\left(\beta \partial_t^2 u\right)\|^2_{H^1(D)}.
\end{align*}
On this point, we assert that $\|\mathcal{N}_G\left[\beta\partial_t^2 u\right]\|_{H^1(D)} \preceq \epsilon^{1/2}$. In order to support this assertion, let us denote $p(x,t) := \mathcal{N}_G\left[\beta \partial_t^2 u\right](x,t)$. Then, $p$ solves the IVP
\begin{align*}
   \begin{cases}
    k_0^{-1}(x) \partial_t ^2 p - \textnormal{div}\left(\rho^{-1}_0(x) \nabla p\right) = \beta(x) \partial_t^2 u(x,t), \, \textnormal{ in } \mathbb{R}^3\times\mathbb{R}, \\
    p(\cdot,0) = \partial_t p(\cdot,0) = 0, \, \textnormal{ in } \mathbb{R}^3.
  \end{cases}   
\end{align*}
Based on the arguments presented in Step I, we observe that
\begin{align*}
    \|p(\cdot,t)\|_{H^2(\Omega)} & \preceq \|\beta(\cdot)\partial_t^2 u(\cdot,t)\|_{L^2(D)}  \preceq \epsilon^{-2} \|\partial_t^2 u(\cdot,t)\|_{L^2(D)} \preceq \epsilon^{-2} \cdot \epsilon^{3/2} = \epsilon^{-1/2}.
\end{align*}
By virtue of Sobolev embedding theorem, we conclude that $\nabla p(\cdot,t) \in L^6(\Omega)$ satisfying 
\begin{align}\label{grad estimate of p_1}
    \|\nabla p(\cdot,t)\|_{L^6(\Omega)} \preceq \|p(\cdot,t)\|_{H^2(\Omega)} \preceq \epsilon^{-1/2}.
\end{align}
Now, we employ H\"older's inequality to obtain
\begin{align*}
    \|\nabla p(\cdot,t)\|^2_{L^2(D)} = \int_{D} |\nabla p(x,t)|^2 \, dx & \le \left(\int_{D} |\nabla p(x,t)|^6 \, dx \right)^{1/3} \left(\int_{D}dx\right)^{2/3}  = \|\nabla p(\cdot,t)\|^{2}_{L^6(D)} |D|^{2/3}
\end{align*}
which, after using the estimate \eqref{grad estimate of p_1}, translates to
\begin{align*}
    \|\nabla p(\cdot,t)\|_{L^2(D)} \preceq \|\nabla p(\cdot,t)\|_{L^6(\Omega)} |D|^{1/3} \preceq \epsilon^{-1/2} \cdot \epsilon = \epsilon^{1/2}.
\end{align*}
Therefore, the assertion is verified and as a consequence, we obtain
\begin{align*}
    \frac{1}{\left( 1+\alpha(z)\rho_0(z)\right)^2} \sum_{n=1}^\infty \left\langle\nabla\mathcal{N}_G\left(\beta \partial_t^2 u\right), e^{(2)}_n\right\rangle^2 + \sum_{n=1}^\infty \frac{1}{\left( 1+\lambda_n^{(3)}\alpha(z)\rho_0(z)\right)^2} \left\langle\nabla\mathcal{N}_G\left(\beta \partial_t^2 u\right), e^{(3)}_n\right\rangle^2 \preceq \epsilon^5.
\end{align*}

\smallskip

To conclude this section, it remains to estimate the error term appearing in \eqref{gradient Lippmann-Schwinger_2}. 
We now elaborate on the structure of the error term. Let us introduce the integral operator 
\begin{align*}
    \mathcal{N}_{G_z}(f)(x,t) = \int_{\mathbb{R}^3 \times \mathbb{R}} G_z(x,t;y,\tau) f(y,\tau) \, dy\, d\tau, \quad (x,t)\in \mathbb{R}^n\times \mathbb{R}, 
\end{align*}
with $G_z$ denoting the fundamental solution of $k_0^{-1}(z)\partial^2_{t} - \rho_0^{-1}(z)\Delta_x$, that is
\begin{align}\label{G_z_expression}
    G_z(x,t;y,\tau) = \rho_0(z) \frac{\delta(t-\tau-c_0^{-1}(z)|x-y|)}{4\pi |x-y|}.
\end{align}
Setting 
\begin{align}\label{defn_m_n}
    m(x,t) = \mathcal{N}_G\left[\alpha\nabla u\right](x,t), \quad n(x,t) = \alpha(z) \mathcal{N}_{G_z}\left[\chi_{D\epsilon}\nabla u\right](x,t),
\end{align}
and we first use Taylor's theorem to see
\begin{align*}
    \nabla u(y,t)-\nabla u(y,t-c_0^{-1}(z)|x-y|) & = c_0^{-1}(z) |x-y| \, \partial_t\nabla u(y,t) \\
    & \quad - c_0^{-1}(z) |x-y| \int_{0}^{1} \partial_t^2 \nabla u(y,t-c_0^{-1}(z) |x-y|\tau)\, d\tau
\end{align*}
and then we proceed to calculate the error term in \eqref{gradient Lippmann-Schwinger_2} as follows
\begin{align} \label{err_est}
    \nonumber \textnormal{err.}(x,t)  & = {\rho_0(z)}{\alpha(z)}\nabla M_\epsilon\left[\nabla u\right](x,t) - \nabla \textnormal{div}\mathcal{N}_G\left[\alpha\nabla u\right](x,t) \\
   \nonumber & = {\rho_0(z)}{\alpha(z)}\nabla \textnormal{div} \int_{D} \frac{\nabla u(y,t)-\nabla u(y,t-c_0^{-1}(z)|x-y|)}{4\pi |x-y|}\, dy + \alpha(z) \nabla\textnormal{div}\mathcal{N}_{G_z}\left[ \chi_{D}\nabla u\right](x,t) \\
  \nonumber & \quad \quad - \nabla \textnormal{div}\mathcal{N}_G\left[\alpha \nabla u\right](x,t)\\
   \nonumber & = - \frac{\rho_0(z)\alpha(z)}{8\pi c_0^2(z)} \nabla \textnormal{div} \int_{D}\int_0^1 |x-y|\,\tau \,\partial^2_t \nabla u(y,t-c_0^{-1}(z) |x-y|+c_0^{-1}(z) |x-y|\tau)\, d\tau\, dy \\
    & \quad - \nabla\textnormal{div}(m-n)(x,t).
\end{align}
To derive \eqref{err_est}, we have used some observations which we discuss next. In view of \eqref{volume op_variable} and \eqref{G_z_expression}, we first added and subtracted the following quantity in the term  $\textnormal{err.}(x,t)$, 
\begin{align*}
    \alpha(z)\nabla \textnormal{div}\mathcal{N}_{G_z}\left[\chi_D\nabla u\right] = {\rho_0(z)}{\alpha(z)}\nabla \textnormal{div} \int_{D} \frac{\nabla u(y,t-c_0^{-1}(z)|x-y|)}{4\pi |x-y|}\, dy, 
\end{align*}
and then used an integral form of Taylor's expansion to see
\begin{align}\label{simplification_err_est}
  \nonumber & \frac{\nabla u(y,t)-\nabla u(y,t-c_0^{-1}(z)|x-y|)}{4\pi |x-y|} \\
  \nonumber  & = \frac{1}{4\pi|x-y|} \left[c^{-1}_0(z)|x-y| \,\partial_t u(y,t) - \frac{|x-y|^2}{2c_0^2(z)}\int_0^1 \,\partial^2_t \nabla u(y,t-c_0^{-1}(z) |x-y|+c_0^{-1}(z) |x-y|\tau)\, d\tau \right] \\
    & = \frac{1}{4c_0(z)} \partial_t u(y,t) - \frac{|x-y|}{8\pi c_0^2(z)}\int_0^1 \tau \,\partial^2_t \nabla u(y,t-c_0^{-1}(z) |x-y|+c_0^{-1}(z) |x-y|\tau)\, d\tau.
\end{align}
Note here that, the first term in \eqref{simplification_err_est} independent of $x$-variable and thus does not survive when we apply $\nabla\textnormal{div}$ as considered in the derivation of \eqref{err_est}.

Let us first focus on the first term of $\textnormal{err.}(x,t)$ which we denote as $\mathcal{R}(x,t)$ i.e. 
\begin{align*}
    \mathcal{R}(x,t) = \frac{\rho_0(z)\alpha(z)}{8\pi c_0^2(z)} \nabla \textnormal{div} \int_{D}\int_0^1 |x-y|\,\tau \,\partial^2_t \nabla u(y,t-c_0^{-1}(z) |x-y|+c_0^{-1}(z) |x-y|\tau)\, d\tau\, dy.
\end{align*}
From a straightforward differentiation, we arrive at the estimation 
\begin{align*}
    \left\vert \mathcal{R}(x,t)\right\vert \preceq \alpha(z) \sup_{0\le \tau \le t}\int_{D} \frac{1}{|x-y|} \sum_{k=1}^4 \lvert\partial_t^k \nabla u\rvert(y,\tau) \, dy.
\end{align*}
At this point, we use the generalized Young's inequality \cite[Lemma 0.10]{Folland_PDE_book} to observe
\begin{align*}
    \|\mathcal{R}(\cdot,t)\|_{L^2(D)} & \preceq \alpha(z)\,\sup_{x\in D}\int_{D} \frac{dy}{|x-y|} \,\sup_{0\le\tau\le t}\sum_{k=1}^4 \|\partial_t^k \nabla u(\cdot,\tau)\|_{L^2(D)} 
\end{align*}
which, due to a scaling argument and Step II in the proof of Theorem \ref{gradient estimate}, reduces to 
\begin{align} \label{err_est_1}
    \|\mathcal{R}(\cdot,t)\|_{L^2(D)} \preceq \epsilon^{-2+2+\frac{3}{2}} = \epsilon^{3/2}.
\end{align}
Next, we show that $\|\nabla\textnormal{div}(m-n)(\cdot,t)\|_{L^2(D)} \preceq \epsilon^{1/2}$, where the functions $m(x,t)$ and $n(x,t)$ are introduced in \eqref{defn_m_n}. In pursuit of doing so, we first define
\begin{align*}
    n_*(x,t) = \mathcal{N}_{G_z}\left[\alpha \nabla u\right](x,t) = n(x,t) \, + \mathcal{N}_{G_z}\left[\left(\alpha(\cdot)-\alpha(z)\right) \nabla u\right](x,t).
\end{align*}
Following the arguments presented earlier, we see that
\begin{align*}
    \|\mathcal{N}_{G_z}\left[\left(\alpha(\cdot)-\alpha(z)\right) \nabla u(\cdot,t)\right]\|_{H^2(D)} \preceq \|\left(\alpha(\cdot)-\alpha(z)\right)\nabla u(\cdot,t)\|_{L^2(D)} \preceq \epsilon \|\nabla u(\cdot,t)\|_{L^2(D)} \preceq \epsilon^{5/2},
\end{align*}
and as a consequence
\begin{align}\label{comparing H^2 estimate}
    \|\nabla\textnormal{div}(m-n)(\cdot,t)\|_{L^2(D)} \preceq \|\nabla\textnormal{div}(m-n_*)(\cdot,t)\|_{L^2(D)} + \epsilon^{5/2}.
\end{align}
To deal with the estimation of $\|\nabla\textnormal{div}(m-n_*)(\cdot,t)\|_{L^2(D)}$, we rely on the following lemma which states 
\begin{lemma}\label{H^2 estimate}
    For $u$ solving \eqref{IVP_u}, let $m(x,t) = \mathcal{N}_G\left[\alpha\nabla u\right](x,t)$, and $n_*(x,t) = \mathcal{N}_{G_z}\left[\alpha\nabla u\right](x,t)$. Then, 
    \begin{align*}
        \|\nabla \textnormal{div}(m-n_*)(\cdot,t)\|_{L^2(D)} \preceq \epsilon^{1/2}, \quad \textnormal{ for a.e. } t \in (0,T). 
    \end{align*}
\end{lemma}
For now, we utilize Lemma \ref{H^2 estimate} and postpone its proof to Section \ref{Appendix}. Combining \eqref{comparing H^2 estimate} and Lemma \ref{H^2 estimate}, where we showed that $\|\nabla\textnormal{div}(m-n_*)(\cdot,t)\|_{L^2(D)} \preceq \epsilon^{\frac{5}{2}}$, we get $\|\nabla\textnormal{div}(m-n)(\cdot,t)\|_{L^2(D)} \preceq \sqrt\epsilon$ for $t\in(0,T)$. It implies
\begin{align*}
     \frac{1}{\left( 1+\alpha(z)\rho_0(z)\right)^2} \sum_{n=1}^\infty \left\langle\textnormal{err}., \,e^{(2)}_n\right\rangle^2 + \sum_{n=1}^\infty \frac{1}{\left( 1+\lambda_n^{(3)}\alpha(z)\rho_0(z)\right)^2} \left\langle\textnormal{err}., \,e^{(3)}_n\right\rangle^2 \, \preceq \epsilon^4 \,\|\textnormal{err}.\|^2_{L^2(D)} \, \preceq \epsilon^5,
\end{align*}
where the last line also uses the relation \eqref{err_est} and estimate \eqref{err_est_1}. This completes the discussion for Step III.


\section{Proof of Theorem \ref{asymptotic expansion_w}}\label{theo_expansion}
As stated in the introduction, Theorem \ref{asymptotic expansion_w} concerns with the asymptotic profile of wave-field $w$ with respect to the parameter $\epsilon$, observed at $x_0\in\partial\Omega$ for a finite time interval. To demonstrate this, we rely on the Lippmann Schwinger representation of $w$ which, using the a-priori estimate previded in Theorem \ref{gradient estimate}, identifies the dominating part to be the average of Neuamnn trace of $u$ on $\partial D$. However, this quantity can be explicitly obtained by solving a second order ODE which arises in consideration of a boundary integral equation; see \eqref{neumann derivative}. Our argument additionally requires some results which we showcase in Section \ref{Appendix}.

We start by reconsidering the Lippmann Schwinger representation \eqref{Lippmann-Schwinger 1}. After that, an integration by parts leads to 
\begin{align}\label{representation of w}
   \nonumber w(x,t) &  = - \mathcal{N}_{G}\left( \beta(\cdot)\,\partial^2_t u(\cdot,\cdot)\right)(x,t) + \textnormal{div}\, \mathcal{N}_{G}\left( \alpha(\cdot) \, \nabla u(\cdot,\cdot)\right)(x,t) \\
   \nonumber & = - \int_{D \times \mathbb{R}} G(x,t;y,\tau) \, \left[\beta(y) \partial^2_t u(y,\tau) - \textnormal{div}(\alpha(y) \nabla u(y,\tau))\right]\, dy \, d\tau \\
    & \quad \quad -  \int_{\partial D \times \mathbb{R}} G(x,t;y,\tau) \alpha(y) \, \partial_\nu u(y,\tau) \, dS_y \, d\tau ,
\end{align}
which, in view of \eqref{IVP_u}, changes to
\begin{align*}
     & = - \int_{D \times \mathbb{R}} G(x,t;y,\tau) \left[\left( 1 - \frac{k_1}{k_0(y)} \right) J(y,\tau) + \left( \rho^{-1}_0 - k_0^{-1} \frac{k_1}{\rho_1} \right)\Delta u(y,\tau) \right]\, dy \,d\tau \\
    & \, - \int_{D \times\mathbb{R}} G(x,t;y,\tau) \nabla(\rho^{-1}_0)\cdot\nabla u(y,\tau) \, dy\,d\tau  + \int_{\partial D \times \mathbb{R}} G(x,t;y,\tau) \alpha(y) \partial_\nu u(y,\tau) \, dS_y \, d\tau. 
\end{align*}
Therefore, we notice 
\begin{align}\label{Lippmann Schwinger_2}
    u(x,t) = v(x,t) - \mathcal{N}_G\left[k \beta J + \gamma \frac{k_1}{\rho_1} \Delta u\right](x,t) -  \mathcal{N}_G[ \chi_{D}\nabla\rho^{-1}_0\cdot\nabla u](x,t) - \mathcal{S}_G\left[\alpha \partial_\nu u\right](x,t).
\end{align}

Now, we treat the revised Lippmann Schwinger equation \eqref{Lippmann Schwinger_2} for points away from $D$, specifically at $x_0\in\partial\Omega$. The structure of $G$ in Theorem \ref{Green function} depicts that there are no singularity present in the kernel of operators in \eqref{Lippmann Schwinger_2} when we are on $\partial\Omega$. That being said, we sort the terms in \eqref{Lippmann Schwinger_2} depending on their order in pointwise sense. With the regularity of $J(\cdot,\cdot)$ and support condition of $g$ in time, we deduce
\begin{align*}
   \big\vert \mathcal{N}_G[k \beta J](x_0,t) \big\vert \preceq & \left\vert \int_{\mathbb{R}^3} \frac{m(x_0,y)}{4\pi\zeta(x_0,y)} k \beta(y) J(y, t -\zeta(x_0,y)) dy \right\vert 
    + \left\vert \int_{\mathbb{R}^3 \times \mathbb{R}} g(x_0,t-\tau;y)\, k \beta(y) J(y,\tau) \, dy d\tau\right\vert \preceq \epsilon^{3}.
\end{align*}
It is worth mentioning that the space integration in the preceding relation is limited to $D$ which is due to the support condition of $\beta$. With Cauchy-Schwarz inequality and the relation \eqref{a-priori estimate} from Theorem \ref{gradient estimate}, we notice
\begin{align*}
    \left\vert \mathcal{N}_G \left[\gamma \frac{k}{\rho} \Delta u\right](x_0,t)\right\vert 
   \preceq \sup_{0\le \tau\le t} \left|\int_{D} \Delta u(y,\tau) \, dy \right| \preceq |D|^{1/2} \sup_{0\le \tau\le t} \|\Delta u(\cdot,\tau)\|_{L^2(D)} \preceq \epsilon^{3},
\end{align*}
and similarly,
\begin{align*}
   \big\vert \mathcal{N}_G[\chi_{D}\nabla\rho_0\cdot\nabla u](x_0,t) \big\vert \preceq \epsilon^{3}.
\end{align*}

\smallskip

In \eqref{Lippmann Schwinger_2}, we are now left with the term $\mathcal{S}_G[\alpha \partial_\nu u](x_0,t)$. This is the term which supplies the dominant terms in the expansion \eqref{asymptotic expansion_2}. From the structure of $G$ in Theorem \ref{Green function}, we expand
\begin{align*}
    & \mathcal{S}_G [\alpha \partial_\nu u](x_0,t)  = \int_{\partial D}  \frac{\alpha(y)\, m(x_0,y)}{4\pi\zeta(x_0,y)}  \partial_\nu u\left(y,t-\zeta(x_0,y)\right) dS_y  + \int_{\partial D \times \mathbb{R}} g(x_0,t-\tau;y) \alpha(y)\,\partial_\nu u(y,\tau)  \\
    & = \frac{\alpha(z)m(x_0,z)}{4\pi \zeta(x_0,z)} \int_{\partial D} \partial_\nu u\left(y,t_z\right) dS_y \, + \int_{\partial D} \left[ \frac{\alpha(y)m(x_0,y)}{4\pi \zeta(x_0,y)}-\frac{\alpha(z)m(x_0,z)}{4\pi \zeta(x_0,z)}\right]\,\partial_\nu u\left(y,t_y\right)\, dS_y \\
    &   \quad + \alpha(z) \int_{\partial D} \int_0^{t_z} g(x_0,t-\tau;z)\partial_\nu u(y,\tau) \, d\tau dS_y + \frac{\alpha(z)m(x_0,z)}{4\pi \zeta(x_0,z)} \, \int_{\partial D} \int^{t_z}_{t_y}\partial_{t} \partial_\nu u(y,\tau) \, d\tau dS_y  \\
    & \quad + \int_{\partial D \times \mathbb{R}} \left( \alpha(y)g(x_0,t-\tau;y) - \alpha(z)g(x_0,t-\tau;z) \right) \partial_\nu u(y,\tau) \, dS_y \,d\tau,
\end{align*}
where $t_y= t -\zeta(x_0,y),$ for $y\in\overline D$. With the aid of Cauchy-Schwarz inequality and \eqref{a-priori estimate} from Theorem \ref{gradient estimate}, we now discard lower order terms appearing in the above expansion. We notice
\begin{align*}
    & \left\vert \alpha(z) \int_{\partial D} \int^{t_y}_{t_z}\partial_{t} \partial_\nu u(y,\tau) \, d\tau dS_y\right\vert \preceq \epsilon^{-2} \left\vert \partial D\right\vert^{1/2} \sup_{0\le \tau \le t}\| \left(\zeta(x_0, \cdot)- \zeta(x_0,z)\right)\partial_t \partial_\nu u (\cdot,\tau)\|_{L^2(\partial D)} \preceq \epsilon^2 ,
\end{align*}
and, similarly
\begin{align*}
    \left\vert\int_{\partial D} \left( \frac{\alpha(y)m(x_0,y)}{4\pi \zeta(x_0,y)}-\frac{\alpha(z)m(x_0,z)}{4\pi \zeta(x_0,z)}\right)\,\partial_\nu u\left(y,t_y\right)\, dS_y \right\vert & \preceq \epsilon^{-2} \left\vert \partial D\right\vert^{1/2} \epsilon \, \sup_{0\le \tau \le t}\|\partial_\nu u (\cdot,\tau)\|_{L^2(\partial D)} \preceq \epsilon^{2}, \\
     \left\vert \int_{\partial D \times \mathbb{R}} \left( \alpha(y)g(x_0,t-\tau;y) - \alpha(z)g(x_0,t-\tau;z) \right) \partial_\nu u(y,\tau) dS_y d\tau \right\vert & \preceq \epsilon^{-2} \left\vert \partial D\right\vert^{1/2} \epsilon \sup_{0\le \tau \le t}\|\partial_\nu u (\cdot,\tau)\|_{L^2(\partial D)} \preceq \epsilon^{2}.
\end{align*}
 In the last line, we have used the Lipschitz regularity of $g(x,t-\tau;y)$ with respect to the $y$ variable inside the characteristic conoid $t>\tau+\zeta(x,y)$. This consideration requires enough regularity for the bulk and density which is already incorporated in Theorem \ref{Green function}. Consequently, we write
\begin{align*}
     \mathcal{S}_G[\alpha \partial_\nu u](x_0,t) = \frac{\alpha(z)m(x_0,z)}{4\pi \zeta(x_0,z)} \int_{\partial D} \partial_\nu u\left(y,t_z\right) dS_y  +  \alpha(z) \int_{\partial D \times \mathbb{R}} g(x_0,t-\tau;z)\partial_\nu u(y,\tau) d\tau dS_y + O(\epsilon^2)
\end{align*}
to be understood pointwise in time. After substituting all these in the relation \eqref{Lippmann Schwinger_2}, we finally have 
\begin{align}\label{asymptotic expansion_1}
 \nonumber u(x_0,t) & = v(x_0,t) + \frac{\alpha(z)m(x_0,z)}{4\pi \zeta(x_0,z)} \int_{\partial D} \partial_\nu u\left(y,t_z\right)dS_y + \alpha(z) \int_{\partial D} \int_0^{t_z} g(x_0,t-\tau;z)\partial_\nu u(y,\tau) d\tau \,dS_y, \\
  & \quad +  O(\epsilon^2)  
\end{align}
which holds in pointwise sense in time. Indeed from straightforward calculations, we see that the terms in \eqref{asymptotic expansion_1} with Neumann averages  are of order $\epsilon$. Having obtained \eqref{asymptotic expansion_1}, it suffices to determine the average of $\partial_\nu u$ over $\partial D$. This is the aspect which we pursue next. It yields the desired  asymptotic expansion \eqref{asymptotic expansion_2}. 
\medskip

Taking the Neumann derivative to \eqref{Lippmann Schwinger_2} on $\partial D$ and then multiplying it by $\rho_0^{-1}(\cdot)$, we obtain
\begin{align}\label{neumann derivative}
   \nonumber \rho_0^{-1}(x) \partial_\nu u(x,t) + \frac{1}{2} \alpha(x)\, \partial_\nu u(x,t) & = \rho_0^{-1}(x) \partial_\nu v(x,t) - \rho_0^{-1}(x)\partial_\nu\mathcal{N}\left[k \beta J + \gamma \frac{k}{\rho} \Delta u\right](x,t)   \\
    & \quad - \rho_0^{-1}(x)\,\partial_\nu\mathcal{N}_G[ \chi_{D}\nabla\rho^{-1}_0\cdot\nabla u](x,t) - \rho_0^{-1}(x)\,\mathcal{K}^t_G\left[\alpha\partial_\nu u\right](x,t)  ,
\end{align}
where we utilize the jump relation 
\begin{align*}
    \lim_{h\to 0-} \partial_\nu \mathcal{S}_G [f](x+h\nu(x),t) = \frac{\rho_0(x)}{2} f(x,t) + \mathcal{K}^t_G [f](x,t), \quad \textnormal{ for } x\in\partial D,\, t \in \mathbb{R}.
\end{align*}
Here, the Neumann-Poincar\'e operator $\mathcal{K}^t_G$ is defined as 
\begin{align}\label{Neumann-poincare}
    \mathcal{K}^t_G [f](x,t) := \int_{\partial D \times \mathbb{R}} \partial_{\nu_x} G(x,t;y,\tau) \, f(y,\tau) \, dS_y\, d\tau,
\end{align}   
which, by virtue of Theorem \ref{Green function}, reduces to 
\begin{align}\label{NP_expression}
  \mathcal{K}^t_G [f](x,t) & = \int_{\partial D} \left(\frac{\partial_\nu m(x,y)}{4\pi\zeta(x,y)} - \frac{m(x,y)\, \partial_\nu \zeta(x,y)}{4\pi\zeta^2(x,y)} \right)\, f(y,t-\zeta(x,y)) \, dS_y \\
     & -  \int_{\partial D} \frac{m(x,y)\, \partial_\nu \zeta(x,y)}{4\pi\zeta(x,y)}\, \partial_t f(y,t-\zeta(x,y)) \, dS_y + \int_{\partial D\times \mathbb{R}} \partial_\nu g(x,t-\tau;y)\, f(y,\tau) \, dS_y\, d\tau,
\end{align}
for $f$ sufficiently smooth in time. Let us now rewrite the equation \eqref{neumann derivative} as 
\begin{align}
    \nonumber \frac{1}{2} \left(\rho_0^{-1}(x) + \rho_1^{-1}\right)\partial_\nu u(x,t) + \underbrace{\rho_0^{-1}(x) \,\mathcal{K}^t_G\left[\alpha\partial_\nu u\right](x,t)}_{T_1} & = \rho_0^{-1}(x) \partial_\nu v(x,t) - \underbrace{\rho_0^{-1}(x)\,\partial_\nu\mathcal{N}_G[ \chi_{D}\nabla\rho^{-1}_0\cdot\nabla u](x,t)}_{T_2} \\
   \label{Neumann Lippman-Schwinger} &  - \underbrace{\rho_0^{-1}(x)\partial_\nu\mathcal{N}\left[k \beta J + \gamma \frac{k_1}{\rho_1} \Delta u\right](x,t)}_{T_3}
\end{align}
which we intend to convert into an ODE for the average of Neumann derivative of $u$ by performing error analysis to each of the terms in \eqref{Neumann Lippman-Schwinger}. For the ease of presentation, we divide the discussion in two steps. In Step I, we deal with the boundary integral of $T_2$ and $T_3$ whereas in Step II, we discuss the same for $T_1$.

\medskip

\textbf{Step I:} Our objective here is to establish that the averages of $T_2$ and $T_3$ are negligible as they are of order $\epsilon^4$ (pointwise in time). Here, we mostly rely on the structure of $G$ in Theorem \ref{Green function} and the estimates from Theorem \ref{gradient estimate}. Let us first deal with the term $T_2$. An integration by parts gives us
\begin{align*}
    & \int_{\partial D} \rho_0^{-1}(x)\,\partial_\nu\mathcal{N}_G[ \chi_{D}\nabla\rho_0\cdot\nabla u](x,t) \, dS_x = \int_{\partial D} \rho_0^{-1}(x) \int_{D \times \mathbb{R}} \nabla\rho_0(y) \cdot\nabla u(y,\tau)\partial_{\nu_x} G(x,t;y,\tau) \, dy \, d\tau \, dS_{x} \\
    & \quad \quad = \int_{D \times \mathbb{R}} \nabla\rho_0^{-1}(y) \cdot\nabla u(y,\tau) \int_{D} \textnormal{div}(\rho_0^{-1}(x)\nabla_x G)(x,t;y,\tau)\,dx \,dy \,d\tau \\
    &  \quad \quad = \int_{D}  \nabla\rho_0^{-1}(x) \cdot\nabla u(x,t) \, dx - \int_{D \times \mathbb{R}} \nabla\rho^{-1}_0(y) \cdot\nabla u(y,\tau) \int_{D} k_0^{-1}(x) \,\partial_t^2 G (x,t;y,\tau)\,dx \,dy \,d\tau 
\end{align*}
 which, due to Theorem \ref{Green function}, reduces to the following 
\begin{align*}    
    & = \int_{D}  \nabla\rho_0^{-1}(x) \cdot\nabla u(x,t) \, dx -\int_{D \times D} \frac{m(x,y)}{4\pi k_0(x,y)\zeta(x,y)} \nabla\rho^{-1}_0(y) \cdot \partial^2_t \nabla u(y,t-\zeta(x,y)) \, dx\, dy \\
    & \quad - \int_{D \times D} \int_{\mathbb{R}} g(x,t-\tau;y) \nabla\rho^{-1}_0(y) \cdot \partial^2_t\nabla u(y,\tau) \, d\tau \, dx \, dy.
\end{align*}
Next we use the Cauchy-Schwarz inequality in combination with \eqref{a-priori estimate} to observe 
\begin{align*}
    \left\vert \int_{D}  \nabla\rho_0^{-1}(x) \cdot\nabla u(x,t) \, dx \right\vert \preceq \|\nabla\rho_0^{-1}(\cdot)\|_{L^2(D)}  \|\nabla u(\cdot,t)\|_{L^2(D)} \preceq \epsilon^{3/2} \cdot \epsilon^{5/2} = \epsilon^{4}.
\end{align*}
Following similar steps, we notice 
\begin{align*}
     \left\vert \int_{D \times D} \frac{m(x,y)  \nabla\rho^{-1}_0(y)}{4\pi k_0(x,y)\zeta(x,y)} \cdot \partial^2_t \nabla u(y,t-\zeta(x,y))  dx dy \right\vert & \preceq \left\vert \int_{D \times D} \frac{dx  dy}{|\zeta(x,y)|^2} \right\vert^{1/2}  \sup_{0\le\tau\le t} \left\vert \int_{D \times D} |\partial^2_t \nabla u(y,\tau|^2 dx  dy\right\vert^{1/2} \\
   & \preceq \epsilon^2\, |D|^{1/2} \sup_{0\le\tau\le t} \left\vert\int_{D} |\partial^2_t \nabla u(y,\tau)|^2\, dy\right\vert^{1/2} \preceq \epsilon^6. 
\end{align*}
 Here, we utilize the estimate for $L^2(D)$ norm of $\partial_t^2 \nabla u(\cdot,t)$ and its supremum in time which have been derived in Theorem \ref{gradient estimate}, see \eqref{a-priori estimate}. Likewise, we deduce
\begin{align*}   
     \left\vert \int_{D \times D} \int_{\mathbb{R}} g(x,t-\tau;y) \nabla\rho^{-1}_0(y) \cdot \partial^2_t\nabla u(y,\tau) \, d\tau \, dx \, dy \right\vert 
    & \preceq |D| \cdot |D|^{1/2} \sup_{0\le\tau\le t} \left\vert\int_{D} |\partial^2_t \nabla u(y,\tau)|^2\, dy\right\vert^{1/2} \preceq \epsilon^7.
\end{align*}
To summarize, we have
\begin{align}\label{T_2}
    \int_{\partial D} \rho_0^{-1}(x)\,\partial_\nu\mathcal{N}_G[ \chi_{D}\nabla\rho_0\cdot\nabla u](x,t) \, dS_x = O(\epsilon^4)
\end{align}
in pointwise sense in time. Following the same strategy as above, we obtain estimates for the boundary integral for $T_3$. But, to identify the constituting terms there, we integrate by parts and obtain
\begin{align}
   \nonumber & \int_{\partial D}\rho_0^{-1}(x)\partial_{\nu_x}\mathcal{N}\left[k_1 \beta J + \gamma \frac{k_1}{\rho_1} \Delta u\right](x,t) \, dS_x = \int_{\partial D} \rho_0^{-1}(x) \int_{D \times \mathbb{R}} \left[k_1 \beta J + \gamma \frac{k}{\rho} \Delta u\right](y,\tau)\, \partial_{\nu_x} G \, dy \, d\tau \, dS_{x} \\
   \nonumber & = \int_{D \times \mathbb{R}} \left[k_1 \beta J + \gamma \frac{k_1}{\rho_1} \Delta u\right](y,\tau) \int_{D} \textnormal{div}\left(\rho_0^{-1}(x) \nabla G\right)(x,t;y,\tau) \, dx \, dy \, d\tau \\
   \nonumber & = - \int_{D} \left[k_1 \beta J + \gamma \frac{k_1}{\rho_1} \Delta u\right](x,t) \, dx + \int_{D \times \mathbb{R}} \left[k_1 \beta J + \gamma \frac{k_1}{\rho_1} \Delta u\right](y,\tau) \int_{D} k_0^{-1}(x) \partial_t^2 G(x,t;y,\tau) \, dx \, dy \, d\tau,
\end{align}
which, after utilizing Theorem \ref{Green function}, reduces to 
\begin{align}   
   \nonumber & = - \int_{D} \left[k_1 \beta J + \gamma \frac{k_1}{\rho_1} \Delta u\right](x,t) \, dx + \int_{D \times D} \frac{m(x,y)}{4\pi k_0(x) \zeta(x,y)} \left[k_1 \beta J_{tt} + \gamma \frac{k_1}{\rho_1} \partial_t^2 \Delta u\right](y,t-\zeta(x,y)) \, dx \, dy \\
   \nonumber & \quad + \int_{D \times \mathbb{R}} \int_{D} g(x,t-\tau;y)\left[k_1 \beta J + \gamma \frac{k_1}{\rho_1} \Delta u\right](y,\tau)\, dx \, dy \, d\tau \\
   \label{T_3} & = - \int_{D} J(x,t)\, dx \, - \gamma(z)\frac{k_1}{\rho_1} \int_{\partial D} \partial_\nu u(x,t) \, dx + O(\epsilon^4).
\end{align}
The last line follows from the fact, $[k_1\beta](x) = 1 + O(\epsilon^2)$, in $D$ and the estimates such as 
\begin{align*}
   \left\vert \int_D \left[\gamma \frac{k_1}{\rho_1} \Delta u\right](x,t) \, dx - \gamma(z)\frac{k_1}{\rho_1} \int_{\partial D} \partial_\nu u (x,t)\, dS_x\right\vert & = \frac{k_1}{\rho_1} \left\vert\int_D \left( \gamma(x) - \gamma(z)\right)\Delta u(x,t) \, dx\right\vert \\
    & \preceq \epsilon\, |D|^{1/2} \, \|\Delta u(\cdot,t)\|_{L^2(D)} \preceq \epsilon^4, \\
     \left\vert \int_{D \times D} \frac{m(x,y)}{4\pi k_0(x) \zeta(x,y)} \left[k_1 \beta J_{tt}\right](y,t-\zeta(x,y)) \, dx \, dy \right\vert 
   & \preceq \left(\int_{D \times D} \frac{dx\, dy}{\zeta^2(x,y)}\right)^{1/2} |D|^{1/2} \sup_{0\le \tau \le t}\|\partial^2_t J(\cdot,\tau)\|_{L^2(D)} \\
   & \preceq \epsilon^{3} \, \left(\int_{D \times D} \frac{dx\, dy}{|x-y|^2}\right)^{1/2} \preceq \epsilon^5,
\end{align*}
and, likewise
\begin{align*}
    \left\vert \frac{k_1}{\rho_1} \int_{D \times D} \frac{m(x,y)}{4\pi k_0(x) \zeta(x,y)}   \left[\gamma\partial_t^2 \Delta u\right](y,t-\zeta(x,y)) \, dx \, dy  \right\vert & \preceq \epsilon^5, \\
   \left\vert \int_{D \times \mathbb{R}} \int_{D} g(x,t-\tau;y)\left[k_1 \beta J + \gamma \frac{k_1}{\rho_1} \Delta u\right](y,\tau)\, dx \, dy \, d\tau \right\vert & \preceq \epsilon^5.
\end{align*}
This concludes the discussion on the averages of $T_2$ and $T_3$. 

\vspace{2mm} 

\textbf{Step II :} Now, we focus on deriving the average of $T_1$ which is of central importance in the context of ODE which will be derived shortly. We will notice there that it constitutes of important terms having order lower than $\epsilon^4$ (point-wise in time). In view of integration by parts, we deduce
\begin{align*}
     \int_{\partial D} \rho_0^{-1}(x) \partial_{\nu_x}\mathcal{S}_G\left[\alpha\partial_{\nu} u\right](x,t) \, dS_x 
    & = \int_{\partial D \times \mathbb{R}} \alpha(y) \partial_{\nu_y} u(y,\tau) \int_{\partial D} \rho_0^{-1}(x) \partial_{\nu_x} G(x,t;y,\tau) \, dS_{x} \, dS_y \, d\tau \\
    & = \int_{\partial D \times \mathbb{R}} \alpha(y) \partial_{\nu_y} u(y,\tau) \int_{D} \textnormal{div}_x (\rho_0^{-1}(x)\nabla G)(x,t;y,\tau) \, dx \, dS_y \, d\tau \\
    & = -  \frac{1}{2} \int_{\partial D} \alpha(x) \partial_{\nu_x} u(x,t) dS_x + \int_{\partial D \times \mathbb{R}} \int_{D}\alpha(y) k_0^{-1}(x) \, \partial_t^2 G(x,t;y,\tau) dx dS_y d\tau 
\end{align*}
which, after an application of Theorem \ref{Green function}, turns into 
\begin{align*}    
    & = -  \frac{1}{2} \int_{\partial D} \alpha(x) \partial_{\nu_x} u(x,t) \, dS_x + \int_{D \times \partial D} \frac{\alpha(y) m(x,y)}{k_0(x) \zeta(x,y)} \partial_t^2 \partial_{\nu_y} u(y,t-\zeta(x,y)) \, dx \, dS_y \\
    & \quad \quad + \int_{\partial D \times \mathbb{R}} \int_{D} \alpha(y) g(x,t-\tau;y) \, \partial_t^2 \partial_{\nu_y} u(y,\tau) \, dx \, dS_y \, d\tau.
\end{align*}
We note the appearance of $\frac{1}{2}$ in third line of the foregoing derivation. In particular, we use here 
\begin{align}\label{delta func_identity}
    \int_{\partial D} \int_D \delta(x-y) f(y) \, dx\, dS_y = \frac{1}{2} \int_{\partial D} f(y)\, dS_y. 
\end{align}
The identity \eqref{delta func_identity} is quite standard and can be justified by using an approximate identity, for instance $\delta_n(x) := \sqrt{\frac{n}{\pi}}\exp{(-n|x|^2)}$. The fraction $\frac{1}{2}$ in \eqref{delta func_identity} is due to the normalization while integrating delta function or an approximate identity. For completeness, we briefly discuss on the identity \eqref{delta func_identity}. Because of the regularity, the boundary $\partial D$ can be straightened. Therefore it suffices to consider \eqref{delta func_identity} for the simpler case when $D= \mathbb{H}:=\{x=(x_1,x_2,x_3)\in\mathbb{R}^3; \, x_3 >0\}$. In light of this, we note that
\begin{align*}
     \int_{\partial \mathbb{H}} \int_\mathbb{H} \delta(x-y) f(y) \, dx\, dS_y & = \lim_{n\to\infty} \int_{\partial \mathbb{H}} \int_\mathbb{H} \delta_n(x-y) f(y) \, dx\, dS_y \\
     & = \lim_{n\to\infty} \int_{\partial \mathbb{H}} \left(\int_0^\infty \delta_n(x_3) \, dx_3 \,\int_\mathbb{R}\delta_n(x_1-y_1) \, dx_1 \,\int_\mathbb{R}\delta_n(x_2-y_2) \, dx_2 \right) f(y)\, dS_y \\
     & = \frac{1}{2} \int_{\partial \mathbb{H}} f(y)\, dS_y.
\end{align*}
Therefore, the identity \eqref{delta func_identity} is true.

Next, to identify the most singular terms in the integral of $T_1$, we rely on the structure of the travel time function $\zeta(\cdot,\cdot)$ in near the point $z$ and obtain 
\begin{align}\label{comparing values_2}
    \frac{\alpha(y) m(x,y)}{k_0(x) \zeta(x,y)} = \frac{1}{\rho_1 c_0^2(z)|x-y|} + \frac{O(\epsilon^{-1})}{|x-y|} + O(\epsilon^{-2}), \quad x\neq y.
\end{align}
The verification of this claim is quite straightforward and similar to the equation (3.3) in \cite{Senapati_Sini_Wang_23}.  It uses the relation \eqref{reln_travel time}, the continuity of functions such as $m(\cdot,\cdot), \, k_0(\cdot)$ and the structure of travel-time function in \eqref{reln_travel time}. A variant of this is discussed in Step I in the proof of Lemma \ref{H^2 estimate}. Furthermore, we express
\begin{align*}
    \partial_t^2 \partial_\nu u(y,t-\zeta(x,y)) = \partial_t^2 \partial_\nu u(y,t) -\int_{0}^{\zeta(x,y)} \partial_t^3 \partial_\nu u(y,t-\tau) \, d\tau
\end{align*}
and utilize the Cauchy-Schwarz inequality and a-priori estimate \eqref{a-priori estimate} to estimate
\begin{align*}
   & \left\vert \int_{D \times \partial D} \frac{\alpha(y) m(x,y)}{k_0(x) \zeta(x,y)} \int_{0}^{\zeta(x,y)} \partial_t^3 \partial_\nu u(y,t-\tau)\, d\tau \, dx \, dS_y\right\vert \preceq \epsilon^{-2} |D| \cdot |\partial D|^{1/2}  \sup_{0\le\tau\le {t}} \, \|\partial^3_t \partial_\nu u(\cdot,\tau)\|_{L^2(\partial D)} 
     \preceq \epsilon^4. 
\end{align*} 
 To elaborate, we apply Theorem \ref{gradient estimate} to bound $L^2(D)$ norms of $\partial_t^3\partial_\nu u$ and its supremum. Therefore, we need $p\ge 7$. In a similar manner, we derive
\begin{align*}
    \left\vert \int_{\partial D \times \mathbb{R}} \int_{D} \alpha(y) g(x,t-\tau;y) \, \partial_t^2 \partial_\nu u(y,\tau) \, dx \, dS_y \, d\tau \right\vert  \preceq \epsilon^{-2} \, |D| \cdot\, |\partial D|^{1/2} \, \sup_{0\le \tau\le t}\|\partial_t^2 \partial_\nu u(\cdot,\tau)\|_{L^2(\partial D)} \preceq \epsilon^{4}. 
\end{align*}
In light of \eqref{comparing values_2} and 
\begin{align*}
   &\int_{D \times \partial D} \frac{\alpha(y) m(x,y)}{k_0(x) \zeta(x,y)} \partial_t^2 \partial_\nu u(y,t-\zeta(x,y)) \, dx \, dS_y \\
   & = \frac{1}{\rho_1\,c_0^2(z)}\int_{D \times \partial D} \frac{\partial_t^2 \partial_\nu u(y,t)}{|x-y|} \, dx \, dS_y \, + \int_{D \times \partial D} \left( \frac{O(\epsilon^{-1})}{|x-y|} + O(\epsilon^{-2})\right)\, \partial_t^2 \partial_\nu u(y,t)\, dx \, dS_y + O(\epsilon^4).
\end{align*}
Likewise in the earlier step, we make use of Theorem \ref{gradient estimate} and arrive at the estimates
\begin{align*}
    \left\vert \int_{D \times \partial D} O(\epsilon^{-1}) \frac{\partial_t^2 \partial_\nu u(y,t)}{|x-y|} \, dx \, dS_y\right\vert  & \preceq \epsilon^{-1} \left(\int_{D \times \partial D} \frac{dx \, dS_y}{|x-y|^2} \right)^{1/2} \left(\int_{D \times \partial D} |\partial_t^2 \partial_\nu u(y,t)|^2 \, dx \, dS_y\right)^{1/2} \\
    & \preceq \epsilon^{-1} \left( |\partial D| \sup_{y\in\partial D}\int_{D}\frac{dx}{|x-y|^2}\right)^{1/2} | D|^{1/2} \|\partial_t^2 \partial_\nu u(\cdot,t)\|_{L^2(\partial D)} \\
    & \preceq \epsilon^{-1} \cdot \epsilon^2 \cdot \epsilon^{3/2} \cdot \epsilon^2 \preceq \epsilon^{9/2}, 
\end{align*}  
and, also 
\begin{align*}
   \left\vert \int_{D \times \partial D} O(\epsilon^{-2}) \, \partial_t^2 \partial_\nu u(y,t) \, dx \, dS_y \right\vert & \preceq \epsilon^{-2}\, |D| \cdot |\partial D|^{1/2} \, \|\partial_t^2 \partial_\nu u(\cdot,t)\|_{L^2(\partial D)} \preceq \epsilon^4.
\end{align*}

\medskip

Furthermore, we write
\begin{align}\label{seperating Neumann averages}
    \nonumber & \frac{1}{\rho_1\,c_0^2(z)}\int_{D \times \partial D} \frac{\partial_t^2 \partial_\nu u(y,t)}{|x-y|} \, dx \, dS_y \\
     & =  \frac{A_{\partial D}}{\rho_1 c_0^2(z)} \int_{\partial D} \partial_t^2 \partial_\nu u(y,t) \, dS_y +  \frac{1}{\rho_1 c_0^2(z)} \int_{\partial D} \left( A(y) - A_{\partial D} \right) \partial_t^2 \partial_\nu u(y,t)\, dS_y.
\end{align}
Here, the constant $A_{\partial D}$ is defined as
\begin{align*}
   A_{\partial D}:= \fint_{\partial D} A(y) \, dS_y, \textnormal{ and }  A(y) = \int_D \frac{dx}{|x-y|} , \quad \textnormal{ for } y \in \partial D.
\end{align*}
Let us first analyze the constant $A_{\partial D}$ in terms of parameter $\epsilon$. We take note of the standard uniform bounds
\begin{align}\label{upper_lower_bound}
    \frac{2\pi}{3}\, \epsilon^2 \le A(y) \le 8\pi\epsilon^2, \quad  \textnormal{for } y \in \partial D ,
\end{align}
which readily implies $A_{\partial D} \simeq \epsilon^2$. Proving the claim \eqref{upper_lower_bound} is elementary since,
$ |x-y| \le |x-z| + |z-y| \le 2\epsilon,$ and $D \subseteq B(y,2\epsilon)$,
for $x,y\in D$.  As an outcome, we obtain
\begin{align*}
   \frac{2\pi}{3}\epsilon^2= \frac{1}{2\epsilon} \int_{D} dx \le A(y) \le \int_{ B(y,2\epsilon)}\frac{dx}{|x-y|} = \int_0^{2\epsilon} 4\pi r \, dr  = 8\pi\epsilon^2.
\end{align*}

\medskip

Now we derive that the second term in right hand side of \eqref{seperating Neumann averages} contributes in the error term and enjoys better estimate than the case when one merely applies Theorem \ref{gradient estimate}. This is manifested in the following Lemma where we discuss the behaviour of the Neumann derivative of $\partial_t^2 u(\cdot,t)$ on the average-zero subspace of $L^2(\partial D)$ defined by 
\begin{align*}
    L^2_0(\partial D) := \left\{f\in L^2(\partial D); \, \fint_{\partial D} f(y) \, dS_y = 0\right\},
\end{align*}
which is endowed with the standard $L^2(\partial D)$ norm. For any $f\in L^2(\partial D)$, we denote $\mathbb{P}f$ to denote the component of $f$ that is projected onto the subspace $L^2_0(\partial D)$. For $f\in L^2(\partial D)$, it is trivial to notice that
\begin{align}\label{projection}
    \mathbb{P}f(x) = f(x) - \fint_{\partial D}f, \, \textnormal{ for }  x\in \partial D, \quad \int_{\partial D} f(x) \phi(x) \, dS_x =\int_{\partial D} \mathbb{P}f(x) \, \phi(x) \, dS_x  , \,\textnormal{ for } \phi \in L^2_0(\partial D).
\end{align} 
Without proof, we state the following Lemma which underscores an important estimate which says 
\begin{lemma}\label{average zero estimate}
    For $u$ solving the IVP \eqref{IVP_u}, we have the estimate
    \begin{align*}
       \left \|\mathbb{P}\left[\partial_t^2 \partial_\nu u (\cdot,t)\right] \right\|_{L^2_0(\partial D)} \preceq \epsilon^3 , \quad \textnormal{ for a.e. } t \in (0,T).
    \end{align*}
\end{lemma}
For a proof of Lemma \ref{average zero estimate}, we refer to Section \ref{Appendix}. Now, on account of Lemma \ref{average zero estimate} and \eqref{projection}, we observe 
\begin{align*}
   \left| \frac{1}{\rho_1 \,c_0^2(z)} \int_{\partial D} \left( A(y) - A_{\partial D} \right) \partial_t^2 \partial_\nu u(y,t)\, dS_y \right| \preceq \epsilon^{-2} \|A(\cdot)\|_{L^2_0(\partial D)} \left\|\mathbb{P}\left[\partial_t^2 \partial_\nu u(\cdot,t)\right]\right\|_{L^2_0(\partial D)} \preceq \epsilon^4.
\end{align*}
Combining all of these observations, we express 
\begin{align}\label{T_1}
    \int_{\partial D} \rho_0^{-1}(x) \partial_\nu\mathcal{S}_G\left[\alpha\partial_\nu u\right](x,t) \, dS_x = -  \frac{1}{2} \int_{\partial D} \alpha(x) \partial_\nu u(x,t) \, dS_x + \frac{A_{\partial D}}{\rho_1 c_0^2(z)} \int_{\partial D} \partial_t^2 \partial_\nu u(y,t) \, dS_y + O(\epsilon^4)
\end{align}
which completes the discussion on averaging the term $T_1$. 

\bigskip

Performing boundary integration to \eqref{Neumann Lippman-Schwinger} and using the relations \eqref{T_2},  \eqref{T_3} and \eqref{T_1}, we see
\begin{align}
   \nonumber  \int_{\partial D} \rho_0^{-1}(x) \partial_\nu u(x,t) \, dS_x + \frac{A_{\partial D}\, \rho_1}{c_0^2(z)} \int_{\partial D} \partial_t^2 \partial_\nu u(x,t) \, dS_x  &
    = \int_{\partial D}\rho_0^{-1}(x) \,\partial_\nu v(x,t) \, dS_x + \int_{D} J(x,t)\, dx \\ 
    \label{pre ODE} & \quad + \gamma(z)\,\frac{ k_1}{\rho_1} \int_{\partial D} \partial_\nu u(x,t)\, dS_x + O(\epsilon^4).
\end{align}
To associate \eqref{pre ODE} with an ODE, let us set $\mathcal{U}(t) = \int_{\partial D} \partial_{\nu} u(y,t) \, dS_y$, and take note of the facts
\begin{align*}
   & \int_{\partial D} \rho_0^{-1}(x) \, \partial_\nu u(x,t) \, dS_x =  \rho_0^{-1}(z) \,\mathcal{U}(t) \, + O(\epsilon^4), \quad\, 1 - \gamma(z) \rho_0(z)\,\frac{k_1}{\rho_1}  = c_0^{-2}(z) \frac{k_1}{\rho_1}, \\
    & \int_{\partial D}\rho_0^{-1}(x) \partial_\nu v(x,t) \, dS_x + \int_{D} J(x,t)\, dx = \frac{|D|}{k_0(z)} v_{tt}(z,t) + O(\epsilon^4).
\end{align*}
Taking the preceding facts into consideration, we arrive at 
\begin{align}\label{ODE_Neumann average}
     \begin{cases}
     \partial^2_t \,\mathcal{U}(t) + \frac{2k_1}{\rho_0(z)A_{\partial D}} \, \mathcal{U}(t) =  \frac{2\rho_1 \, |D|}{\rho_0(z) A_{\partial D}} \, v_{tt}(z,t) \, + \, O(\epsilon^4), \\
     \mathcal{U}(0) = 0, \, \partial_t\mathcal{U}(0)= 0,
     \end{cases}
\end{align}
where the causality of the wave-field $u$ allows us to incorporate zero initial conditions to $\mathcal{U}$.  Solving \eqref{ODE_Neumann average}, we notice that 
\begin{align*}
  \mathcal{U}(t) = \frac{2\rho_1 |D|}{\rho_0(z) A_{\partial D}} \left(\frac{2k_1}{\rho_0(z)A_{\partial D}}\right)^{-1/2} \int_0^{t} \sin\left(\sqrt{\frac{2k_1}{\rho_0(z)A_{\partial D}}}(t-\tau)\right) v_{tt}(z,\tau) \, d\tau \, + \, O(\epsilon^4),
\end{align*}
in point-wise sense. Making use of integration by parts twice in time, we obtain
\begin{align}    \label{Neumann average}
   \mathcal{U}(t) = \frac{2\rho_1 |D|}{\rho_0(z) A_{\partial D}} v(z,t) - \frac{2\rho_1 |D|}{\rho_0(z) A_{\partial D}} \sqrt{\frac{2k_1}{\rho_0(z)\,A_{\partial D}}} \int_0^{t} \sin\left(\sqrt{\frac{2k_1}{\rho_0(z)A_{\partial D}}}(t-\tau)\right) v(z,\tau) \, d\tau\, + \, O(\epsilon^4) . 
\end{align}

\medskip

Now, to obtain the asymptotic expansion \eqref{asymptotic expansion_2}, we substitute \eqref{Neumann average} in the relation \eqref{asymptotic expansion_1} and introduce a time-shifted version of the wave-field $v$ denoted by
\[ \tilde v(x,t) = v(x, t-\zeta(x_0,z)), \quad (x,t) \in Q.\]
Due to the support condition of $g$ and a change of variables, we notice
\begin{align}\label{change of variables_1}
     \int_{0}^t g(x,t-\tau;z) v(z,\tau)\, d\tau  & = \int_{0}^{t-\zeta(x,z)} g(x,t-\tau;z) v(z,\tau)\, d\tau 
     = \int_{0}^{t} g(x,t-\tau+\zeta(x,z);z) \tilde v(z,\tau)\, d\tau, 
\end{align} 
and, similarly 
\begin{align}
   \label{change of variables_2} \int_0^{t-\zeta(x_0,z)} \sin\left(\sqrt{\frac{2k}{\rho_0(z)A_{\partial D}}}(t-\zeta(x_0,z)-\tau)\right) v(z,\tau) \, d\tau = \int_0^{t} \sin\left(\sqrt{\frac{2k}{\rho_0(z)A_{\partial D}}}(t-\tau)\right) \tilde v(z,\tau) \, d\tau.
\end{align} 
Moreover, the same consideration yields
\begin{align}
    \nonumber & \int_{0}^{t-\zeta(x_0,z)} g(x_0,t-\tau;z) \int_{0}^{\tau} \sin\left(\sqrt{\frac{2k}{\rho_0(z)A_{\partial D}}}(\tau-s)\right) v(z,s) \, ds\, d\tau \\
    \label{change of variables_3} & = \int_{0}^{t} \left(\int_{0}^{t-s} g(x_0,\tau+\zeta(x_0,z);z) \sin\left(\sqrt{\frac{2k}{\rho_0(z)A_{\partial D}}}(t-s-\tau)\right) \, d\tau\right) \tilde v(z,s) \, ds.
\end{align}
In essence, employing the relations \eqref{change of variables_1}-\eqref{change of variables_3} into the expansion \eqref{asymptotic expansion_2}, we represent
\begin{align}\label{volterra form_1}
    w(x_0,t) = \alpha(x_0,z) \, \tilde v(z,t) +\mathbb{K} \tilde v(z,\cdot)(t) + O(\epsilon^2), 
\end{align}
where, we denote 
\begin{align}\label{alpha}
    \alpha(x_0,z) = \frac{|D|}{2\pi A_{\partial D}} \frac{m(x_0,z)}{\rho_0(z) \zeta(x_0,z)}, \quad  \mathbb{K}f(t) = \int_0^t K(t-\tau) f(\tau) \, d\tau, \quad \textnormal{ for } f\in L^2(0,T).
\end{align}
That is $\mathbb{K}$ represents a convolution operator having the kernel $K(\cdot)$ given by
\begin{align}\label{kernel}
    \nonumber K(t) & = \frac{2|D|}{\rho_0(z) A_{\partial D}}\, g(x_0,t+\zeta(x_0,z);z) - \frac{\sqrt{k_1}|D|}{\pi \sqrt{2} A^{3/2}_{\partial D}} \frac{m(x_0,z)}{\rho^{3/2}_0(z)\zeta(x,z)} \, \sin\left(t\sqrt{\frac{2k_1}{\rho_0(z)A_{\partial D}}}\right) \\
    & \quad  - \frac{2\sqrt{2k_1}|D| }{A^{3/2}_{\partial D}} \frac{1}{\rho^{3/2}_0(z)} \, \int_0^t g(x_0,\tau+\zeta(x_0,z);z) \sin\left((t-\tau)\sqrt{\frac{2k_1}{\rho_0(z)A_{\partial D}}}\right)\, d\tau.
\end{align}
From the definition of $\alpha(x_0,z)$ and the function $K(\cdot)$, it is clear that both of them scale linearly in terms of the asymptotic parameter $\epsilon$ and hence, the expression 
\begin{align}\label{dominant term_measurement}
    w_d(x_0,t) = \alpha(x_0,z) \, \tilde v(z,t) + \int_0^t K(t-\tau) \,\tilde v(z,\tau) \, d\tau,
\end{align}
represents the dominant term of the wave-field $w(x_0,\cdot)$ when considered point-wise in time.

\section{Proof of Theorem \ref{inverse problem}}\label{section_inverse problem}
Let us continue the discussion with an equivalent expression of \eqref{volterra form_1} that is 
\begin{align}\label{data_compact form}
    w(x_0,t) = \mathbb{K}_{\alpha}\tilde v (z,t) + O(\epsilon^2),
\end{align}
where $\mathbb{K}_{\alpha}:= \alpha \mathbb{I} + \mathbb{K}$, with $\alpha$ denoting the constant $\alpha(x_0,z)$. In parallel to \eqref{volterra form_1} and \eqref{dominant term_measurement}, the dominant part of $w(x_0,\cdot)$ is expressed as $ w_d(x_0,t) =  \mathbb{K}_{\alpha}\tilde v (z,t)$. On account of the causality of the wave-field $v$ and linearity of $\mathbb{K}_\alpha$, we observe that $w_d(x_0,t) =0$ for $t\le \zeta(x_0,z)$. Now, we intend to demonstrate that $w_d(x_0,\cdot)$ does not trivially vanish  after the time-level $t=\zeta(x_0,z)$. To establish this, we may argue as in \cite{Senapati_Sini_Wang_23} which requires only the non-vanishing character of the source function at $t=0$. However, for the purpose of solving the full inverse problem and determining the travel-time function with better resolution, we work with a more restrictive class of source functions which are smooth and admits the non-vanishing assumption on the source function as in Theorem \ref{inverse problem}. Our argument involves  two steps. First, we demonstrate that $\mathbb{K}_{\alpha}: L^2(0,T) \to L^2(0,T)$ is an invertible operator and then we show that $\tilde v(z,\cdot)$ does not identically vanish after $t=\zeta(x_0,z)$. Moreover, we prove that $\tilde v(z,t) \simeq (t-\zeta(x_0,z))^{p+3}$ for $t > \zeta(x_0,z)$. At this point, we emphasize that it is the last step which requires the non-vanishing assumption of the source locally near $z$. Combining these steps, we arrive at the desired result and interpret the quantity $\zeta(x_0,z)$ as elucidated in Section \ref{introduction}.   

\bigskip

From the definition of $\mathbb{K}$, it is apparent that
\begin{align*}
     |\mathbb{K}^nf(t)| = \left\vert \int_0^ t K(t-\tau) \, \mathbb{K}^{n-1}f(\tau)\, d\tau \right\vert \le \|K\|_{\infty} \|\mathbb{K}^{n-1}f\|_{L^2(0,t)} \sqrt{t}, \quad \textnormal{ for } t \in (0,T), \, n \in \mathbb{N},
\end{align*}
which results into
\begin{align}\label{ineq for coposite op}
    \|\mathbb{K}^n f\|_{L^2(0,t)} \le \frac{\|K\|_{\infty}^n \|f\|_{L^2(0,t)}}{2^{n/2}\sqrt{n!}}\, t^n , \quad \textnormal{ for } t \in (0,T), \, n \in \mathbb{N}.
\end{align}
Likewise \cite{Senapati_Sini_Wang_23}, we adapt similar arguments to verify \eqref{ineq for coposite op}. For completion, we briefly sketch the details in the following. In order to establish \eqref{ineq for coposite op}, we rely on the principle of mathematical induction. 

For $n=1$, it is immediate that
\begin{align*}
    \|\mathbb{K}f\|_{L^2(0,t)} = \left(\int_0^t |\mathbb{K}f(s)|^2 \, ds \right)^{1/2} \le \left(\int_0^t \|K\|^2_{\infty} \|f\|^2_{L^2(0,s)} \, s\, ds \right)^{1/2} \le \frac{\|K\|_{\infty} \|f\|_{L^2(0,t)}}{\sqrt{2}} \, t.
\end{align*}
Hence, the inequality \eqref{ineq for coposite op} holds true for $n=1$. Let us suppose, \eqref{ineq for coposite op} is valid for $n=m$, that is
\begin{align*}
    \|\mathbb{K}^m f\|_{L^2(0,t)} \le \frac{\|K\|_{\infty}^m \|f\|_{L^2(0,t)}}{2^{m/2}\sqrt{m!}}\, t^m , \quad \textnormal{ for } t \in (0,T),
\end{align*}
and we aim to show \eqref{ineq for coposite op} for $n=m+1$. In this regard, the Cauchy-Schwarz inequality implies 
\begin{align*}
    \|\mathbb{K}^{m+1} f\|_{L^2(0,t)} = \left(\int_0^t |\mathbb{K}^{m+1} f(s)|^2 \, ds \right)^{1/2} & \le \|K\|_{\infty} \left(\int_0^t \|\mathbb{K}^{m}f\|^2_{L^2(0,s)}\, s\, ds \right)^{1/2} \\
    & \le \frac{\|K\|^{m+1}_{\infty}\|f\|_{L^2(0,t)}}{2^{m/2}\sqrt{m!}} \left(\int_0^t s^{2m+1}\, ds \right)^{1/2} \\
    & = \frac{\|K\|^{m+1}_{\infty}\|f\|_{L^2(0,t)}}{2^{\frac{m+1}{2}}\sqrt{(m+1)!}} \, t^{m+1},
\end{align*}
which verifies \eqref{ineq for coposite op} for $n=m+1$. Therefore, we ensure the validity of \eqref{ineq for coposite op} and hence
\begin{align*}
    \|\mathbb{K}^{n}\|^{1/n}_{L^2(0,T)} \le \frac{\|K\|_{\infty}\,T}{\sqrt{2}\left({\sqrt{n!}}\right)^{1/n}} \, , \quad \textnormal{ for }  n \in \mathbb{N}.
\end{align*}
With the help of Gelfand's theorem on spectral radius, we write 
\begin{align*}
    \rho(\mathbb{K}) := \liminf\limits_{n\to\infty} \|\mathbb{K}^n\|_{L^2(0,T)}^{1/n} \preceq \lim_{n \to \infty} \frac{1}{\left({\sqrt{n!}}\right)^{1/n}} = 0
\end{align*}
which dictates the invertibility of $\mathbb{K}_{\alpha}$ in $L^2(0,T)$. Therefore, from \eqref{data_compact form}, we deduce that
\begin{align*}
    v(z,t) = \mathbb{A} w (x_0,\cdot) \left(t+\zeta(x_0,z)\right) + O(\epsilon),
\end{align*}
to be understood point-wise in $t\in(0,T)$, where we denote $\mathbb{A}:= (\alpha \mathbb{I}+ \mathbb{K})^{-1}$ which can be equivalently expressed as Neumann series of the convolution operator $\mathbb{K}$. In addition, let us underline that there is a scaling factor $\epsilon$ present in $\alpha$ and the kernel of $\mathbb{K}$. This can be readily seen from \eqref{alpha} and \eqref{kernel}. On account of this, the remainder in \eqref{data_compact form} translates to the term $O(\epsilon)$ in the expansion of $v(z,\cdot)$. 

\medskip

Now, we address the second component in our argument which verifies $\tilde v(z,t) \simeq (t-\zeta(x_0,z))^{p+3}$ for $t>\zeta(x_0,z)$. In essence, this follows from the Taylor's expansion. From the assumption on the source function in Theorem \ref{inverse problem}, we recall that $p\in\mathbb{N}$ is the smallest integer such that  $\partial^{p+1}_t J(z,0) \neq 0$. Without loss of generality, we proceed with the consideration that $\partial^{p+1}_t J(z,0) > 0$. The other case can be treated similarly. The smoothness of $J(\cdot,\cdot)$ and the assumption $\partial^{p+1}_t J(z,0) > 0$, lead us to the consideration that there exists $\alpha,\, \delta,\, r >0$, such that we have
\[\partial^{p+1}_t J(y,t) \ge \alpha > 0, \textnormal{ for } (y,t)\in\mathcal{B}_{r,\delta}:= B(z,r)\times [0,\delta].\]
The Taylor's expansion of $J(y,\cdot)$ at $t=0$ yields
\begin{align*}
    J(y,t) = \sum_{k=1}^{p} \frac{\partial^{k}_t J(y,0)}{k!} t^k + \frac{ \partial^{p+1}_t J(y,0)}{(p+1)!} t^{p+1} + \frac{\partial^{p+2}_t J(y,\theta t)}{(p+2)!} t^{p+2} = \frac{\partial^{p+1}_t J(y,0)}{(p+1)!} t^{p+1} + \frac{\partial^{p+2}_t J(y,\theta t)}{(p+2)!} t^{p+2}.
\end{align*}
Revising the choice $\delta$ if necessary, the Taylor's expansion of $J(y,\cdot)$ above indicates
\begin{align}\label{Taylor's expansion}
    J(y,t) \simeq t^{p+1} \textnormal{ for } (y,t) \in \mathcal{B}_{r,\delta}.
\end{align}
Here, we have also utilized the space regularity of $J$ and its time-derivatives. From the expression of $G$ in \eqref{G_expression}, we represent the wave-field $v$ as 
\begin{align}\label{representation of initial wave field}
    v(z,t) =  \frac{1}{4\pi} \int_{\Omega_t} \frac{m(z,y)}{\zeta(z,\,y)} J(y,t-\zeta(z,y)) \, dy \, + \int_{0}^t \int_{\Omega_{t-\tau}} g(z,t-\tau;y) J(y,\tau) \, dy d\tau
\end{align}
where $\Omega_t:= \{y\in\Omega; \, \zeta(y,z) < t\}$ for $t\in\mathbb{R}_{+}$. It is evident that there exists $d_1,\, d_2 > 0$ such that
\begin{align}\label{travel-time bound}
   d_1 |x-y| < \zeta(x,y) < d_2 |x-y|, \quad  B(z,t/{d_2}) \subseteq \Omega_t \subseteq B(z,t/{d_1}).
\end{align}
For $t>0$ sufficiently small, we use \eqref{Taylor's expansion} and \eqref{travel-time bound} to calculate
\begin{align}
 \nonumber \int_{\Omega_t} \frac{m(z,y)}{\zeta(z,\,y)} J(y,t-\zeta(z,y)) \, dy &\succeq \int_{\Omega_t} \frac{\left(t-\zeta(z,y)\right)^{p+1}}{\zeta(z,y)}\, dy\\
  \nonumber & \succeq \int_{B(z,t/{d_2})} \frac{\left(t-d_2|y-z|\right)^{p+1}}{d_2|y-z|} \, dy \\
 \label{first term} & \simeq \int_{0}^t r(t-r)^{p+1} \, dr = \frac{t^{p+3}}{(p+2)(p+3)}
\end{align}
and similarly, 
\begin{align}\label{second term}
        \left|\int_{0}^t \int_{\Omega_{t-\tau}} g(z,t-\tau;y) J(y,\tau) \, dy d\tau \right| \preceq \int_{0}^{t} \int_{B(z,\frac{t-\tau}{d_1})} \tau^{p+1} \, d y\, d\tau \preceq \int_{0}^{t} (t-\tau)^3\, \tau^{p+1} d\tau \simeq t^{p+5}.
\end{align}
As a consequence of \eqref{first term} and \eqref{second term}, the representation of $v$ in \eqref{representation of initial wave field} entails 
\begin{align}\label{estimate for initial wave field}
    v(z,t) \simeq t^{p+3}, \quad t \in (0,\delta),
\end{align}
for $\delta>0$ sufficiently small. Also, we notice
\begin{align}\label{second term_dominant term}
   \left\vert \int_0^t K(t-\tau) \tilde v(z,\tau) \, d\tau \right\vert \preceq \int_{\zeta(x_0,z)}^t \left(\tau-\zeta(x_0,z)\right)^{p+3}\, d\tau \simeq \left(t-\zeta(x_0,z)\right)^{p+4}.
\end{align}
It is readily apparent from \eqref{estimate for initial wave field} and \eqref{second term_dominant term}, that the dominant term in the measurement given by \eqref{dominant term_measurement} satisfies $ w_d(x_0,t) \simeq (t-\zeta(x_0,z))^{p+3}, \, \textnormal{ for } t > \zeta(x_0,z)$. Now, we discuss the aspect of determining the travel-time function from the full data i.e. $w(x_0,\cdot)$. Taking into account the asymptotic parameter $\epsilon$ present in the operator $\mathbb{K}_\alpha$, we notice 
\begin{align*}
    w_d(x_0,t) \simeq \epsilon (t-\zeta(x_0,z))^{p+3}, \, \textnormal{ for } t > \zeta(x_0,z),
\end{align*}
which translates to 
\begin{align}\label{time-jump in the full data}
    w(x_0,t) \succeq  \epsilon \left(t-\zeta(x_0,z)\right)^{p+3} + O(\epsilon^2), \, \textnormal{ for } t > \zeta(x_0,z).
\end{align}
Let us recall the representation \eqref{representation of w} which leads to the fact that, $w(x_0,\cdot)$ vanishes before the time level $t_{*}= \inf_{y\in\overline D} \zeta(x_0,y) = \zeta(x_0,z) + O(\epsilon)$. On account of \eqref{time-jump in the full data}, we arrive at the point that the function $w(x_0,\cdot)$ is strictly positive from the time-level $t_{**}= \zeta(x_0,z) + \epsilon^{1/(p+3)}$ onwards. Therefore, the jump in the graph of $w(x_0,\cdot)$ with respect to time variables determines the quantity $\zeta(x_0,z)$ upto an error having order $\epsilon^{{1}/{(p+3)}}$.  

\section{Reconstruction of the medium properties and the passive source}
Here, we discuss on the reconstruction procedure of the medium properties $\rho_0(\cdot)$ and $k_0(\cdot)$ along with the source function $J(\cdot,\cdot)$ simultaneously. Continuing the discussion in Section \ref{section_inverse problem}, we notice that the travel time function $\zeta(x_0,\cdot)$ can be recovered in $\Omega$ by varying $z$ within $\Omega$ i.e. by injecting the contrast agent at different points in $\Omega$ and then analyzing the non-vanishing aspect (in time) of the respective measurement collected at $x_0\in\partial\Omega$. With this knowledge of $\zeta(x_0,\cdot)$, we make use of the Eikonal equation \eqref{Eikonal eqn} to recover the sound speed $c_0(\cdot)$ in $\Omega$. Other than the errors in determining $\zeta(x_0,\cdot)$, we should take note of the point  that the reconstruction of $c_0(\cdot)$ will also be susceptible to errors arising from numerical differentiation (through solving the Eikonal equation \eqref{Eikonal eqn}). In order to reconstruct $\rho_0(\cdot)$ and $k_0(\cdot)$, we need to separate them out from the wave-speed $c_0(\cdot)$. We concentrate on this next. We comprehend this reconstruction procedure by discussing it into two different cases.

\subsection{A partial case : One of the medium properties being known}\label{partial case}

In this case, we assume one of the material parameters to be known. Without loss of generality, we assume the mass density $\rho_0(\cdot)$ is known and therefore aim to reconstruct the bulk modulus $k_0$ and the source $J$. Before proceeding further, we notice that the inverse problem, with regard to the parameter count, is yet formally determined as both the data and unknown functions are four dimensional.  

\medskip

 The discussion in Section \ref{section_inverse problem} implies that we can recover the wave-speed $c_0$ from the data modulo an error. $\rho_0(\cdot)$ being known from our assumption, it is therefore immediate that we can also recover $k_0(\cdot)$. This is due to the fact that $c_0(x)= \sqrt{\frac{k_0(x)}{\rho_0(x)}}, \, x\in\Omega$. In order to complete the desired reconstruction process, it remains to recover the source $J$. Following Theorem \ref{section_inverse problem}, we see that the operator $\mathbb{K}_\alpha$ is invertible and depends explicitly on the coefficients $\rho_0(\cdot)$ and $c_0(\cdot)$. From $\rho_0(\cdot)$ and $c_0(\cdot)$, one can calculate all the geodesics corresponding to the Riemmain metric \eqref{metric}. This in turn will lead us to the knowledge of functions $m(\cdot,\cdot)$ from \eqref{defn_g_m} and also of $g(\cdot,\cdot;\cdot)$ by solving the forward problem as discussed in \cite[Theorem 4.1]{Romanov_book_1987}. As a result, we can explicitly determine the kernel of $\mathbb{K}$ and the coefficient $\alpha$. Next we can rely on Born approximation to compute the inverse of $\mathbb{K}_\alpha$ i.e. $\mathbb{A}$. Here, it is important to recall from Section \ref{section_inverse problem} that we have already showed the spectral radius of $\mathbb{K}$ to be zero. Therefore, we can recover $v$ with an error term from the relation \eqref{data_compact form} as follows
\begin{align*}
    v(z,t) = \mathbb{A} w(x_0,\cdot)(t+\zeta(x_0,z)) + O(\epsilon), \quad \textnormal{ for } t\in (0,T).
\end{align*}
Here, the error term is of order $\epsilon$ in point-wise sense with respect to time variable. As before, we vary $z\in\Omega$ and thus determine the initial wave-field $v(\cdot,\cdot)$. As it is now established that we know the functions $k_0,\, \rho_0$ and $v$, we employ the PDE \eqref{IVP_v} to recover the source function $J(\cdot,\cdot)$. Needless to say, numerical differentiation error also emerges in reconstructing $J$. This completes the discussion on reconstructing the material properties and the source function under the assumption that one of $k_0$ and $\rho_0$ is known a-priori.

\subsection{The general case}\label{general case}
In this case, we address the full resolution of the reconstruction problem modulo the assumption that the spatial function $\left.\partial_t^{p+1} J(\cdot,0)\right\vert_{\Omega}$ is known and never vanishes in $\Omega$. Also here, we revise the data and work with $w_d(x_0,\cdot)$ which, as we have seen earlier, represents the dominant part of $w(x_0,\cdot)$. 
To recover both the medium properties, it suffices to recover one of them since, we have already identified their ratio in terms of the wave speed $c_0$. 

\medskip

To this aim, we start by recalling the dominant part of the data which we currently take as measurement
\begin{align*}
    w_d(x_0,t) = \alpha \, \tilde v(z,t) + \int_0^t K(t-\tau) \tilde v(z,\tau)\, d\tau
\end{align*}
and calculate $(p+3)$-th time-derivative of this revised measurement in the following 
\begin{align}\label{recovering k_0}
    \lim_{t\to\zeta(x_0,z)+}\partial^{p+3}_t w_d(x_0,t) =  \frac{\alpha\, m(z,z)}{4\pi}\,\partial^{p+1}_t J(z,0) \,\lim_{t\to 0+}\partial_t^{p+3} \int_{\Omega_t} \frac{\left(t-\zeta(z,y)\right)^{p+1}}{\zeta(z,y)} \, dy.
\end{align}
While the relation \eqref{recovering k_0} help us determine $\alpha \,m(z,z)$ from the revised data as $\partial^{p+1}_t J(\cdot,0)\vert_{\Omega}$ is known, it further requires to prove that the quantity 
\begin{align}\label{A}
    A= \lim_{t\to 0+}\partial_t^{p+3} \int_{\Omega_t} \frac{\left(t-\zeta(z,y)\right)^{p+1}}{\zeta(z,y)} \, dy,
\end{align}
is a positive constant which depends only on the wave speed $c_0(\cdot)$ that has been recovered earlier. The bulk modulus $k_0(\cdot)$ is encoded within the quantity $\alpha\, m(z,z)$ since
\begin{align}\label{alpha and m}
   \nonumber \alpha \, m(z,z) = \frac{|D|}{2\pi A_{\partial D}} \frac{\sigma(x_0,z) \, \sigma(z,z) \, \left(k_0(z) c_0^{-3}(z)\right)^2}{\zeta(x_0,z) \, \rho_0(z)} & = \frac{|D|}{2\pi A_{\partial D}} \frac{\sigma(x_0,z) \, \sigma(z,z)}{c_0^4(z) \,\zeta(x_0,z)} \, k_0(z) \\
    & = \frac{|D|}{2\pi A_{\partial D}} \frac{|\nabla_x \eta (x_0,z)|^{1/2}}{c_0^4(z) \,\zeta(x_0,z)} \, k^{3/2}_0(z).
\end{align}
In view of \eqref{recovering k_0} and \eqref{alpha and m}, the quantity $k_0(z)$ will be determined. This can be realized through the point we can compute the quantities $\nabla_x\,\eta(x_0,z)$ and $\zeta(x_0,z)$ from the knowledge of wave-speed $c_0(\cdot)$. Therefore, it is enough to show that $A>0$. To accomplish this, we argue as in \eqref{first term} and \eqref{second term} 
\begin{align*}
    \int_{\Omega_t} \frac{\left(t-\zeta(z,y)\right)^{p+1}}{\zeta(z,y)} \, dy \simeq \int_0^t r \, (t-r)^{p+2} \, dr \simeq \, t^{p+3} 
\end{align*}
which in turn implies that $A>0$. Therefore, the only task left is to verify  \eqref{recovering k_0}.
From the regularity of the wave-field $v$ and the estimate \eqref{second term_dominant term}, it is immediate that
\begin{align}\label{time derivative_third term_revised measurement}
   \lim_{t\to\zeta(x_0,z)+}  \partial^{p+3}_t \int_0^t K(t-\tau) \, \tilde v(z,\tau) \, d\tau =  \lim_{t\to 0+} \partial^{p+3}_t \int_{0}^{t} \int_{\Omega_{t-\tau}} g(z,t-\tau;y) J(y,\tau) \, dy \, d\tau = 0.
\end{align}
With reference to \eqref{time derivative_third term_revised measurement}, it boils down to compute the $(p+3)$-th time-derivative of the revised data as
\begin{align}
   \nonumber \lim_{t\to\zeta(x_0,z)+}\partial^{p+3}_t w_d(x_0,t) & = \frac{\alpha}{4\pi} \, \lim_{t\to\zeta(x_0,z)+}\partial^{p+3}_t \int_{\Omega_{t-\zeta(x_0,z)}} \frac{m(z,y)}{\zeta(z,\,y)} J(y,t-\zeta(x_0,z)-\zeta(z,y)) \, dy  \\
     \label{general case_1} & = \frac{\alpha}{4\pi} \,\lim_{t\to 0+}\partial^{p+3}_t \int_{\Omega_t} \frac{m(z,y)}{\zeta(z,\,y)} J(y,t-\zeta(z,y)) \, dy.
\end{align}
Meanwhile, the Taylor's expansion of $J(y,\cdot)$ at $t=0$ yields
\begin{align}
    \nonumber \int_{\Omega_t} \frac{m(z,y)}{\zeta(z,\,y)} J(y,t-\zeta(z,y)) \, dy & = \int_{\Omega_t} \frac{m(z,y)}{\zeta(z,\,y)} \left[\frac{ \partial^{p+1}_t J(y,0)}{(p+1)!} \left(t-\zeta(z,y)\right)^{p+1} \right.\\
   \label{derivative_taylor expansion} & \quad \quad \quad \left. \quad \quad + \, \, \frac{\partial^{p+2}_t J(y,\theta t-\theta \zeta(z,y))}{(p+2)!} \left(t-\zeta(z,y)\right)^{p+2} \right] \, dy.
\end{align}
Similar to \eqref{first term}, we can show that
\begin{align}\label{general case_2}
    \frac{1}{(p+2)!}\lim_{t\to 0+}\partial^{p+3}_t \int_{\Omega_t} \frac{m(z,y)}{\zeta(z,\,y)} \partial^{p+2}_t J(y,\theta t-\theta \zeta(z,y)) \left(t-\zeta(z,y)\right)^{p+2} \, dy= 0.
\end{align}
We utilize the smoothness of $m(z,\cdot)$ and $J(\cdot,0)$ to deal with the first part of right hand side of \eqref{derivative_taylor expansion}.
\begin{align*}
     \int_{\Omega_t} \frac{m(z,y)}{\zeta(z,\,y)} \partial^{p+1}_t J(y,0) \left(t-\zeta(z,y)\right)^{p+1} \, dy & = m(z,z)\,\partial^{p+1}_t J(z,0) \int_{\Omega_t} \frac{ \left(t-\zeta(z,y)\right)^{p+1}}{\zeta(z,y)} \, dy \\
     & \quad + \int_{\Omega_t} O(|y-z|)\frac{ \left(t-\zeta(z,y)\right)^{p+1}}{\zeta(z,y)} \, dy.
\end{align*}
Following same steps as before, we note that $\lim\limits_{t\to 0+}\partial^{p+3}_t  \int_{\Omega_t} O(|y-z|)\frac{ \left(t-\zeta(z,y)\right)^{p+1}}{\zeta(z,y)} \, dy = 0$,
and hence, conclude
\begin{align}\label{general case_3}
    \frac{\alpha}{4\pi} \lim_{t\to 0+}\partial_t^{p+3} \, \int_{\Omega_t} \frac{m(z,y)}{\zeta(z,\,y)} J(y,t-\zeta(z,y)) \, dy = \frac{\alpha\, m(z,z)}{4\pi}\,\partial^{p+1}_t J(z,0) \, A. 
\end{align}
Thus, the verification of \eqref{recovering k_0} is complete. 

To conclude the discussion on the recovery of medium properties, let us make a few observations. In \eqref{general case_3}, the constant $A$ can be explicitly calculated from \eqref{A} as we know the wave speed $c_0(\cdot)$ and hence the travel time function $\zeta(\cdot,\cdot)$  solving \eqref{Eikonal eqn}. Also, the non-vanishing function $\partial^{p+1}_t J(\cdot,0)$ is known by our assumption and the function $\eta(\cdot,\cdot)$ can be calculated from the knowledge of $\zeta(\cdot,\cdot)$ via \eqref{Normal coordiante}. With reference to these observation, we find $k_0(z)$ is the only unknown in the right hand side of \eqref{recovering k_0} or \eqref{general case_3}, while the left hand side of the same equations relates to the measurement. Therefore, we recover the bulk modulus $k_0(\cdot)$. As a consequence of this, we are also able to determine the mass density $\rho_0(\cdot)$ and source function $J(\cdot,\cdot)$ following the arguments in Section \ref{partial case}.

\section{Reconstruction of the medium properties using active sources}\label{active source}
    In this section, we discuss the case of reconstructing medium properties using an active source and injecting microscaled agents at different points of $\Omega$. As earlier, our dateset here consists of functions $\left[v(x_0,t), \, u(x_0,t;z)\right]_{z\in\Omega,\, t\in(0,T)}$. Quite naturally, we no longer need Theorem \ref{inverse problem} in this case as it concerns the aspect of determining the source function which is assumed to be known here. Let us briefly discuss on this as most of the arguments are already present in Section \ref{section_inverse problem}, Subsection \ref{partial case} and \ref{general case}. 

\subsection{Use of regular sources}
    First, we describe a family of source functions to be used here that admits sufficient regularity and non-vanishing assumption at all points of $\Omega$ as stated in Theorem \ref{inverse problem}. We consider
    \begin{align*}
    J(x,t) = \begin{cases}
             t^{p+1} \psi(t)\, \phi(x), \quad t\in\mathbb{R}_+, \\
              0, \quad t\le 0,
             \end{cases}
    \end{align*}
    where $\phi(\cdot)$ is an smooth function (at least $C^{p+3}$) in $\mathbb{R}^3$ which does not vanish in $\Omega$ and $\psi(\cdot)$ is a compactly supported and smooth temporal function (at least $C^{p+3}$) satisfying 
    \begin{align*}
    \psi(t) = 1, \quad \textnormal{for } t\in[-T,T].
    \end{align*}
    In view of this choice, we notice that $\partial_t^{p+1} J(z,0) \neq 0$ for all $z\in\Omega$. As discussed in the second part of Section \ref{section_inverse problem}, the initial wave-field does not trivially vanish, see \eqref{estimate for initial wave field}. This results in the time-jump of the measurement $w(x_0,\cdot)$ at $t=\zeta(x_0,z)$, which in turn determines the travel-time function $\zeta(x_0,\cdot)$ and hence the wave speed $c_0(\cdot)$ in $\Omega$ after a numerical differentiation via the Eikonal equation \eqref{Eikonal eqn}. Now, to separate out the bulk modulus and mass density from the wave-speed, we rely on the relation \eqref{recovering k_0} and the discussion presented in Subsection \ref{general case} which primarily takes into account the time-slope of the measured data.  


\subsection{Use of point sources}
    Another way of proceeding is to use point sources which are localized near a fixed point located outside $\Omega$. In view of \eqref{IVP_v}, we can describe this by considering $J(x,t) = \lambda(t)\,\delta(x-x_0)$. Here, $\lambda(\cdot)$ is a sufficiently regular causal function and $x_0\in\partial\Omega$ denotes the source point where we set the point source. Also, we assume that the source and receivers are placed at two different points lying within close distance. The reason behind this will be clear in due time. Let us consider that we record temporal measurements at $x_c \in \partial\Omega$ before and after injecting micro-scaled agents into the medium. For $x \neq x_0$, we use the structure of the Green's function in Theorem \ref{Green function} to express the initial wave-field as 
    \begin{align*}
        v(x,t) & = \int_{\mathbb{R}^3\times\mathbb{R}_+} G(x,t;y,\tau) \,\lambda(\tau)\, \delta(y-x_0) \, dy \, d\tau \\
       & = \frac{m(x,x_0)}{4\pi\zeta(x,x_0)} \lambda(t-\zeta(x,x_0)) + \int_0^{t-\zeta(x,x_0)} g(x,t-\tau;x_0) \lambda(\tau)\, d\tau
    \end{align*}
    where we also utilized the support condition of $g$ in time in last line. Therefore, it is immediate that $v(x,t)=0$ for $t\le \zeta(x,x_0)$. Without deriving the a-priori estimates for $v$ and $u$ as done in Theorem \ref{gradient estimate}, we formally discuss the reconstruction algorithm in this context. Note here that, our arguments in Theorem \ref{gradient estimate} mostly relied on the analysis of 
    the IVP satisfied by $w:= u-v$, as 
\begin{align*}
    \begin{cases}
     k_0^{-1}(x) w_{tt} - \textnormal{div}\left(\rho_0^{-1}(x) \nabla w\right) = \left(\frac{1}{k_0(x)} - \frac{1}{k(x)}\right) u_{tt} - \textnormal{div}\left(\left(\frac{1}{\rho_0(x)} - \frac{1}{\rho(x)}\right) \nabla u\right), \, \textnormal{ in } \, \mathbb{R}^3\times \mathbb{R}_+, \\
     w(\cdot,0) = w_t(\cdot,0) = 0, \textnormal{ in } \mathbb{R}^3.
    \end{cases}   
\end{align*}
Now, if $\lambda(\cdot)$ does not trivially vanish after $t=0$, we can conclude that $v(x,\cdot)$ also does not trivially vanish just after the time-level $t=\zeta(x,x_0)$. This is because 
\begin{align}\label{v_point source_rep}
    v(x,t) = \frac{m(x,x_0)}{4\pi\zeta(x,x_0)} \tilde\lambda(t) + \int_{0}^{t} g(x,t-\tau+\zeta(x,x_0);x_0) \, \tilde\lambda(\tau) \, d\tau
\end{align}
which is a second-kind Volterra operator acting on the function $\tilde \lambda(t):= \lambda(t-\zeta(x,x_0)$. On account of the invertibility proven for such integral operators in Section \ref{section_inverse problem}, we conclude that $v(x,\cdot)$ experiences a jump (in time) at $t=\zeta(x,x_0)$. This underlines the distinctive feature of the present case from the one considered earlier. It is due to the fact that the wave initiated at $x_0\in\partial\Omega$ requires some amount of time to arrive at $x\in\Omega$. However, in the case considered earlier, the function $v(x,\cdot)$ experiences jump right after $t=0$ as the source function was assumed to be initially non-zero throughout $\Omega$. As mentioned in Section \ref{section_inverse problem}, we need to consider $ \lambda_0:=\partial_t^{p+1}\lambda(0) \neq 0$ for a better resolution in determining the wave speed. The constant $\lambda_0$ is further assumed to be known in order to separate out $\rho_0(\cdot)$ and $k_0(\cdot)$ from the wave-speed. Due to this non-vanishing assumption, Taylor's series expansion implies
\begin{align*}
    \lambda(t) = \frac{\lambda_0}{(p+1)!}\, t^{p+1} + O(t^{p+2}), \textnormal{ for } t\to 0+,
\end{align*}
which in view of \eqref{v_point source_rep} results in 
\begin{align}\label{v_point source}
    v(x,t) \simeq (t-\zeta(x,x_0))^{p+1}, \, \textnormal{ for } \, \zeta(x,x_0) < t < \zeta(x,x_0) + \delta
\end{align}
where $\delta>0$ is sufficiently small. Similar to the representation \eqref{volterra form_1}, we can also derive 
\begin{align}\label{w_point source}
    w(x_c,t) = \alpha(x_c,z) \, \tilde v(z,t) +\mathbb{K} \tilde v(z,\cdot)(t) + O(\epsilon^2)
\end{align}
which is to be understood point-wise in $t\in(0,T)$ and $\tilde v(z,t) = v (z,t-\zeta(x_c,z))$. Hence, in consideration of \eqref{v_point source} and \eqref{w_point source}, it is clear that the dominant part of $w(x_0,\cdot)$ i.e. 
\begin{align*}
    w_d(x_c,t) := \alpha(x_c,z) \, \tilde v(z,t) +\mathbb{K} \tilde v(z,\cdot)(t)
\end{align*}
experiences jump at the time-level $t= \zeta(x_c,z) + \zeta(z,x_0)$. Varying the point of insertion of contrast agent in $\Omega$, we first determine the travel-time function $\zeta(x_0,\cdot)$ in $\Omega$ (upto an error that is proportional to the distance between $x_0$ and $x_c$). This is because
\begin{align*}
    \zeta(x_c,z) + \zeta(z,x_0) = 2 \zeta(z,x_0) + O(|x_c-x_0|).
\end{align*}
Assuming $x_c$ and $x_0$ are close, we therefore recover travel-time function $\zeta(z,x_0)$ upto a small error of order $|x_c-x_0|$. Consequently, we have the wave-speed $c_0(\cdot)$ by means of the Eikonal equation \eqref{Eikonal eqn}. To separate out the bulk and mass density, we take into account the time-slope of the measured data as before and then vary $z$ within $\Omega$. We omit the details in this regard as it is similar to the one discussed in Subsection \ref{general case}. 

%

\section{Proofs of Lemma \ref{H^2 estimate} and Lemma \ref{average zero estimate}}\label{Appendix}
\subsection{Proof of Lemma \ref{H^2 estimate}} 
Let us recall that our objective here is to prove the estimate 
\begin{equation}\label{grad-div-estimate}
\|\nabla \textnormal{div}(m-n_*)(\cdot,t)\|_{L^2(D)} \preceq \sqrt\epsilon, \quad \textnormal{ for a.e. } t\in(0,T).
\end{equation}
However, we derive this estimate for all possible second order spatial derivatives for the function $(m-n_*)$. That is, we show $\|\partial^2_{ij}(m-n_*)(\cdot,t)\|_{L^2(D)} \preceq \sqrt\epsilon$ for a.e. $t\in(0,T)$ and all $i, \, j \in\{1,2,3\}$. In order to achieve this, we first note from Theorem \ref{Green function} that we have 
\begin{align*}
    (m-n_*)(x,t) = \int_{\mathbb{R}^3 \times \mathbb{R}} K(x,t;y,\tau) \, \alpha(y) \nabla u(y,\tau) \, dy\, d\tau  = \sum_{k=1}^3 T_k [\alpha\nabla u](x,t)
\end{align*}
where the operators $\{T_k\}_{1\le k \le 3}$ are defined in the following
\begin{align*}
   T_1 [\alpha\nabla u](x,t) & = \int_{\mathbb{R}^3} \left(\frac{m(x,y)}{4\pi\zeta(x,y)} - \frac{\rho_0(z)}{4\pi |x-y|}\right) \, \alpha(y) \, \nabla u(y,t-\zeta(x,y)) \, dy, \\
T_2 [\alpha\nabla u](x,t) & = \int_{\mathbb{R}^3} \frac{\rho_0(z)\, \alpha(y)}{4\pi |x-y|} \, \int_{c_0^{-1}(z)|x-y|}^{\zeta(x,y)} \partial_t \nabla u(y,t-\tau)\, dy\, d\tau, \\
T_3 [\alpha\nabla u](x,t) & = \int_{\mathbb{R}^3 \times \mathbb{R}} g(x,t-\tau;y) \,  \alpha(y)\, \nabla u(y,\tau) \, dy \, d\tau.
\end{align*}
Therefore, to prove Lemma \ref{H^2 estimate}, it suffices to show that 
\begin{align}\label{T estimates}
    \|\partial_{ij}^2 T_k[\alpha\nabla u](\cdot,t)\|_{L^2(D)} \preceq \sqrt\epsilon, \quad \textnormal{ for a.e. } t\in(0,T), \textnormal{ and }\, i,j,k\in\{1,2,3\}.
\end{align}
The key idea in proving the estimation \eqref{T estimates} is identifying the kernels present in the operators $T_k$ when exposed to spatial differentiation. Upon identifying the most singular kernel, we resort to the use of mapping properties of Newtonian potential to estimate the related integral operators. However, we use a generalized Young's inequality \cite[Lemma 0.10]{Folland_PDE_book} to take care of the less singular operators.  Needless to say, this argument will require a comparison between derivatives of reciprocal of the travel-time function $\zeta(x,y)$ and its euclidean analog in $D$. 

To enhance clarity, we divide the proof into three steps. In Step I, we derive pointwise comparison between several quantities related to $\zeta(x,y)$ and its derivatives with their Euclidean analog. These relations will be utilized in Step II to work out the estimates for $\|\partial_{ij}^2 T_1[\alpha\nabla u](\cdot,t)\|_{L^2(D)}$. In Step III, we derive the estimates for $\|\partial_{ij}^2 T_k[\alpha\nabla u](\cdot,t)\|_{L^2(D)}$ for $k=2,\,3$.

\medskip

\textbf{Step I:} We start with describing a slightly more general expansion of the travel time-function discussed in \eqref{reln_travel time}. The bicharecteristics to \eqref{Eikonal eqn} i.e. $(x(t),\, p(t))$ are given by (see also \cite[Chapter 3]{Romanov_book_1987}) the system:
\begin{align*}
    \begin{cases}
        \frac{dx}{dt} = c_0^2(x) \, p, \\
        x(0) = y,
    \end{cases}
    \begin{cases}
        \frac{dp}{dt} = - \nabla_x \ln c_0(x),\\
        p(0) = p_0,
    \end{cases}
\end{align*}
and, as a consequence,
\begin{align*}
    \left.\frac{d^2x}{dt^2}\right\vert_{t=0} & = c_0^2(x) \, \left.\frac{dp}{dt}\right\vert_{t=0} - 2 c_0(x) p \, \left.\nabla_x c_0(x) \cdot \frac{dx}{dt}\right\vert_{t=0} \\
    &  = - c^2_0(y) \nabla_x \ln c_0(y) - 2 c^2_0(y) p_0 \, \nabla_x \ln c_0(y) \cdot c_0^2(y) p_0.
\end{align*}
In terms of the Riemannian coordinate $\xi = c_0^2(y) p_0 \, t $, we can therefore express
\begin{align*}
    x = f(\xi,y) & = y + \left.\frac{dx}{dt}\right\vert_{t=0} t + \frac{1}{2}\left.\frac{d^2x}{dt^2}\right\vert_{t=0} t^2 + O(t^3) \\ 
    & = y + c_0^2(y) p_0 \, t - \frac{1}{2} c^2_0(y) \,t^2 \nabla_x \ln c_0(y) - c^2_0(y) \,p_0 t \, \nabla_x \ln c_0(y) \cdot c_0^2(y) \,p_0 t + O(t^3) \\ 
    &=  y + \xi - \frac{1}{2} \nabla_x \ln c_0(y) |\xi|^2 - \xi \nabla_x \ln c_0(y) \cdot \xi + O(|\xi|^3).
\end{align*}
The above Taylor series expansion of $f(\cdot,y)$ near $\xi=0$ ensures us that $\nabla_\xi f(0,y)= \mathbb{I}$, and then due to the Inverse function theorem, we notice that $f(\cdot,y)$ admits a unique inverse function $\eta$ to $f$ near $\xi=0$ which satisfies $\eta(y,y)=0$ and $\nabla_x \eta(y,y)= \mathbb{Id}$. The travel-time function is defined as
\begin{align*}
    \zeta(x,y) = c_0^{-1}(y) |\eta(x,y)|.
\end{align*}
At this stage, let us also mention the estimates on higher order derivatives of $\zeta(\cdot,\cdot)$ from  \cite[Lemma 3.4]{Romanov_book_1987} 
\begin{align}\label{differentiating travel time}
   |\partial_i \zeta(x,y)| \preceq 1, \quad |\partial^2_{ij} \zeta(x,y)| \preceq |x-y|^{-1}, \quad x \neq y.
\end{align}
Similar to the equation (3.3) in \cite{Senapati_Sini_Wang_23}, we derive from the continuity of $m(\cdot,\cdot)$ at $z$ 
\begin{align}\label{comparing the value}
    \frac{m(x,y)}{4\pi\zeta(x,y)} = \frac{\rho_0(z)}{4\pi |x-y|} + \frac{O(\epsilon)}{|x-y|} + O(1).
\end{align}
For $m,i,j\in\{1,2,3\}$, we further note
\begin{align}\label{second derivative_space points}
   \left.\partial^2_{ij} f_m(\xi,y)\right\vert_{\xi=0} = - \delta_{ij}\partial_m \ln c_0(y)- \delta_{jm}\partial_i \ln c_0(y) - \delta_{im}\partial_j \ln c_0(y), 
\end{align}
and now employing the chain rule successively to the equality $\eta(f)(\xi)=\xi$, for $\xi$ small, we find 
\begin{align*}
    \sum_k \partial_k \eta_m \, \partial_i f_k = \delta_{mi}, \quad \quad \sum_{k,l} \partial^2_{kl} \eta_m \, \partial_j f_l \, \partial_i f_k + \sum_k \partial_k \eta_m \partial^2_{ij} f_k = 0
\end{align*}
which results in $a^m_{ij}(y):=\partial^2_{ij}\eta_m = - \partial^2_{ij}f_m$ at $\xi=0$ i.e. $x=y$. The later is calculated in \eqref{second derivative_space points}. The function $g$ being sufficiently smooth, we have
\begin{align}\label{g_expansion}
    \eta(x,y) = x-y + \sum_{|\alpha|=2} \frac{D^\alpha \eta(y,y)}{\alpha!} (x-y)^\alpha + \sum_{|\beta|=3} R_\beta(x,y) \, (x-y)^\beta
\end{align}
with $R_\beta(x,y) \to 0$ as $x\to y$ for all multi-indices $\beta$. Said differently, we have, for $m \in\{1,2,3\}$,
\begin{align*}
    \eta_m(x,y) = x_m -y_m + \sum_{p,q} \frac{a^m_{pq}(y)}{2}\, (x_p-y_p)(x_q-y_q) + \sum_{p,q,r} b^m_{pqr}(x,y)\, (x_p-y_p)(x_q-y_q)(x_r-y_r).
\end{align*}
which along with the symmetry in right hand side of \eqref{second derivative_space points} leads to 
\begin{align}\label{g_sum}
   \sum_{m=1}^{3}\eta_m(x,y) \, \partial_j \eta_m(x,y)  =
      x_j - y_j + \frac{3}{2} \sum_{p,q} a^j_{pq}(y) \, (x_p-y_p)(x_q-y_q) + O(|x-y|^3).
\end{align}
Now, we utilize \eqref{g_sum} and the definition of $\zeta$ to deduce, for $i,j \in \{1,2,3\}$, that 
\begin{align}
    \label{travel time_1st derivative} \partial_j \left( \frac{1}{|\eta(x,y)|^3}\right) & = -\frac{3}{|\eta(x,y)|^5} \sum_{m=1}^{3} \eta_m(x,y) \, \partial_j \eta_m(x,y) 
    = \partial_j \left( \frac{1}{|x-y|^3}\right) + O(|x-y|^{-3}), \quad x\neq y.
\end{align}    
Making use of the regularity of $m$ and the relations \eqref{g_sum} and \eqref{travel time_1st derivative}, we derive 
\begin{align}\label{P_Q_R}
    m(x,y) c_0(y) & = \rho_0(z) + P(x,y),  \quad \quad \quad \sum_{m=1}^{3} \eta_m(x,y) \, \partial_j \eta_m(x,y) = x_j - y_j + R(x,y), \\
   \nonumber |\eta(x,y)|^{-3} & = |x-y|^{-3} + Q(x,y),
\end{align}
where $P,\, Q$ and $R$ satisfy the point-wise estimates for $x,y \in D, \textnormal{ and } x\neq y$,
\begin{align*}
\begin{cases}
 |P(x,y)|\preceq \epsilon, \quad |Q(x,y)|\preceq 1, \quad |R(x,y)| \preceq |x-y|^2, \\
 |\nabla_x P(x,y)| \preceq 1, \quad |\nabla_x Q(x,y)| \preceq |x-y|^{-3}, \quad |\nabla_x R(x,y)| \preceq |x-y|.
\end{cases}
\end{align*}
In view of these relations, we calculate
\begin{align*}
    \partial_j \left(\frac{m(x,y)}{\zeta(x,y)}\right)  & = - \frac{m(x,y)c_0(y)}{|\eta(x,y)|^3}\sum_{m=1}^{3} \eta_m(x,y) \, \partial_j \eta_m(x,y) + \frac{c_0(y)\,\partial_j m(x,y)}{|\eta(x,y)|} \\
    & = - \left(\rho_0(z)+P(x,y)\right)\, \left(\frac{1}{|x-y|^3}+Q(x,y)\right)\, \left(x_j - y_j + R(x,y)\right) + \frac{c_0(y)\,\partial_j m(x,y)}{|\eta(x,y)|} 
\end{align*}
implying
\begin{align}\label{comparing derivatives}
    \nonumber \partial_j \left(\frac{m(x,y)}{\zeta(x,y)}- \frac{\rho_0(z)}{|x-y|}\right) & = \frac{c_0(y)\,\partial_j m(x,y)}{|\eta(x,y)|} - \rho_0(z) \left\{(x_j-y_j)Q + \frac{R}{|x-y|^3} + QR\right\} \\
    \nonumber & \quad \quad + \left\{\frac{x_j-y_j}{|x-y|^3} P + PQ(x_j-y_j) + \frac{PR}{|x-y|^3} + PQ R  \right\} \\
    & = \frac{x_j-y_j}{|x-y|^3} P(x,y) + \widetilde P(x,y)
\end{align}
where $\widetilde P(x,y)$ satisfies the estimate 
\begin{align}\label{P_tilde}
  \left|{\widetilde P}(x,y)\right| \preceq |x-y|^{-1}, \quad  \left|\nabla_x {\widetilde P}(x,y)\right| \preceq |x-y|^{-2} , \quad \textnormal{ for } x,y\in D \, \textnormal{ and } x \neq y.
\end{align}
Here we would like to remark that the relation \eqref{second derivative_space points} is of explicit interest specially when deriving a relation similar to \eqref{comparing derivatives} but for second derivatives. 

\medskip 

\textbf{Step II:} Let us recall the estimate \eqref{a-priori estimate} which we will often use in the following discussion. For $i,j \in \{1,2,3\}$, let us note that 
\begin{align}\label{first term_grad div est}
   \nonumber \partial^2_{ij}T_1[\alpha\nabla u](x,t) & = \partial_{x_i}\int_{\mathbb{R}^3}  \alpha(y) \, \nabla u(y,t-\zeta(x,y)) \,\partial_{x_j} \left[\frac{m(x,y)}{4\pi\zeta(x,y)} - \frac{\rho_0(z)}{4\pi |x-y|}\right] \, dy \\
   \nonumber & \quad \quad - \int_{\mathbb{R}^3}  \alpha(y) \, \partial_t\nabla u(y,t-\zeta(x,y)) \, \left[\frac{m(x,y)}{4\pi\zeta(x,y)} - \frac{\rho_0(z)}{4\pi |x-y|}\right] \,\partial^2_{x_i x_j}\zeta(x,y) \, dy\\
   \nonumber & \quad \quad + \int_{\mathbb{R}^3} \alpha(y) \, \partial^2_t \nabla u(y,t-\zeta(x,y)) \, \partial_{x_i}\zeta(x,y)\,\partial_{x_j}\zeta(x,y)\,\left[\frac{m(x,y)}{4\pi\zeta(x,y)} - \frac{\rho_0(z)}{4\pi |x-y|}\right]\, dy \\
    & \quad \quad - \int_{\mathbb{R}^3} \alpha(y) \, \partial_t \nabla u(y,t-\zeta(x,y)) \, \partial_{x_i}\zeta(x,y)\,\partial_{x_j}\left[\frac{m(x,y)}{4\pi\zeta(x,y)} - \frac{\rho_0(z)}{4\pi |x-y|}\right]\, dy.
\end{align}
We simplify all the terms in \eqref{first term_grad div est} other than the first one. With regard to \eqref{differentiating travel time} and \eqref{comparing the value}, we find
\begin{align*}
    & \left| \int_{\mathbb{R}^3} \alpha(y) \, \partial_t\nabla u(y,t-\zeta(x,y)) \, \left(\frac{m(x,y)}{4\pi\zeta(x,y)} - \frac{\rho_0(z)}{4\pi |x-y|}\right) \,\partial^2_{x_i x_j}\zeta(x,y) \, dy\right| \\
    & \preceq \epsilon^{-1} \sup\limits_{0\le\tau\le t} \int_{D} \frac{|\partial_t \nabla u(y,\tau)|}{|x-y|^2} \, dy + \epsilon^{-2} \sup\limits_{0\le\tau\le t} \int_{D} \frac{|\partial_t \nabla u(y,\tau)|}{|x-y|} \, dy.
\end{align*}
Now, we employ generalized Young's inequality \cite[Lemma 0.10]{Folland_PDE_book} which guaranties 
\begin{align*}
  \sup\limits_{0\le\tau\le t} \left\|\int_{D} \frac{|\partial_t \nabla u(y,\tau)|}{|\cdot-y|^2} \, dy\right\|_{L^2(D)} & \le \left(\sup\limits_{x\in D}\int_{D}\frac{dy}{|x-y|^2}\right) \sup\limits_{0\le\tau\le t}\|\partial_t \nabla u(\cdot,\tau)\|_{L^2(D)} \\
  & \le \epsilon^{3/2} \, \sup\limits_{x\in D}\int_{B(x,2\epsilon)}\frac{dy}{|x-y|^2} \, \preceq \epsilon^{5/2},
\end{align*}
and as a consequence,  for $t\in(0,T)$, we derive
\begin{align*}
   \epsilon^{-1} \, \sup\limits_{0\le\tau\le t} \left\|\int_{D} \frac{|\partial_t \nabla u(y,\tau)|}{|\cdot-y|^2} \, dy\right\|_{L^2(D)} \preceq \epsilon^{3/2}.
\end{align*}
From the Cauchy-Schwarz inequality and a scaling argument, we also have 
\begin{align*}
    \left\|\int_{D}\frac{f(y)}{|\cdot-y|}\, dy \right\|_{L^2(D)} \le \left( \int_{D} |f(y)|^2\, dy\int_{D \times D} \frac{dx\, dy}{|x-y|^2}\right)^{1/2} \le \epsilon^{2} \|f\|_{L^2(D)},
\end{align*}
which implies
\begin{align*}
    \left\|\epsilon^{-2}\, \sup\limits_{0\le\tau\le t}\int_{D}\frac{|\partial_t \nabla u(y,\tau)|}{|\cdot-y|}\, dy \right\|_{L^2(D)} \preceq \sup\limits_{0\le\tau\le t} \|\partial_t \nabla u(\cdot,\tau)\|_{L^2(D)} \preceq \epsilon^{3/2}. 
\end{align*}
In view of the above relations, it follows that
\begin{align*}
    \left| \int_{\mathbb{R}^3} \alpha(y) \, \partial_t\nabla u(y,t-\zeta(x,y)) \, \left(\frac{m(x,y)}{4\pi\zeta(x,y)} - \frac{\rho_0(z)}{4\pi |x-y|}\right) \,\partial^2_{x_i x_j}\zeta(x,y) \, dy\right| \preceq \epsilon^{3/2}. 
\end{align*}
In an essentially similar way, we can also derive
\begin{align*}
    \left\vert\int_{\mathbb{R}^3} \alpha(y)\, \partial^2_t \nabla u(y,t-\zeta(x,y)) \, \partial_{x_i}\zeta(x,y)\,\partial_{x_j}\zeta(x,y)\,\left[\frac{m(x,y)}{4\pi\zeta(x,y)} - \frac{\rho_0(z)}{4\pi |x-y|}\right]\, dy\right\vert \preceq \epsilon^{3/2},
\end{align*}
and 
\begin{align*}
    \left\vert \int_{\mathbb{R}^3} \alpha(y) \, \partial_t \nabla u(y,t-\zeta(x,y)) \, \partial_{x_j}\zeta(x,y)\,\partial_{x_i}\left[\frac{m(x,y)}{4\pi\zeta(x,y)} - \frac{\rho_0(z)}{4\pi |x-y|}\right]\, dy \right\vert \preceq \epsilon^{3/2}.
\end{align*}
where we utilize \eqref{comparing derivatives} in the last estimation. Therefore, it remains to estimate the integral 
\begin{align*}
    I(x,t) & = \partial_{x_i}\int_{\mathbb{R}^3} \alpha(y) \, \nabla u(y,t-\zeta(x,y)) \,\partial_{x_j} \left[\frac{m(x,y)}{4\pi\zeta(x,y)} - \frac{\rho_0(z)}{4\pi |x-y|}\right] \, dy \\
    & = I_*(x,t) + \int_{D} \, \partial_{x_i} \widetilde{P}(x,y) \,\alpha(y) \,\nabla u(y,t-\zeta(x,y))  \, dy,
\end{align*}
where the function $\widetilde P(x,y)$ has been introduced in \eqref{comparing derivatives}. Here, $I_* (x,t) = M(x,t) + N(x,t)$ and we denote
\begin{align*}
    M(x,t) & = \partial_{x_j} \int_D \frac{x_i-y_i}{|x-y|^3}\, P(x,y) \left[\alpha(y)\nabla u(y,t-\zeta(x,y))- \alpha(y)\, \nabla u(y,t)\right] \, dy, \\
    N(x,t) & = \partial_{x_j}  \int_D \frac{x_i-y_i}{|x-y|^3}\, P(x,y) \alpha(y) \nabla u(y,t) \, dy.
\end{align*}
To estimate $\|I(\cdot,t)\|_{L^2(D)}$, we use Young's inequality along with \eqref{P_tilde} to find
\begin{align*}
  \left\|  \int_{D} \, \partial_{x_i} \widetilde{P}(\cdot,y) \,\alpha(y)\, \nabla u(y,t-\zeta(x,y))  \, dy \right\|_{L^2(D)} \preceq \epsilon^{1/2}, \quad t \in (0,T).
\end{align*}
Now, we derive an estimation for $\|I_*(\cdot,t)\|_{L^2(D)}$. For $P$ defined in \eqref{P_Q_R}, we write 
$$P(x,y) = P(y,y) + P_*(x,y)$$
and take note of the following point-wise bounds
\begin{align}\label{further division of P}
    |P(y,y)| \le \epsilon, \quad |P_*(x,y)| \preceq |x-y|, \quad |\partial_i P_*(x,y)| \preceq 1, \quad \textnormal{ for } x ,y \in D.
\end{align}
From straightforward calculations, it is clear that, for $i,j \in \{1,2,3\}$ 
\begin{align*}
    M(x,t) & = \int_D \left( \frac{\delta_{ij}}{|x-y|^3} - \frac{3(x_i-y_i)(x_j-y_j)}{|x-y|^5} \right) P(x,y) \alpha(y)\left[\nabla u(y,t-\zeta(x,y)) - \nabla u(y,t)\right]\, dy \\
    & + \int_D \frac{x_i-y_i}{|x-y|^3} \alpha(y)\,\Big\{\partial_i P(x,y)\, \left[\nabla u(y,t-\zeta(x,y)) - \nabla u(y,t)\right] + P(x,y)\zeta_j(x,y)\,\nabla u(y,t-\zeta(x,y)) \Big\} \, dy 
\end{align*}
and 
\begin{align}\label{N term}
   \nonumber N(x,t) & =  \int_D \left( P_*(x,y)\partial_j \left(\frac{x_i-y_i}{|x-y|^3}\right)+ \frac{(x_i-y_i) \partial_j P_*(x,y)}{|x-y|^3} \right) \alpha(y)\nabla u(y,t) \, dy \\
    &  \quad + \partial^2_{ij} \int_D \frac{1}{|x-y|} \, P(y,y) \alpha(y) \nabla u(y,t)\, dy.
\end{align}
Making use of \eqref{G_expression} and the pointwise estimates of $P,\, Q,\, R$ and their gradients obtained in Step I above, we observe that
\begin{align*}
    |M(x,t)| \preceq \sup_{0\le\tau\le t}\int_{D} \frac{1}{|x-y|} \, \alpha(y) |\nabla u (y,\tau)|\, dy, \quad \textnormal{ for } (x,t)\in D\times (0,T).
\end{align*}
In addition, one can also establish similar bounds for all the terms of $N(x,t)$ other than the last one. For all these terms, we use Young's inequality and  Theorem \ref{gradient estimate}. To exemplify this, we have estimates such as 
\begin{align*}
    \|M(\cdot,t)\|_{L^2(D)} \preceq \sqrt{\epsilon}, \quad \textnormal{ for } t\in(0,T).
\end{align*}
We can derive similar estimate for the first term in \eqref{N term}. Now, to deal with the last expression in \eqref{N term}, we use \eqref{further division of P}, the scaling and mapping properties of Newtonian potential. For $i,j\in\{1,2,3\}$, we see
\begin{align*}
    \left\|\partial^2_{ij} \int_D \frac{1}{|\cdot-y|} \, P(y,y) \alpha(y) \nabla u(y,t)\, dy\right\|_{L^2(D)} \preceq \|P(\cdot,\cdot)\alpha(\cdot)\nabla u(\cdot,t)\|_{L^2(D)} \preceq \epsilon^{1/2} , \quad \textnormal{ for } t \in (0,T).
\end{align*}
Summarizing all of the above discussion, we write
\begin{align*}
    \|\partial^2_{ij} T_1[\alpha\nabla u](\cdot,t)\|_{L^2(D)} \preceq \epsilon^{1/2}, \quad \textnormal{ for } t \in (0,T).
\end{align*}

\medskip

\textbf{Step III:} Now, we aim to derive the estimate \eqref{T estimates} for the operators $T_2$ and $T_3$. Similar to \eqref{first term_grad div est}, we differentiate the expression of $T_2[\alpha\nabla u](x,t)$ twice with respect to space variables to obtain
\begin{align}\label{T2}
   \nonumber \partial^2_{ij} T_2[\alpha\nabla u](x,t) & = \int_{\mathbb{R}^3} \partial_{x_i}\partial_{x_j}\left[\frac{\alpha(y)\,\rho_0(z)}{4\pi |x-y|} \right]\, \int_{c_0^{-1}(z)|x-y|}^{\zeta(x,y)} \partial_t \nabla u(y,t-\tau)\, dy\, d\tau \\
   \nonumber  & \quad \quad +   \int_{\mathbb{R}^3} \frac{\alpha(y)\,\rho_0(z)}{4\pi |x-y|}\, \partial_{x_i}\partial_{x_j}\left[\nabla u(y,t-\zeta(x,y)) - \nabla u(y,t-c_0^{-1}(z)|x-y|)\right] \, dy \\
   \nonumber  & \quad \quad +  \int_{\mathbb{R}^3} \partial_{x_i}\left[\frac{\alpha(y)}\,\rho_0(z){4\pi |x-y|}\right]\, \partial_{x_j}\left[\nabla u(y,t-\zeta(x,y)) - \nabla u(y,t-c_0^{-1}(z)|x-y|)\right] \, dy \\
     & \quad \quad +  \int_{\mathbb{R}^3} \partial_{x_j}\left[\frac{\alpha(y)}\,\rho_0(z){4\pi |x-y|}\right]\, \partial_{x_i}\left[\nabla u(y,t-\zeta(x,y)) - \nabla u(y,t-c_0^{-1}(z)|x-y|)\right] \, dy.
\end{align}
From \eqref{reln_travel time}, it is immediate that
\begin{align}\label{Step 3_1}
    \left\vert \partial_{x_i}\partial_{x_j}\left[\frac{\alpha(y)}\,\rho_0(z){4\pi |x-y|} \right]\, \int_{c_0^{-1}(z)|x-y|}^{\zeta(x,y)} \partial_t \nabla u(y,t-\tau)\, d\tau \right\vert \preceq \frac{\alpha(y)}{|x-y|^2} \sup_{0\le \tau\le t} |\partial_t\nabla u(y,\tau)|, \quad \textnormal{for } x \neq y.
\end{align}
Here we underline that, as $p\ge 3$, the time derivative of $\alpha\nabla u(x,t)$ upto second order and its supremum make sense. For $x\neq y$, a straightforward differentiation and the relation \eqref{differentiating travel time} will lead to the following implication 
\begin{align}\label{Step 3_2}
   \notag \left\vert \frac{\alpha(y)}\,\rho_0(z){4\pi |x-y|}\, \partial_{x_i}\partial_{x_j}\left[\nabla u(y,t-\zeta(x,y)) - \nabla u(y,t-c_0^{-1}(z)|x-y|)\right]  \right\vert \\
    \preceq \frac{\alpha(y)}{|x-y|^2} \sup_{0\le \tau\le t} \left(|\partial_t\nabla u(y,\tau)| + |\partial^2_t\nabla u(y,\tau)|\right) 
\end{align}
and also
\begin{align}
   \label{Step 3_3_1}  \left\vert \partial_{x_i}\left[\frac{\alpha(y)}\,\rho_0(z){4\pi |x-y|}\right]\, \partial_{x_j}\left[\nabla u(y,t-\zeta(x,y)) - \nabla u(y,t-c_0^{-1}(z)|x-y|)\right] \right\vert \preceq \frac{\alpha(y)}{|x-y|^2} \sup_{0\le \tau\le t} \left\vert \partial_t \nabla u(y,\tau)\right\vert, \\
   \label{Step 3_3_2}  \left\vert \partial_{x_j}\left[\frac{\alpha(y)}\,\rho_0(z){4\pi |x-y|}\right]\, \partial_{x_i}\left[\nabla u(y,t-\zeta(x,y)) - \nabla u(y,t-c_0^{-1}(z)|x-y|)\right] \right\vert \preceq \frac{\alpha(y)}{|x-y|^2} \sup_{0\le \tau\le t} \left\vert \partial_t \nabla u(y,\tau)\right\vert.
\end{align}
Now the support condition of $\alpha$ and the relations \eqref{Step 3_1}-\eqref{Step 3_3_2} readily imply 
\begin{align*}
    \left|\partial^2_{ij} T_2[\alpha\nabla u](x,t)\right| \preceq \sup_{0\le \tau \le t} \int_D \frac{\alpha(y)}{|x-y|^2}\left(|\partial_t\nabla u(y,\tau)| + |\partial^2_t\nabla u(y,\tau)|\right)\, dy , \quad \textnormal{ for } (x,t) \in D\times (0,T).
\end{align*}
As before, we refer to the generalized Young's inequality from \cite[Lemma 0.10]{Folland_PDE_book} and use Theorem \ref{gradient estimate} to derive the following estimate
\begin{align*}
    \|\partial^2_{ij}T_2[\alpha\nabla u](\cdot,t)\|_{L^2(D)} \preceq \epsilon^{-1} \sup_{0\le \tau \le t} \left(\|\partial_t \nabla u(\cdot,t)\|_{L^2(D)} + \|\partial^2_t \nabla u(\cdot,t)\|_{L^2(D)}\right) \preceq \epsilon^{1/2}.
\end{align*}
Therefore, we are done with estimating the term $\|\partial^2_{ij}T_2[\alpha\nabla u](\cdot,t)\|_{L^2(D)}$ for $t\in(0,T)$. It only remains to derive an estimate for the term $\|\partial^2_{ij}T_3[\alpha\nabla u](\cdot,t)\|_{L^2(D)}$. Needless to say, the regularity and support condition of $g$ given in Theorem \ref{Green function} are important to derive such estimate. We briefly discuss on the computation. With regard to the regularity and support condition of $g$ in Theorem \ref{Green function}, we take space derivatives to 
\begin{align*}
    T_3[\alpha\nabla u](x,t) = \int_{0}^{t-\zeta(x,y)} \int_D g(x,t-\tau;y) \,\alpha(y)\nabla u(y,\tau) \, dy\,d\tau.
\end{align*}
A straightforward computation implies
\begin{align*}
    & \partial^2_{ij} T_3[\alpha\nabla u](x,t) \\
    &= \int_{0}^{t-\zeta(x,y)} \partial^2_{ij} g(x,t-\tau;y) \,\alpha(y)\nabla u(y,\tau) \, dy\,d\tau - \int_{D} \partial^2_{ij}\zeta(x,y)  g(x,\zeta(x,y); y) \, \alpha(y)\nabla u(y,t-\zeta(x,y)) \, dy \\
    & - \int_{D} \partial_{i}\zeta(x,y) \left[\partial_j g+ \partial_{j}\zeta\,\partial_t g\right](x,\zeta(x,y); y) \, \alpha(y)\nabla u(y,t-\zeta(x,y)) \, dy \\
    & + \int_{D} \alpha(y) \,\partial_i \zeta(x,y) \partial_j \zeta(x,y) \,g(x,\zeta(x,y);y)\, \partial_t \nabla u(y,t-\zeta(x,y)) \, dy,
\end{align*}
which, in view of \eqref{differentiating travel time} and Theorem \ref{Green function}, further gives us the following estimation
\begin{align*}
   \left\vert \partial^2_{ij} T_3[\alpha\nabla u](x,t) \right\vert \preceq \sup\limits_{0\le\tau\le t} \int_D \frac{\alpha(y)}{|x-y|^2} \,\left( |\partial_t\nabla(y,\tau)| + |\nabla u(y,\tau)|\right) \, dy.
\end{align*}
As before, we employ the generalized Young's inequality alongwtih Theorem \ref{gradient estimate} to obtain 
\begin{align*}
    \|\partial^2_{ij}T_3[\alpha\nabla u](\cdot,t)\|_{L^2(D)} \preceq \epsilon^{1/2}, \quad \textnormal{for } t\in (0,T),
\end{align*}
which thus concludes the proof of Lemma \ref{H^2 estimate}.

\subsection{Proof of Lemma \ref{average zero estimate}}
    We demonstrate this only for $\partial_\nu u$, as one can consider time-derivatives to the IVP \eqref{IVP_u} and follow similar steps to be presented here to deal with higher order time derivatives of $\partial_\nu u$. In acknowledgment of this, let us first rewrite the boundary integral equation \eqref{Neumann Lippman-Schwinger} as
    \begin{align}
      \nonumber  \partial_\nu u(x,t) + \frac{\rho_0(z)\,\alpha(z)}{2}  \, \partial_\nu u(x,t) + \alpha(z)\mathcal{K}^t_G [\partial_\nu u](x,t) & = F(x,t) - \mathcal{K}^t_G\left[(\alpha(\cdot)-\alpha(z))\partial_\nu u\right](x,t) \\
      \label{Neumann Lippman-Schwinger_2}  & \quad - \frac{1}{2}\left( \rho_0(x)\alpha(x)-\rho(z)\alpha(z)\right) \partial_\nu u(x,t) 
    \end{align}
   with $F$ denoting the expression 
   \begin{align*}
       F(x,t) := \partial_\nu v(x,t) - \partial_\nu \mathcal{N}_G \left[k\beta J + \gamma\frac{k_1}{\rho_1}\Delta u\right](x,t) - \partial_\nu \mathcal{N}_G \left[\chi_D \nabla\rho_0^{-1}\cdot\nabla u\right](x,t)
   \end{align*}
   which satisfies the bound     
   \begin{align}\label{F_estimate}
       \|F(\cdot,t)\|_{L^2(\partial D)} \preceq  \epsilon, \quad \textnormal{  for a.e } t\in(0,T).
   \end{align}
    To verify the estimation \eqref{F_estimate}, we only concentrate on the function $\partial_\nu \mathcal{N}_G \left[\gamma\frac{k_1}{\rho_1}\Delta u\right](x,t)$ and demonstrate the bound
   \begin{align}\label{F_estimate_1}
       \left\|\partial_\nu \mathcal{N}_G \left[\gamma\frac{k_1}{\rho_1}\Delta u\right](\cdot,t)\right\|_{L^2(\partial D)} \preceq \epsilon.
   \end{align}
   This is due to the fact that either the other terms can be estimated similarly or, 
   the estimate $\|\partial_\nu v(\cdot,t)\|_{L^2(\partial D)} \preceq \epsilon$ is immediate due to spatial boundedness of $v(\cdot,t)$. Let us proceed as in Section \ref{theo_expansion} to notice
   \begin{align*}
       \partial_\nu \mathcal{N}_G \left[\gamma\frac{k_1}{\rho_1}\Delta u\right](x,t) & = \int_D \left( \frac{\partial_\nu m(x,y)}{4\pi\zeta(x,y)} - \frac{m(x,y)\,\partial_\nu \zeta(x,y)}{4\pi\zeta^2(x,y)}\right) \frac{k_1}{\rho_1}\,\gamma(y)\Delta u(y,t-\zeta(x,y))\, dy \\
       & \quad - \int_D \frac{m(x,y)\partial_\nu \zeta(x,y)}{4\pi\zeta(x,y)} \, \frac{k_1}{\rho_1}\,\gamma(y)\, \partial_t\Delta u(y,t-\zeta(x,y))\, dy \\
       & \quad + \int_{D \times \mathbb{R}} \partial_\nu g(x,t-\tau;y) \, \frac{k_1}{\rho_1}\,\gamma(y)\Delta u(y,\tau) \, dy \, d\tau.
   \end{align*}
   Looking at the most singular term from above, we calculate
   \begin{align*}
      & \left( \int_{\partial D} \left|\int_D \frac{m(x,y)\,\partial_\nu \zeta(x,y)}{4\pi\zeta^2(x,y)} \, \frac{k_1}{\rho_1}\,\gamma(y)\Delta u(y,t-\zeta(x,y))\, dy\right|^2 \, dS_x \right)^{1/2} \\
       & \quad \preceq \sup_{0\le\tau\le t} \|\Delta u(\cdot,\tau)\|_{L^2(D)} \, \left(\int_{D\times \partial D} \frac{dy\, dS_x}{|x-y|^4} \right)^{1/2}  \preceq \epsilon^2.
   \end{align*}
   Similar arguments can be employed to estimate the other terms in $\partial_\nu \mathcal{N}_G \left[\gamma\frac{k_1}{\rho_1}\Delta u\right](x,t)$ and conclude \eqref{F_estimate_1}.
   
   Now we identify the dominating part of the operator $\mathcal{K}^t_{G}$ as $\mathcal{K}^t_{G_z}$. In view of this, we will utilize the expansion \eqref{comparing derivatives} and the estimates of $\partial_\nu u$ and its time derivatives as derived in Lemma \ref{gradient estimate}. Here we derive
   \begin{align}\label{principal NP op}
       \mathcal{K}^t_{G}[\partial_\nu u](x,t) =  \mathcal{K}^t_{G_z}[\partial_\nu u](x,t) + O(\epsilon^3)
   \end{align}
   to be understood poin-twise in time and $L^2$-wise in space. To prove \eqref{principal NP op}, we first note that 
   \begin{align*}
       &\left(\mathcal{K}^t_{G} - \mathcal{K}^t_{G_z} \right)[\partial_\nu u](x,t) = \int_{\partial D} \partial_{\nu_x} \left(\frac{m(x,y)}{\zeta(x,y)} - \frac{\rho_0(z)}{|x-y|} \right) \partial_{\nu_y} u(y,t-\zeta(x,y)) \, dS_y \\
       & \quad - \int_{\partial D} \frac{m(x,y)\, \partial_{\nu_x}\zeta(x,y)}{\zeta(x,y)} \,\partial_t\partial_{\nu_y} u(y,t-\zeta(x,y)) \, dS_y + \int_{\partial D \times \mathbb{R}} \partial_{\nu_x} g(x,t-\tau;y) \, \partial_\nu u(y,\tau) \, dS_y\, d\tau. \\
       & \quad \quad + \frac{\rho_0(z)}{c_0(z)}\int_{\partial D} \left(\frac{\nu_x\cdot(x-y)}{|x-y|^2} + \frac{\partial_{\nu_x}\zeta(x,y)}{|x-y|} \right) \partial_t\partial_{\nu_y} u(y,t-c_0^{-1}(z)|x-y|)\, dS_y .
   \end{align*}
   Making use of the smoothness of $\partial D$ as in \cite[Lemma 3.15]{Folland_PDE_book} and a scaling argument, we notice that
   \begin{align}\label{regularity of boundary}
       (x-y)\cdot \nu(x) \le C \epsilon^{-1}|x-y|^2, \quad \textnormal{ for } \,x, \,y \in \partial D,
   \end{align}
   for some $C>0$. From \eqref{differentiating travel time}, \eqref{comparing derivatives} and \eqref{regularity of boundary}, it is immediate that
   \begin{align*}
       \left\vert \left(\mathcal{K}^t_{G} - \mathcal{K}^t_{G_z} \right)[\partial_\nu u](x,t)\right\vert \preceq \, \sup_{0\le \tau\le t}\int_{\partial D} \frac{\left\vert\partial_{\nu} u(y,\tau)\right\vert + \left\vert\partial_t\partial_{\nu} u(y,\tau)\right\vert }{|x-y|} \, d S_y.
   \end{align*}
   Again using the generalized Young's inequality from \cite[Lemma 0.10]{Folland_PDE_book} and the estimates from Lemma \ref{gradient estimate}, we obtain 
   \begin{align*}
       \|\left(\mathcal{K}^t_{G} - \mathcal{K}^t_{G_z} \right)[\partial_\nu u](\cdot,t)\|_{L^2(\partial D)} & \preceq \left(\sup_{x\in\partial D} \int_{\partial D} \frac{dS_y}{|x-y|}\right) \sup_{0\le \tau\le t} \left( \|\partial_\nu u(\cdot,t)\|_{L^2(\partial D)} + \|\partial_t\partial_\nu u(\cdot,t)\|_{L^2(\partial D)} \right) \\
       & \preceq \epsilon^{1+2} = \epsilon^3. 
   \end{align*}
   Therefore the relation \eqref{principal NP op} holds true.

   In light of \eqref{principal NP op}, we recast \eqref{Neumann Lippman-Schwinger_2} in the form
   \begin{align}\label{boundary operator eqn}
      \left[\left(\alpha^{-1}(z)+ \frac{1}{2}\right) \mathbb{I}\textnormal{d} + \mathcal{K}^t_{z} \right][\partial_\nu u](x,t) = F_*(x,t) 
   \end{align}
   where we have introduced the boundary operator $\mathcal{K}^t_z$ defined as follows
   \begin{align*}
       \mathcal{K}^t_{z} [f](x) = -\rho_0(z)\,\int_{\partial D} \frac{(x-y)\cdot \nu(x)}{|x-y|^3} \, f(y) \, dS_y, \quad x\in\partial D
   \end{align*}
   and the term $F_*(x,t)$ is given by 
   \begin{align*}
       F_*(x,t) & := \rho_0(z)\,\int_{\partial D} \frac{(x-y)\cdot \nu(x)}{|x-y|^3} \, \left[\partial_\nu u(y, t - c_0^{-1}(z)|x-y|)-\partial_\nu u(y,t)\right] \, dS_y + \alpha^{-1}(z)F(x,t) \\
     \nonumber & \quad \quad - \frac{1}{\alpha(z)} \,\mathcal{K}_G\left[(\alpha(\cdot)-\alpha(z))\partial_\nu u\right](x,t) - \frac{1}{2\alpha(z)} \left(\rho_0(x)\alpha(x)-\rho_0(z)\alpha(z)\right)\partial_\nu u(x,t) + O(\epsilon^3).
   \end{align*}

   Now, we refer to \cite{Ammari_Kang_book} to conclude the important bound
   \begin{align*}
       \left\|\left[\left(\alpha^{-1}(z)+ \frac{\rho_0(z)}{2}\right) \mathbb{I}\textnormal{d} + \mathcal{K}^t_{z} \right]^{-1}\right\|_{L^2_0(\partial D)} = O(1) 
   \end{align*}
   and as a consequence, for a.e $t\in(0,T)$, we obtain
   \begin{align}\label{average zero_estiamte}
      \nonumber \|\mathbb{P}\left[\partial_\nu u(\cdot,t)\right]\|_{L_0^2(\partial D)} & = \sup_{\|\phi\|_{L_0^2(\partial D)} \le 1} \left| \int_{\partial D} \partial_\nu u(x,t)\, \phi(x) \, dS_x\right| \\
      \nonumber & = \sup_{\|\phi\|_{L_0^2(\partial D)} \le 1} \left|\int_{\partial D} F_*(x,t) \, \left[\left(\alpha^{-1}(z)+ \frac{\rho_0(z)}{2}\right) \mathbb{I}\textnormal{d} + \mathcal{K}_{z} \right]^{-1}\phi (x) \right|\, dS_x\\
       & \preceq \|F_*(\cdot,t)\|_{L^2(\partial D)} .
   \end{align}
   In order to complete the proof, we need to calculate $L^2$ norm of the terms in right hand side of \eqref{boundary operator eqn}. Having the relation \eqref{regularity of boundary} and the time-regularity of $\partial_\nu u$ in hand, we derive
   \begin{align}\label{average zero_estiamte 1}
       \nonumber \left\|\int_{\partial D} \frac{(\cdot-y)\cdot\nu(\cdot)}{|\cdot-y|^3} \, \left[\partial_\nu u(y, t - c_0^{-1}(z)|\cdot-y|)-\partial_\nu u(y,t)\right] \, dS_y \right\|_{L^2(\partial D)} & \preceq \, \sup_{0\le\tau\le t} \int_{\partial D} |\partial_t\partial_\nu u(y,\tau)|\, dS_y \\
      \nonumber & \preceq |\partial D|^{1/2} \, \sup_{0\le\tau\le t} \|\partial_t\partial_\nu u(\cdot,\tau)\|_{L^2(\partial D)} \\
      & = \epsilon^3.
   \end{align}
   From \eqref{F_estimate} and the continuity of $\alpha$ in $\overline D$, we note that 
   \begin{align}\label{average zero_estiamte 2}
       \alpha^{-1}(z) \left\| F(\cdot,t)\right\|_{L^2(\partial D)} \preceq \epsilon^3 , \quad \alpha^{-1}(z)\|\left(\rho_0(\cdot)\alpha(\cdot)- \rho_0(z)\alpha(z)\right) \partial_\nu u(\cdot,t)\|_{L^2(\partial D)}\preceq \epsilon^3, \quad \textnormal{ for } t\in (0,T).
   \end{align}
   Therefore, it remains to estimate the term 
   $$\alpha^{-1}(z)\, \|\mathcal{K}^t_G \left[(\alpha(\cdot)-\alpha(z))\partial_\nu u\right](\cdot,t)\|_{L^2(\partial D)}$$
   for a.e. $t\in(0,T)$. For that purpose, we recall the definition of Neumann-Poincar\'e operator in \eqref{Neumann-poincare} and identify its most singular part in the expression \eqref{NP_expression} which we represent as the following operator
   \begin{align*}
       \mathcal{Q}[f](x,t) = \int_{\partial D} \frac{m(x,y) \partial_{\nu_x} \zeta(x,y)}{4\pi\zeta^2(x,y)}f(y,t-\zeta(x,y))\, dS_y. 
   \end{align*}
   To estimate the function $\mathcal{K}^t_{G}[f]$ in $L^2(\partial D)$ in space and pointwise in time, we first deal with the estimation $\mathcal{Q}[f]$. In this regard, we observe that  
   \begin{align*}
       \frac{\partial_{\nu_x} \zeta(x,y)}{\zeta^2(x,y)} = c_0^{-1}(y) \,\frac{\nu_x\cdot(x-y) + O(|x-y|^2)}{|\eta(x,y)|\,\zeta^2(x,y)} =  c_0^{-1}(y) \,\frac{\nu_x\cdot(x-y) }{|\eta(x,y)|\,\zeta^2(x,y)} + O(|x-y|^{-1}), \quad \textnormal{ for } x \neq y,
   \end{align*}
   which follows from the relations \eqref{reln_travel time} and \eqref{g_sum}. Now an application of the estimate \eqref{regularity of boundary} further leads us to the observation that
   \begin{align*}
       \frac{\partial_{\nu_x} \zeta(x,y)}{\zeta^2(x,y)} \preceq \epsilon^{-1} \frac{1}{|x-y|}
   \end{align*}
   which immediately implies
   \begin{align*}
       \left\vert \mathcal{Q}[f](x,t)\right\vert \preceq \epsilon^{-1} \sup_{0\le\tau\le t}\int_{\partial D} \frac{|f(y,\tau)|}{|x-y|}\, dS_y.
   \end{align*}
   As before, we appeal to the generalized Young's inequality \cite[Lemma 0.10]{Folland_PDE_book} to derive
   \begin{align*}
       \|\mathcal{Q}[f](\cdot,t)\|_{L^2(\partial D)} \preceq \left(\epsilon^{-1}\sup_{x\in\partial D} \int_{\partial D} \frac{dS_y}{|x-y|}\right) \sup_{0\le\tau\le t} \|f(\cdot,t)\|_{L^2(\partial D)} \preceq \sup_{0\le\tau\le t} \|f(\cdot,t)\|_{L^2(\partial D)}.
   \end{align*}
   With similar considerations in mind, we can estimate other terms in  \eqref{Neumann-poincare} and finally arrive at
   \begin{align*}
       \|\mathcal K^t_G[f](\cdot,t)\|_{L^2(\partial D)} \preceq  \sup_{0\le\tau\le t} \left(\|f(\cdot,\tau)\|_{L^2(\partial D)} + \|\partial_t f(\cdot,\tau)\|_{L^2(\partial D)}\right),
   \end{align*}
   which, in view of Lemma \ref{gradient estimate}, implies the bound
   \begin{align}\label{average zero_estiamte 3}
     \|  \alpha^{-1}(z) \,\mathcal{K}^t_G \left[(\alpha(\cdot)-\alpha(z))\,\partial_\nu u\right](\cdot,t) \|_{L^2(\partial D)} \preceq \epsilon^{2+1+2} = \epsilon^5.
   \end{align}
   Here, we have also used the smoothness of the function $\alpha(\cdot)$ in $\overline D$. With regard to the bounds \eqref{average zero_estiamte 1}, \eqref{average zero_estiamte 2} and \eqref{average zero_estiamte 3}, the estimate \eqref{average zero_estiamte} finally implies
   \begin{align*}
       \|\mathbb{P}\left[\partial_\nu u(\cdot,t)\right]\|_{L_0^2(\partial D)} \preceq \epsilon^3, \quad \textnormal{ for a.e. } t \in (0,T).
   \end{align*}
    To derive similar estimate for $\partial_t^2 \partial_\nu u$, we need the assumption $p \ge 7$.
\bigskip

{\bf Acknowledgement:}  The authors would like to thank the referees for their critics and suggestions that improved much the work.

\section{Funding and/or Conflicts of interests/Competing interests}
The research leading to these results received funding from the Austrian Science Fund FWF under Grant Agreement No 32660.
The authors have no relevant financial or non-financial interests to disclose.

\end{document}